\documentclass[a4paper, 11pt]{article}
\pdfoutput=1 

\usepackage[T1]{fontenc} 
\usepackage[utf8]{inputenc} 

\usepackage{amsfonts,amsmath,amssymb,mathtools,dsfont} 
\usepackage{microtype} 
\usepackage[usenames, dvipsnames]{xcolor} 

\usepackage{fancyhdr}

\usepackage[labelsep = period, labelfont = bf, justification = centering]{caption}
\usepackage{float,graphicx,subcaption} 

\usepackage{booktabs,makecell} 
\usepackage{enumitem,mdwlist} 
\setlist[itemize]{topsep=0ex,itemsep=0ex,parsep=0.4ex}
\setlist[enumerate]{topsep=0ex,itemsep=0ex,parsep=0.4ex}

\usepackage{etoolbox}
\usepackage[super]{nth} 
\allowdisplaybreaks

\usepackage[left = 2cm, right = 2cm, top = 2.5cm, bottom = 2.5cm, headsep = 12pt, headheight = 15pt]{geometry}
\usepackage{parskip}

\usepackage[noBBpl]{mathpazo} 
\DeclareFontFamily{U}{matha}{\hyphenchar\font45}
\DeclareFontShape{U}{matha}{m}{n}{
	<5> <6> <7> <8> <9> <10> gen * matha
	<10.95> matha10 <12> <14.4> <17.28> <20.74> <24.88> matha12
}{}
\DeclareSymbolFont{matha}{U}{matha}{m}{n}

\DeclareMathSymbol{\specialuparrow}{\mathrel}{matha}{"D2}
\DeclareMathSymbol{\specialrightarrow}{\mathrel}{matha}{"D1}

\DeclareFontFamily{U} {cmmi}{}

\DeclareFontShape{U}{cmmi}{m}{n}{
	<-6> cmmi5
	<6-7> cmmi6
	<7-8> cmmi7
	<8-9> cmmi8
	<9-10> cmmi9
	<10-12> cmmi10
	<12-> cmmi12}{}

\DeclareSymbolFont{Xcmmi} {U} {cmmi}{m}{n}

\DeclareMathSymbol{\mu}{\mathord}{Xcmmi}{'026}
\DeclareMathSymbol{\rho}{\mathord}{Xcmmi}{'032}
\DeclareMathSymbol{\varphi}{\mathord}{Xcmmi}{'047}

\DeclareFontFamily{U} {cmr}{}

\DeclareFontShape{U}{cmr}{m}{n}{
	<-6> cmr5
	<6-7> cmr6
	<7-8> cmr7
	<8-9> cmr8
	<9-10> cmr9
	<10-12> cmr10
	<12-> cmr12}{}

\DeclareSymbolFont{Xcmr} {U} {cmr}{m}{n}
\DeclareMathSymbol{\Delta}{\mathord}{Xcmr}{'001}
\DeclareMathSymbol{\Upsilon}{\mathord}{Xcmr}{'007}
\DeclareMathSymbol{\Omega}{\mathord}{Xcmr}{'012}

\usepackage[hyphens]{url} 
\usepackage[linktoc = all, hidelinks, colorlinks, unicode=true, pagebackref=true]{hyperref} 
\usepackage[capitalise, compress, nameinlink, noabbrev]{cleveref} 

\hypersetup{linkcolor={blue!70!black}, citecolor={green!70!black}, urlcolor={blue!70!black}}

\crefname{section}{\S}{\S\S} 
\Crefname{section}{Section}{Sections} 
\crefname{subsection}{\S}{\S\S} 
\Crefname{subsection}{Subsection}{Subsections} 
\crefname{page}{page}{pages}
\Crefname{page}{Page}{Pages}

\usepackage{pgf,tikz}
\usepackage{mathrsfs}
\usetikzlibrary{decorations.pathreplacing}
\usetikzlibrary[patterns]
\usetikzlibrary{arrows.meta}
\usepackage{tkz-euclide}

\tikzset{font={\fontsize{10pt}{12}\selectfont}}
\tkzSetUpPoint[fill=black, size = 3pt]
\tkzSetUpLine[color=black, line width=0.6pt]

\usepackage{amsthm,thmtools,thm-restate}

\declaretheoremstyle[
spaceabove = .5\baselineskip\@plus.2\baselineskip\@minus.2\baselineskip, 
spacebelow = .2\baselineskip\@plus.2\baselineskip\@minus.2\baselineskip,
headfont = \normalfont\bfseries,
notefont = \mdseries, 
notebraces = {(}{)},
bodyfont = \normalfont\itshape,
postheadspace = .5em,
headpunct = .
]{bolditalic}

\declaretheoremstyle[
spaceabove = .5\baselineskip\@plus.2\baselineskip\@minus.2\baselineskip, 
spacebelow = .2\baselineskip\@plus.2\baselineskip\@minus.2\baselineskip,
headfont = \normalfont\bfseries,
notefont = \mdseries, 
notebraces = {(}{)},
bodyfont = \normalfont,
postheadspace = .5em,
headpunct = .
]{boldnormal}

\declaretheoremstyle[
spaceabove = .2\baselineskip\@plus.2\baselineskip\@minus.2\baselineskip, 
spacebelow = .5\baselineskip\@plus.2\baselineskip\@minus.2\baselineskip,
headfont = \normalfont\itshape,
notefont = \mdseries, 
notebraces = {}{},
bodyfont = \normalfont,
postheadspace = .5em,
headpunct = .,
qed = \qedsymbol
]{proofstyle}

\declaretheoremstyle[
spaceabove = .5\baselineskip\@plus.2\baselineskip\@minus.2\baselineskip, 
spacebelow = .2\baselineskip\@plus.2\baselineskip\@minus.2\baselineskip,
headfont = \normalfont\bfseries,
notefont = \mdseries, 
notebraces = {(}{)},
bodyfont = \normalfont,
postheadspace = .5em,
headpunct = .,
qed = \qedsymbol
]{solutionstyle}

\renewcommand*{\backref}[1]{}
\renewcommand*{\backrefalt}[4]{
	\ifcase #1 Not cited.%
	\or $\specialuparrow$#2%
	\else $\specialuparrow$#2%
	\fi%
}

\declaretheorem[name = Conjecture, numberwithin = section, style = bolditalic, refname = {Conjecture,Conjectures}, Refname = {Conjecture,Conjectures}]{conjecture}
\declaretheorem[name = Corollary, numberlike = conjecture, style = bolditalic, refname = {Corollary,Corollaries}, Refname = {Corollary,Corollaries}]{corollary}
\declaretheorem[name = Definition, numberlike = conjecture, style = bolditalic, refname = {Definition,Definitions}, Refname = {Definition,Definitions}]{definition}

\declaretheorem[name = Lemma, numberlike = conjecture, style = bolditalic, refname = {Lemma,Lemmas}, Refname = {Lemma,Lemmas}]{lemma}

\declaretheorem[name = Proposition, numberlike = conjecture, style = bolditalic, refname = {Proposition,Propositions}, Refname = {Proposition,Propositions}]{proposition}

\declaretheorem[name = Theorem, numberlike = conjecture, style = bolditalic, refname = {Theorem,Theorems}, Refname = {Theorem,Theorems}]{theorem}
\declaretheorem[name = Theorem, numbered = no, style = bolditalic, refname = {Theorem,Theorems}, Refname = {Theorem,Theorems}]{theorem*}
\declaretheorem[name = Question, numberlike = conjecture, style = bolditalic, refname = {Question,Questions}, Refname = {Question,Questions}]{question}
\declaretheorem[name = Claim, numberwithin = conjecture, style = bolditalic, refname = {Claim,Claims}, Refname = {Claim,Claims}]{claim}

\declaretheorem[name = Remark, numberlike = conjecture, style = boldnormal, refname = {Remark,Remarks}, Refname = {Remark,Remarks}]{remark}

\declaretheorem[name = Proof, numbered = no, style = proofstyle, refname = {Proof,Proofs}, Refname = {Proof,Proofs}]{Proof}

\renewcommand{\epsilon}{\varepsilon}

\renewcommand{\emptyset}{\varnothing}

\DeclarePairedDelimiter{\abs}{\lvert}{\rvert}

\DeclarePairedDelimiter{\set}{\{}{\}}

\newcommand*{\ZZN}{\mathbb{Z}_{\geqslant 0}}

\newcommand*{\cD}{\mathcal{D}}
\newcommand*{\cF}{\mathcal{F}}

\DeclareMathOperator{\KG}{KG}
\DeclareMathOperator{\SG}{SG}

\newcommand{\defn}[1]{\textcolor{Maroon}{\emph{#1}}}

\begin{document}

\author{Freddie Illingworth\footnotemark[1]}

\title{\bf The chromatic profile of locally colourable graphs}

\date{10 May 2022}

\maketitle

\begin{abstract}
	The classical Andr\'{a}sfai-Erd\H{o}s-S\'{o}s theorem considers the chromatic number of $K_{r + 1}$-free graphs with large minimum degree, and in the case $r = 2$ says that any $n$-vertex triangle-free graph with minimum degree greater than $2/5 \cdot n$ is bipartite. This began the study of the chromatic profile of triangle-free graphs: for each $k$, what minimum degree guarantees that a triangle-free graph is $k$-colourable? The chromatic profile has been extensively studied and was finally determined by Brandt and Thomass\'{e}.
	
	Triangle-free graphs are exactly those in which each neighbourhood is one-colourable. As a natural variant, \L uczak and Thomass\'{e} introduced the notion of a locally bipartite graph in which each neighbourhood is 2-colourable. Here we study the chromatic profile of the family of graphs in which every neighbourhood is $b$-colourable (locally $b$-partite graphs) as well as the family where the common neighbourhood of every $a$-clique is $b$-colourable. Our results include the chromatic thresholds of these families as well as showing that every $n$-vertex locally $b$-partite graph with minimum degree greater than $(1 - 1/(b + 1/7)) \cdot n$ is $(b + 1)$-colourable. Understanding these locally colourable graphs is crucial for extending the Andr\'{a}sfai-Erd\H{o}s-S\'{o}s theorem to non-complete graphs, which we develop elsewhere.
\end{abstract}

\renewcommand{\thefootnote}{\fnsymbol{footnote}} 


\footnotetext[1]{Mathematical Institute, University of Oxford (\textsf{illingworth@maths.ox.ac.uk}). Research carried out while at DPMMS, University of Cambridge. Research supported by EPSRC grant 2114463.}

\renewcommand{\thefootnote}{\arabic{footnote}} 

\section{Introduction}\label{sec:intro}

In 1973, Erd\H{o}s and Simonovits~\cite{ErdosSimonovits1973} asked the following question: for each graph $H$ and positive integer~$k$, what $\delta$ guarantees that every $n$-vertex $H$-free graph with minimum degree greater than $\delta n$ is $k$-colourable? The values of $\delta$, as $k$ varies, form the chromatic profile of $H$-free graphs. More generally, for a family of graphs $\cF$, the \defn{chromatic profile} of $\cF$ is the sequence of values $\delta_{\chi}(\cF, k)$ where
\begin{equation*}
	\delta_{\chi}(\cF, k) = \inf\set{d \colon \text{if }\delta(G) \geqslant d \abs{G} \text{ and } G \in \cF, \text{then } \chi(G) \leqslant k}.
\end{equation*}
In the case where $\cF$ is the family of $H$-free graphs we write $\delta_{\chi}(H, k)$ for $\delta_{\chi}(\cF, k)$. The question of determining $\delta_{\chi}(H, k)$, first asked in~\cite{ErdosSimonovits1973}, was re-emphasised by Allen, B\"{o}ttcher, Griffiths, Kohayakawa, and Morris~\cite{ABGKM2013}. For general $H$ very little is known about the chromatic profile and indeed Erd\H{o}s and Simonovits described it as ``too complicated''.

There has been much greater success with the chromatic threshold. The \defn{chromatic threshold} of a family $\cF$ is the limit (and so infimum) of the sequence $\delta_{\chi}(\cF, k)$: that is,
\begin{align*}
	\delta_{\chi}(\cF) & = \inf_{k} \delta_{\chi}(\cF, k) \\
	& = \inf \set{d \colon \exists C = C(\cF, d) \text{ such that if } \delta(G) \geqslant d \abs{G} \text{ and } G \in \cF, \text{then } \chi(G) \leqslant C}.
\end{align*}
Allen, B\"{o}ttcher, Griffiths, Kohayakawa, and Morris~\cite{ABGKM2013} determined the chromatic threshold of $H$-free graphs for every graph $H$. They further obtained the chromatic threshold of \defn{locally bipartite graphs} -- the family of graphs in which each neighbourhood is bipartite -- confirming the conjecture of \L uczak and Thomass\'{e}~\cite{LuczakThomasse2010} that this threshold is $1/2$.

Understanding the chromatic profile of locally bipartite graphs is a very natural local to global colouring problem: every neighbourhood is 2-colourable and large, so it seems possible that the whole graph might have small chromatic number. A further motive (that originally brought them to our attention) is to extend the Andr\'{a}sfai-Erd\H{o}s-S\'{o}s theorem~\cite{AES1974} -- that theorem gives the first interesting value in the chromatic profile of complete graphs, namely
\begin{equation*}
	\delta_{\chi}(K_{r + 1}, r) = 1 - \frac{1}{r - 1/3}.
\end{equation*}
This theorem can be seen as a minimum degree analogue of Erd\H{o}s and Simonovits's stability theorem~\cite{Erdos1967,Erdos1968,Simonovits1968} for the structure of $K_{r + 1}$-free graphs with close to the maximum number of edges. Erd\H{o}s and Simonovits's result says that any $H$-free graph with $(1 - 1/r - o(1)) \binom{n}{2}$ edges (where $r + 1$ is the chromatic number of $H$) can be made $r$-partite by deleting $o(n^{2})$ edges. The natural minimum degree analogue asks what minimum degree guarantees that an $H$-free graph can be made $r$-partite by deleting $o(n^{2})$ edges. The theorem of Andr\'{a}sfai, Erd\H{o}s, and S\'{o}s answers this for cliques. It turns out that, in order to extend the theorem to non-complete $H$ (as in~\cite{Illingworth2023mindegstab}), it is crucial to better understand locally colourable graphs (as defined next), and that is our purpose in this paper.
\begin{definition}
	A graph is \defn{$a$-locally $b$-partite} if the common neighbourhood of every $a$-clique is $b$-colourable. A graph is \defn{locally $b$-partite} if it is 1-locally $b$-partite: the neighbourhood of every vertex is $b$-colourable. We use \defn{$\cF_{a, b}$} to denote the family of $a$-locally $b$-partite graphs.
\end{definition}%
The family of locally bipartite graphs is $\cF_{1, 2}$. Note that, in general, 
\begin{equation*}
	\cF_{1, \ell} \subset \cF_{2, \ell - 1} \subset \dotsb \subset \cF_{\ell, 1} = \set{G \colon G \text{ is } K_{\ell + 2}\text{-free}}.
\end{equation*}
In particular, $\cF_{a, b}$ is a subfamily of $K_{\ell + 2}$-free graphs where $\ell = a + b - 1$.

As mentioned, the Andr\'{a}sfai-Erd\H{o}s-S\'{o}s theorem gives the first interesting value in the chromatic profile of $K_{\ell + 2}$-free graphs, that is, of the family $\cF_{\ell, 1}$. For the basic case of triangle-free graphs, this shows $\delta_{\chi}(K_{3}, 2) = 2/5$. Hajnal (see~\cite{ErdosSimonovits1973}) showed that $\delta_{\chi}(K_{3}) \geqslant 1/3$. H\"{a}ggkvist~\cite{Haggkvist1982} showed that $\delta_{\chi}(K_{3}, 3) \geqslant 10/29$ and Jin~\cite{Jin1995} proved that equality holds. Finally, building on work of Chen, Jin, and Koh~\cite{CJK1997}, Brandt and Thomass\'{e}~\cite{BrandtThomasse2005} fully described triangle-free graphs with minimum degree above the chromatic threshold (showing that, in fact, $\delta_{\chi}(K_{3}) = \delta_{\chi}(K_{3}, 4) = 1/3$). Moreover, Goddard and Lyle~\cite{GoddardLyle2010}, and Nikiforov~\cite{Nikiforov2010} independently showed that the chromatic profile of $K_{r + 2}$-free graphs can be derived straightforwardly from that of triangle-free graphs. The chromatic profile of the $K_{\ell + 2}$-free graphs, that is, of $\cF_{\ell, 1}$, has thus been completely determined.

The chromatic profile of the family $\cF_{a, b}$ for $b \geqslant 2$ displays substantially different behaviour to the case $b = 1$, however. In \cref{sec:chromthresholds} we determine the chromatic threshold of the family of $a$-locally $b$-partite graphs. This result in the special case $a = 1$ and $b = 2$ was conjectured (with a construction for the lower bound) by \L uczak and Thomass\'{e}~\cite{LuczakThomasse2010} and proved by Allen, B\"{o}ttcher, Griffiths, Kohayakawa, and Morris~\cite{ABGKM2013}.

\begin{theorem}\label{chromthreFab}
	Let $a$ and $b$ be positive integers with $b \geqslant 2$ and let $\ell = a + b - 1$. Then
	\begin{equation*}
		\delta_{\chi}(\cF_{a, b}) = 1 - \frac{1}{\ell},
	\end{equation*}
	and, in particular, the chromatic threshold of locally $b$-partite graphs, $\delta_{\chi}(\cF_{1, b})$, is $1 - 1/b$.
\end{theorem}

For comparison, the chromatic threshold of the family $\cF_{\ell, 1}$ of $K_{\ell + 2}$-free graphs is $1 - 1/(\ell + 1/2)$.

Now for the chromatic profile. All $a$-locally $b$-partite graphs are $K_{a + b + 1}$-free and furthermore, the $n$-vertex $K_{a + b + 1}$-free graph with highest minimum degree (and most edges), the Tur\'{a}n graph~\cite{Turan1941}, $T_{a + b}(n)$ is $a$-locally $b$-partite and $(a + b)$-chromatic. In particular,
\begin{equation*}
	\delta_{\chi}(\cF_{a, b}, k) = 1 - \frac{1}{a + b} \text{ for } k = 1, 2, \dotsc, \ell \coloneqq a + b - 1,
\end{equation*}
since there are no $K_{a + b + 1}$-free graphs (and so no $a$-locally $b$-partite) graphs with $\delta(G) > (1 - 1/(a + b)) \cdot \abs{G}$. In particular, the first interesting value in the chromatic profile of $a$-locally $b$-partite graphs is $\delta_{\chi}(\cF_{a, b}, \ell + 1)$. An upper bound for this is given by Theorems \ref{main4alocalbpart} ($b \geqslant 3$) and \ref{main4alocalbipart} ($b = 2$). In \cref{sec:localbpart} we prove \cref{main4alocalbpart} for locally $b$-partite graphs (that is, for $a = 1$) and, in \cref{sec:alocalbpart}, we extend it to all $a$ and prove \cref{main4alocalbipart}.
\begin{theorem}\label{main4alocalbpart}
	Let $a$ and $b$ be positive integers with $b \geqslant 3$ and let $\ell = a + b - 1$. Then
	\begin{equation*}
		\delta_{\chi}(\cF_{a, b}, \ell + 1) \leqslant 1 - \frac{1}{\ell + 1/7},
	\end{equation*}
	and, in particular, every locally $b$-partite graph $G$ with $\delta(G) > (1 - 1/(b + 1/7)) \cdot \abs{G}$ is $(b + 1)$-colourable.
\end{theorem}

\begin{theorem}\label{main4alocalbipart}
	Let $a$ be a positive integer and $\ell = a + 1$. Then
	\begin{equation*}
		\delta_{\chi}(\cF_{a, 2}, \ell + 1) = 1 - \frac{1}{\ell + 1/3}.
	\end{equation*}
\end{theorem}
Note that the chromatic threshold of the family of $K_{\ell + 2}$-free graphs is $1 - 1/(\ell + 1/2)$. In particular, all chromatic profile values for $K_{\ell + 2}$-free graphs are greater than the first interesting value for $a$-locally $b$-partite ones ($b \geqslant 2$).

To extend the chromatic profile of triangle-free graphs to $K_{r + 1}$-free graphs as mentioned above, Goddard and Lyle~\cite{GoddardLyle2010} and Nikiforov~\cite{Nikiforov2010} showed that every $n$-vertex maximal $K_{r + 1}$-free graph with minimum degree greater than $\delta_{\chi}(K_{r + 1}) \cdot n$ consists of an independent set joined to a $K_{r}$-free graph. That is, maximal graphs of $\cF_{r, 1}$ with sufficiently large minimum degree are obtained from those in $\cF_{r - 1, 1}$ by joining an independent set. A simple induction then converts the structure of triangle-free graphs to the structure of $K_{r + 1}$-free graphs. It is natural to ask whether something similar can be done to convert between different families $\cF_{a, b}$. Firstly, there does not seem to be an easy way to convert between $\cF_{a, b - 1}$ and $\cF_{a, b}$ (certainly joining on an independent set fails). Although we obtain the upper bound for $\delta_{\chi}(\cF_{1, b}, b + 1)$ in \cref{main4alocalbpart} using knowledge of locally bipartite graphs, it is not a straightforward induction. 

Joining on an independent set to a graph in $\cF_{a - 1, b}$ does give a graph in $\cF_{a, b}$, but it is not clear that all maximal graphs in $\cF_{a, b}$ of large minimum degree are obtained in this way -- the lower value of $\delta_{\chi}(\cF_{a, b})$ for $b \geqslant 2$ means the structural lemma of Goddard, Lyle, and Nikiforov does not apply. While our arguments extending results from locally $b$-partite graphs to $a$-locally $b$-partite graphs are simpler, they interestingly do require knowledge of locally $b'$-partite graphs for all $b' \leqslant a + b$. It seems that the crux of understanding locally colourable graphs is understanding locally $b$-partite graphs, a sentiment we will crystallise in \cref{sec:alocalbpart}.

Some of the required knowledge of locally bipartite graphs was obtained by the author in~\cite{Illingworth2022localbipart}, which included results such as
\begin{equation*}
	\delta_{\chi}(\cF_{1, 2}, 3) = 4/7, \quad \text{and} \quad \delta_{\chi}(\cF_{1, 2}, 4) \leqslant 6/11.
\end{equation*}
The results in~\cite{Illingworth2022localbipart} give detailed information about $n$-vertex locally bipartite graphs with minimum degree greater than $6/11 \cdot n$, including homomorphism properties. However, in order to prove the upper bound for $\delta_{\chi}(\cF_{1, b}, b + 1)$ in \cref{main4alocalbpart}, we need to extend our structural knowledge of locally bipartite graphs down to $8/15$, though, fortunately, we do not need any homomorphism properties. We do this in \cref{sec:localbipart} and our full knowledge of locally bipartite graphs is summarised there in \cref{spec4localbip}.

\subsection{Notation}\label{sec:notation}

Let $G$ be a graph and $X \subset V(G)$. We write \defn{$\Gamma(X)$} for $\cap_{v \in X} \Gamma(v)$ (the common neighbourhood of the vertices of $X$) and \defn{$d(X)$} for $\abs{\Gamma(X)}$. We often omit the parentheses so $\Gamma(u, v) = \Gamma(u) \cap \Gamma(v)$ and $d(u, v) = \abs{\Gamma(u, v)}$. We write \defn{$G_{X}$} for $G[\Gamma(X)]$ so, for example, $G_{u, v}$ is the induced graph on the common neighbourhood of vertices $u$ and $v$. Note that $G$ being $a$-locally $b$-partite is equivalent to $\chi(G_{K}) \leqslant b$ for every $a$-clique $K$ of $G$. We make frequent use of the fact that for two vertices $u$ and $v$ of $G$
\begin{equation*}
	d(u, v) = d(u) + d(v) - \abs{\Gamma(u) \cup \Gamma(v)} \geqslant d(u) + d(v) - \abs{G} \geqslant 2 \delta(G) - \abs{G}.
\end{equation*}
Given a set of vertices $X \subset V(G)$, we write \defn{$e(X, G)$} for the number of ordered pairs of vertices $(x, v)$ with $x \in X$, $v \in G$, and $xv$ an edge in $G$. In particular, $e(X, G)$ counts each edge in $G[X]$ twice and each edge from $X$ to $G - X$ once and satisfies
\begin{equation*}
	e(X, G) = \sum_{x \in X} d(x) = \sum_{v \in G} \abs{\Gamma(v) \cap X}.
\end{equation*}
We generalise this notation to vertex weightings which will appear in many of our arguments. We will take a set of vertices $X \subset V(G)$ and assign weights $\omega \colon X \to \ZZN$ to the vertices of $X$. Then we define
\begin{equation*}
	\omega(X, G) = \sum_{x \in X} \omega(x) d(x) = \sum_{v \in G} \text{Total weight of the neighbours of } v \text{ in } X.
\end{equation*}
We will often use the word \defn{circuit} (as opposed to cycle) in our arguments. A circuit is a sequence of (not necessarily distinct) vertices $v_{1}, v_{2}, \dotsc, v_{\ell}$ with $\ell > 1$, $v_{i}$ adjacent to $v_{i + 1}$ (for $i = 1, 2, \dotsc, \ell - 1$) and $v_{\ell}$ adjacent to $v_{1}$. Note that in a locally bipartite graph the neighbourhood of any vertex does not contain an odd circuit (and of course does not contain an odd cycle). We use circuit to avoid considering whether some pairs of vertices are distinct when it is unnecessary to do so.

Given a graph $G$, a \defn{blow-up} of $G$ is a graph obtained by replacing each vertex $v$ of $G$ by a non-empty independent set $I_{v}$ and each edge $uv$ by a complete bipartite graph between classes $I_{u}$ and $I_{v}$. We say we have blown-up a vertex $v$ by $n$ if $\abs{I_{v}} = n$. It is often helpful to think of this as weighting vertex $v$ by $n$.

A blow-up is \defn{balanced} if the independent sets $(I_{v})_{v \in G}$ are as equal in size as possible. We use \defn{$G(t)$} to denote the graph obtained by blowing-up each vertex of $G$ by $t$: $G(t)$ is the balanced blow-up of $G$ on $t \abs{G}$ vertices. Note, for example, that the balanced blow-ups of the clique, $K_{r}$, are exactly the Tur\'{a}n graphs, $T_{r}(n)$. We note in passing that a graph has the same chromatic and clique number as any of its blow-ups. Furthermore, if $H$ is a blow-up of $G$, then $G$ is $a$-locally $b$-partite if and only if $H$ is $a$-locally $b$-partite.

Given two graphs $G$ and $H$, the \defn{join} of $G$ and $H$, denoted \defn{$G + H$}, is the graph obtained by taking disjoint copies of $G$ and $H$ and joining each vertex of the copy of $G$ to each vertex of the copy of $H$. Note that the chromatic and clique numbers of $G + H$ are the sum of the chromatic and clique numbers of $G$ and $H$.

A graph $G$ is \defn{homomorphic} to a graph $H$, written \defn{$G \to H$}, if there is a map $\varphi \colon V(G) \to V(H)$ such that for any edge $uv$ of $G$, $\varphi(u)\varphi(v)$ is an edge of $H$. A graph $G$ is homomorphic to a graph $H$ if and only if $G$ is a subgraph of some blow-up of $H$. In particular, if $G \to H$, then $\chi(G) \leqslant \chi(H)$ and, moreover, if $H$ is $a$-locally $b$-partite, then $G$ is also.

\subsection{Importance to extending the \texorpdfstring{Andr\'{a}sfai-Erd\H{o}s-S\'{o}s}{Andrasfai-Erdos-Sos} theorem}

As previously mentioned, our initial motivation for studying locally colourable graphs came from trying to extend the Andr\'{a}sfai-Erd\H{o}s-S\'{o}s theorem~\cite{AES1974} to non-complete graphs and obtain a minimum degree analogue of Erd\H{o}s and Simonovits's classical stability theorem~\cite{Erdos1967,Erdos1968,Simonovits1968}. In particular, we wish to determine, for each $(r + 1)$-chromatic $H$, the value of
\begin{equation*}
	\begin{split}
		\delta_{H} = \inf \set{c \colon &\textnormal{if }\abs{G} = n, \, \delta(G) \geqslant c n, \textnormal{ and } G \textnormal{ is } H\textnormal{-free}, \\
		& \textnormal{then } G \textnormal{ can be made } r\textnormal{-partite by deleting } o(n^{2}) \textnormal{ edges}}.
	\end{split}
\end{equation*}
Here we sketch the link between this and locally colourable graphs deferring a thorough treatment to~\cite{Illingworth2023mindegstab}. When $r = 2$ the situation is particularly clean.
\begin{theorem}
	Let $H$ be a 3-chromatic graph. There is a smallest positive integer $g$ such that $H$ is not homomorphic to $C_{2g + 1}$. Then
	\begin{equation*}
		\delta_{H} = \frac{2}{2g + 1}.
	\end{equation*}
\end{theorem}
The graph that shows $\delta_{H} \geqslant 2/(2g + 1)$ is a balanced blow-up of the cycle $C_{2g + 1}$. The intuitive explanation for the theorem is that the main obstacle for being close to (that is, within $o(n^{2})$ edges of) bipartite is containing some blow-up of an odd cycle. That odd cycle must be consistent with being $H$-free (in particular, the blow-up of the odd cycle must be $H$-free) and hence it is the first odd cycle to which $H$ is not homomorphic that determines $\delta_{H}$.

Understanding locally colourable graphs becomes crucial when extending this theorem to $r \geqslant 3$. For concreteness consider $r = 3$: we are interested in which graphs' blow-ups are the main obstacles for being close to tripartite. Given the importance of odd cycles when $r = 2$ it seems natural that odd wheels (the join of a single vertex and an odd cycle) would be obstacles here and indeed they are. However, there are further obstacles that do not contain odd wheels and so are locally bipartite. These observations suggest we should pay attention to 4-chromatic locally bipartite graphs as these may be obstacles for being close to tripartite. Many of the graphs appearing in our theorems here also play a leading role for minimum degree stability.

\section{Chromatic thresholds}\label{sec:chromthresholds}

The $(r - 1)$-locally 1-partite graphs are exactly the $K_{r + 1}$-free graphs, and their chromatic threshold was determined by Goddard and Lyle~\cite{GoddardLyle2010}, and Nikiforov~\cite{Nikiforov2010}.
\begin{equation*}
	\delta_{\chi}(\cF_{r - 1, 1}) = \delta_{\chi}(K_{r + 1}) = 1 - \frac{1}{r - 1/2}.
\end{equation*}
In this section we prove \cref{chromthreFab}, showing that $\delta_{\chi}(\cF_{a, b}) = 1 - 1/(a + b - 1)$ for all $a \geqslant 1$ and $b \geqslant 2$. The upper bound follows from the work of Allen, B\"{o}ttcher, Griffiths, Kohayakawa, and Morris~\cite{ABGKM2013}.

\begin{Proof}[of upper bound in \cref{chromthreFab}]
	Fix $a$ and $b$ positive integers with $b \geqslant 2$ and let $\ell = a + b - 1$. Let $d > 1 - 1/\ell$ and let $G \in \cF_{a, b}$ with $\delta(G) \geqslant d \abs{G}$. Now $K_{\ell - 1} + C_{5}$ has an $a$-clique whose common neighbourhoods contains $K_{b - 2} + C_{5}$ so is not $b$-colourable. In particular, $G$ is $(K_{\ell - 1} + C_{5})$-free.
	
	Now, in the language of~\cite{ABGKM2013}, $K_{\ell - 1} + C_{5}$ is $(\ell + 2)$-near-acyclic, so the chromatic threshold of the family of $(K_{\ell - 1} + C_{5})$-free graphs is $1 - 1/\ell$. In particular, there is a constant $C$ depending only upon $d$ and $\ell$ (and not on $G$) such that $\chi(G) \leqslant C$. Hence,
	\begin{equation*}
		\delta_{\chi}(\cF_{a, b}) \leqslant 1 - \frac{1}{\ell} = 1 - \frac{1}{a + b - 1}. \qedhere
	\end{equation*}
\end{Proof}

For the lower bound it suffices to give examples of graphs $G \in \cF_{a, b}$ with $\delta(G) \geqslant (1 - 1/(a + b - 1) - o(1)) \cdot \abs{G}$ which have arbitrarily large chromatic number. Allen, B\"{o}ttcher, Griffiths, Kohayakawa, and Morris~\cite{ABGKM2013} used graphs of large girth and high chromatic number as well as Borsuk-Hajnal graphs for their lower bounds. However, these depend upon the forbidden subgraph $H$ and do not seem to be applicable here where the collection of forbidden subgraphs is infinite. \L uczak and Thomass\'{e}~\cite{LuczakThomasse2010} modified a Borsuk graph to give the lower bound of 1/2 for the chromatic threshold of locally bipartite graphs. We will give a somewhat simpler example which gives the lower bound for all $\cF_{a, b}$ ($a \geqslant 1, b \geqslant 2$).

Our example is based upon the classical Schrijver graph~\cite{Schrijver1978}. The \defn{Kneser graph}, $\KG(n, k)$, is the graph whose vertex set is all $k$-subsets of $\set{1, 2, \dotsc, n}$ with two vertices adjacent if the corresponding $k$-sets are disjoint. In 1955, Kneser~\cite{Kneser1955} conjectured that the chromatic number of $\KG(n, k)$ is $n - 2k + 2$. This conjecture remained open for two decades and was first proved by Lov\'{a}sz~\cite{Lovasz1978} using homotopy theory (see also B\'{a}r\'{a}ny~\cite{Barany1978} and Greene~\cite{Greene2002} for very short proofs).

The \defn{Schrijver graph}, $\SG(n, k)$, is the graph whose vertex set is all $k$-subsets of $\set{1, 2, \dotsc, n}$ which do not contain both $i$ and $i + 1$ for any $i = 1, 2, \dotsc, n - 1$ and do not contain both $n$ and 1. Again vertices are adjacent if the corresponding $k$-sets are disjoint. Put another way, $\SG(n, k)$ is the induced subgraph of $\KG(n, k)$ obtained by deleting all vertices whose corresponding sets are supersets of any of $\set{1, 2}$, $\set{2, 3}$, $\set{3, 4}$, \ldots, $\set{n - 1, n}$, $\set{n, 1}$. Schrijver~\cite{Schrijver1978} showed that $\SG(n, k)$ is vertex-critical with chromatic number $\chi(\SG(n, k)) = \chi(\KG(n, k)) = n - 2k + 2$.

\begin{Proof}[of lower bound in \cref{chromthreFab}]
	Fix $a$ and $b$ positive integers with $b \geqslant 2$ and let $\ell = a + b - 1$. Fix $k$ and let $n = 2k + f(k)$ where $f(k)$ is a non-negative integer less than $k$, and both $f(k) \to \infty$, $f(k)/k \to 0$ as $k \to \infty$. Since $k > n/3$, the graph $\SG(n, k)$ is triangle-free. We will eventually consider a blow-up of the graph shown in \cref{fig:chromthresh}.
	
	\begin{figure}[H]
		\centering
		\begin{tikzpicture}
			\foreach \i in {-1,1}{
				\foreach \j in {-1,1}{
					\tkzDefPoint(4*\i - 1.6,1.5*\j){v_{\i,\j,n,1}}
					\tkzDefPoint(4*\i - 0.8,1.5*\j){v_{\i,\j,1,2}}
					\tkzDefPoint(4*\i,1.5*\j){v_{\i,\j,2,3}}
					\tkzDefPoint(4*\i + 0.8,1.5*\j){blob_{\i,\j}}
					\tkzDefPoint(4*\i + 1.6,1.5*\j){v_{\i,\j,n-1,n}}
					\tkzDrawPoints(v_{\i,\j,1,2},v_{\i,\j,2,3},v_{\i,\j,n-1,n},v_{\i,\j,n,1})
					\tkzLabelPoint[below](v_{\i,\j,1,2}){$v_{1,2}$}
					\tkzLabelPoint[below](v_{\i,\j,2,3}){$v_{2,3}$}
					\tkzLabelPoint[below](v_{\i,\j,n-1,n}){$v_{n-1,n}$}
					\tkzLabelPoint[below](v_{\i,\j,n,1}){$v_{n,1}$}
					\tkzLabelPoint[above = -4pt](blob_{\i,\j}){$\dotsc$}
					\tkzDefPoint(4*\i - 2.2,1.5*\j - 0.5){A_{\i,\j}}
					\tkzDefPoint(4*\i + 2.2,1.5*\j - 0.5){B_{\i,\j}}
					\tkzDefPoint(4*\i + 2.2,1.5*\j + 0.5){C_{\i,\j}}
					\tkzDefPoint(4*\i - 2.2,1.5*\j + 0.5){D_{\i,\j}}
					\tkzDefMidPoint(A_{\i,\j},B_{\i,\j}) \tkzGetPoint{AB_{\i,\j}}
					\tkzDefMidPoint(B_{\i,\j},C_{\i,\j}) \tkzGetPoint{BC_{\i,\j}}
					\tkzDefMidPoint(C_{\i,\j},D_{\i,\j}) \tkzGetPoint{CD_{\i,\j}}
					\tkzDefMidPoint(D_{\i,\j},A_{\i,\j}) \tkzGetPoint{DA_{\i,\j}}
					\tkzDrawPolySeg(A_{\i,\j},B_{\i,\j},C_{\i,\j},D_{\i,\j},A_{\i,\j})
				}
			}	
			\tkzDefPoint(10,0){v}
			\tkzDrawPoint(v)
			\tkzLabelPoint[right](v){$v$}
			
			\node(SG) at (0,3.8){$\SG(n, k)$};
			
			\tkzDrawSegments(AB_{-1,1},CD_{-1,-1} AB_{1,1},CD_{1,-1} BC_{-1,1},DA_{1,1} BC_{-1,-1},DA_{1,-1} B_{-1,1},D_{1,-1} C_{-1,-1},A_{1,1})
			\tkzDrawSegments(v,B_{-1,1} v,B_{1,1} v,C_{-1,-1} v,C_{1,-1})
			\draw[dashed, line width = 0.6pt] (SG) -- (C_{-1,1});
			\draw[dashed, line width = 0.6pt] (SG) -- (D_{1,1});
			\draw[dashed, line width = 0.6pt] (SG) -- (C_{-1,-1});
			\draw[dashed, line width = 0.6pt] (SG) -- (D_{1,-1});
		\end{tikzpicture}
		\caption{In this diagram $\ell = 5$.}\label{fig:chromthresh}
	\end{figure}
	
	\begin{itemize}[noitemsep]
		\item Each rectangle is an independent set of $n$ vertices: $v_{1, 2}$, $v_{2, 3}$, \ldots, $v_{n - 1, n}$, $v_{n, 1}$ -- we always consider indices modulo $n$.
		\item There are $\ell - 1$ rectangles and the vertices in rectangles form a complete $(\ell - 1)$-partite graph with $n$ vertices in each part.
		\item The vertex $v$ is adjacent to all of the $v_{i, i + 1}$ but has no neighbours in the copy of $\SG(n, k)$ -- in particular, the rectangles together with $v$ form a complete $\ell$-partite graph.
		\item Finally $A \in \SG(n, k)$ is adjacent to $v_{i, i + 1}$ if either $i$ or $i + 1$ is in $A$ (note, by the definition of the Schrijver graph, that it is impossible for both $i$ and $i + 1$ to be in $A$).
	\end{itemize}
	
	We first check that the graph, $G$, shown in \cref{fig:chromthresh} is locally $\ell$-partite. Fix a vertex $u$ of $G$. If $u = v$, then $G_{u}$ consists of the $\ell - 1$ rectangles and so is $(\ell - 1)$-colourable. Next suppose that $u$ is a vertex in the copy of $\SG(n, k)$, so $v \not \in \Gamma(u)$. As $\SG(n, k)$ is triangle-free, $\Gamma(u)$ consists of an independent set in $\SG(n, k)$ together with some vertices from the $\ell - 1$ rectangles. In particular, $G_{u}$ is $\ell$-colourable. Finally suppose that $u$ lies in the $\ell - 1$ rectangles. By symmetry, we may take $u$ to be a $v_{1, 2}$. Then $\Gamma(u)$ consists of two independent sets in $\SG(n, k)$ (sets containing 1 and sets containing 2), $\ell - 2$ rectangles and $v$. There are no edges from $v$ to the copy of $\SG(n, k)$, so $G_{u}$ is $\ell$-colourable.
	
	Now we consider a suitable blow-up of $G$. Let $s$ be a multiple of $n$ that is much larger than $\abs{\SG(n, k)}$. We will blow-up each $v_{i, i + 1}$ by $s/n$ so that all rectangles contain $s$ vertices and we will blow up $v$ by $s$ also. Note that $v$ together with the rectangles form a complete $\ell$-partite graph with $s$ vertices in each part. Keep the copy of $\SG(n, k)$ as it is. Call the resulting graph $G'$.
	\begin{itemize}[noitemsep]
		\item $G'$ has $\ell s + \abs{\SG(n, k)}$ vertices.
		\item Any $A \in \SG(n, k)$ is adjacent to $2k/n$ proportion of the $v_{i, i + 1}$ so has degree at least $(\ell - 1) s (2k)/n = (\ell - 1)s/(1 + f(k)/(2k))$.
		\item Any other vertex has degree at least $(\ell - 1) s$.
		\item $\chi(G') \geqslant \chi(\SG(n, k)) = n - 2k + 2 = f(k) + 2$.
	\end{itemize}
	Given any $C, \varepsilon > 0$, we may choose $k$ large enough so that $f(k) \geqslant C$ and $(\ell - 1)/(1 + f(k)/(2k)) \geqslant \ell - 1 - \varepsilon$, and then choose $s$ large enough so that $\ell s + \abs{\SG(n, k)} \leqslant (\ell + \varepsilon) s$. The resulting graph $G'$ has chromatic number at least $f(k) + 2 \geqslant C$ and
	\begin{equation*}
		\frac{\delta(G')}{\abs{G'}} \geqslant \frac{\ell - 1 -\varepsilon}{\ell + \varepsilon} \geqslant 1 - \frac{1}{\ell} - \varepsilon.
	\end{equation*}
	Being $a$-locally $b$-partite is preserved when taking blow-ups. Now $G$ is locally $\ell$-partite and so $G'$ is too. As $\cF_{1, \ell} \subset \cF_{a, b}$, $G'$ is $a$-locally $b$-partite. Thus,
	\begin{equation*}
		\delta_{\chi}(\cF_{a, b}) \geqslant 1 - \frac{1}{\ell} - \varepsilon = 1 - \frac{1}{a + b - 1} - \varepsilon,
	\end{equation*}
	but $\varepsilon > 0$ was arbitrary and so we have the required result.
\end{Proof}

\section{Structure of locally bipartite graphs down to \texorpdfstring{$8/15$}{8/15}}\label{sec:localbipart}

Our understanding of the structure of locally bipartite graphs can be summarised as follows. The graphs $\overline{C}_{7}$, $H_{2}^{+}$, $H_{2}$, etc.\ can be seen in \cref{fig:graphs} where they are discussed more thoroughly.

\begin{theorem}[locally bipartite graphs]\label{spec4localbip}
	Let $G$ be a locally bipartite graph.
	\begin{enumerate}[noitemsep, label=\textit{\alph{*}}., ref=\textit{\alph{*}}]
		\item If $\delta(G) > 4/7 \cdot \abs{G}$, then $G$ is 3-colourable. \label{a}
		\item If $\delta(G) > 5/9 \cdot \abs{G}$, then $G$ is homomorphic to $\overline{C}_{7}$. Also, $G$ is either 3-colourable or contains $\overline{C}_{7}$. \label{b}
		\item There is an absolute constant $\varepsilon > 0$ such that if $\delta(G) > (5/9 - \varepsilon) \cdot \abs{G}$, then $G$ is homomorphic to either $\overline{C}_{7}$ or $H_{2}^{+}$. \label{c}
		\item If $\delta(G) > 6/11 \cdot \abs{G}$, then $G$ is 4-colourable. Also, $G$ is either 3-colourable or contains $\overline{C}_{7}$ or $H_{2}^{+}$. In the first two cases, $G$ is homomorphic to $\overline{C}_{7}$. \label{d}
		\item If $\delta(G) > 7/13 \cdot \abs{G}$, then $G$ is either 3-colourable or contains $H_{2}$. \label{e}
		\item If $\delta(G) > 8/15 \cdot \abs{G}$, then $G$ is either 3-colourable or contains $H_{2}$ or $T_{0}$. \label{f}
	\end{enumerate}
\end{theorem}

In~\cite{Illingworth2022localbipart} we considered the structure of locally bipartite graphs with minimum degree down to $6/11$ and, indeed, the first four parts of \cref{spec4localbip} are proved there. The purpose of this section is to prove the final two parts. In \cref{sec:tight} we motivate this proof and also show that many of the constants in the theorem are tight.

\Cref{fig:graphs} displays the graphs of \cref{spec4localbip} as well as some that appear in the proof.

\begin{figure}[ht]
	\centering
	\begin{subfigure}{.24\textwidth}
		\centering
		\begin{tikzpicture}
			\foreach \pt in {0,1,...,6} 
			{
				\tkzDefPoint(\pt*360/7 + 90:1){v_\pt}
			} 
			\tkzDrawPolySeg(v_0,v_1,v_2,v_3,v_4,v_5,v_6,v_0) 
			\tkzDrawPolySeg(v_5,v_0,v_2)
			\tkzDrawPolySeg(v_1,v_3)
			\tkzDrawPolySeg(v_4,v_6)
			\tkzDrawPoints(v_0,v_...,v_6)
		\end{tikzpicture}
		\subcaption*{$H_{0}$}
	\end{subfigure}
	\begin{subfigure}{.24\textwidth}
		\centering
		\begin{tikzpicture}
			\foreach \pt in {0,1,...,6} 
			{
				\tkzDefPoint(\pt*360/7 + 90:1){v_\pt}
			} 
			\tkzDrawPolySeg(v_0,v_1,v_2,v_3,v_4,v_5,v_6,v_0) 
			\tkzDrawPolySeg(v_3,v_5,v_0,v_2,v_4)
			\tkzDrawPolySeg(v_6,v_1)
			\tkzDrawPoints(v_0,v_...,v_6)
		\end{tikzpicture}
		\subcaption*{$H_{1}$}
	\end{subfigure}
	\begin{subfigure}{.24\textwidth}
		\centering
		\begin{tikzpicture}
			\foreach \pt in {0,1,...,6} 
			{
				\tkzDefPoint(\pt*360/7 + 90:1){v_\pt}
			} 
			\tkzDrawPolySeg(v_0,v_1,v_2,v_3,v_4,v_5,v_6,v_0)
			\tkzDrawPolySeg(v_1,v_3,v_5,v_0, v_2,v_4,v_6)
			\tkzDrawPoints(v_0,v_...,v_6)
		\end{tikzpicture}
		\subcaption*{$H_{2}$}
	\end{subfigure}
	\begin{subfigure}{.24\textwidth}
		\centering
		\begin{tikzpicture}
			\foreach \pt in {0,1,...,6} 
			{
				\tkzDefPoint(\pt*360/7 + 90:1){v_\pt}
			} 
			\tkzDrawPolySeg(v_0,v_1,v_2,v_3,v_4,v_5,v_6,v_0) 
			\tkzDrawPolySeg(v_0,v_2,v_4,v_6,v_1,v_3,v_5,v_0)
			\tkzDrawPoints(v_0,v_...,v_6)
		\end{tikzpicture}
		\subcaption*{$C_{7}^{2} = \overline{C}_{7}$}
	\end{subfigure}
	
	\bigskip
	
	\begin{subfigure}{.24\textwidth}
		\centering
		\begin{tikzpicture}
			\foreach \pt in {0,1,...,6} 
			{
				\tkzDefPoint(\pt*360/7 + 90:1){v_\pt}
			} 
			\tkzDefPoint(0,0){u}
			\tkzDrawPolySeg(v_0,v_1,v_2,v_3,v_4,v_5,v_6,v_0)
			\tkzDrawPolySeg(v_1,v_3,v_5,v_0, v_2,v_4,v_6)
			\tkzDrawSegments(u,v_0 u,v_2 u,v_5)
			\tkzDrawPoints(v_0,v_...,v_6)
			\tkzDrawPoint(u)
		\end{tikzpicture}
		\subcaption*{$H_{2}^{+}$}
	\end{subfigure}
	\begin{subfigure}{.24\textwidth}
		\centering
		\begin{tikzpicture}
			\foreach \pt in {0,1,...,6} 
			{
				\tkzDefPoint(\pt*360/7 + 90:1){v_\pt}
			}
			\tkzDefPoint(0,0.33){t}
			\tkzDefPoint(0.46,-0.4){u_1}
			\tkzDefPoint(-0.46,-0.4){u_6}
			\tkzDrawPolySeg(v_0,v_1,v_2,v_3,v_4,v_5,v_6,v_0) 
			\tkzDrawSegments(t,v_0 t,u_1 t,u_6)
			\foreach \pt in {0,2,3,4,5,6}
			{
				\tkzDrawSegment(u_1,v_\pt)
			}
			\foreach \pt in {0,1,2,3,4,5}
			{
				\tkzDrawSegment(u_6,v_\pt)
			}
			\tkzDrawPoints(v_0,v_...,v_6)
			\tkzDrawPoints(t,u_1,u_6)
		\end{tikzpicture}
		\subcaption*{$T_{0}$}
	\end{subfigure}
	\begin{subfigure}{.24\textwidth}
		\centering
		\begin{tikzpicture}
			\foreach \pt in {0,1,...,6} 
			{
				\tkzDefPoint(\pt*360/7 + 90:1){v_\pt}
			}
			\tkzDefPoint(-0.33,0){ul}
			\tkzDefPoint(0.33,0){ur}
			
			\tkzDrawPolySeg(v_0,v_1,v_2,v_3,v_4,v_5,v_6,v_0) 
			\tkzDrawPolySeg(v_3,v_5,v_0,v_2,v_4)
			\tkzDrawPolySeg(v_6,v_1)
			\tkzDrawSegments(ul,v_0 ul,v_2 ul,v_3 ur,v_0 ur,v_5 ur,v_3)
			\tkzDrawPoints(v_0,v_...,v_6)
			\tkzDrawPoints(ul,ur)
		\end{tikzpicture}
		\subcaption*{$H_{1}^{++}$}
	\end{subfigure}
	\begin{subfigure}{.24\textwidth}
		\centering
		\begin{tikzpicture}
			\foreach \pt in {0,1,...,6} 
			{
				\tkzDefPoint(\pt*360/7 + 90:1){v_\pt}
			} 
			\tkzDefPoint(0,0){u}
			\tkzDrawPolySeg(v_0,v_1,v_2,v_3,v_4,v_5,v_6,v_0) 
			\foreach \pt in {0,1,...,6} 
			{
				\tkzDrawSegment(u,v_\pt)
			} 
			\tkzDrawPoints(v_0,v_...,v_6)
			\tkzDrawPoint(u)
		\end{tikzpicture}
		\subcaption*{$W_{7}$}
	\end{subfigure}
	\caption{The graphs of \cref{spec4localbip}.}\label{fig:graphs}
\end{figure}

\begin{itemize}[noitemsep]
	\item All graphs shown are 4-chromatic and all bar $W_{7}$ are locally bipartite.
	\item The graph $H_{0}$ is isomorphic to the \defn{Moser Spindle} -- the smallest 4-chromatic unit distance graph. $H_{0}$ is also the smallest 4-chromatic locally bipartite graph and so it is natural that it should play such an integral part in many of our results. The graph $\overline{C}_{7}$ is the complement (and also the square) of the 7-cycle.
	\item Adding a single edge to $H_{0}$ while maintaining local bipartiteness can give rise to two non-isomorphic graphs, one of which is $H_{1}$. The other will appear fleetingly in the proof of \cref{claim:H02H1}. Adding a single edge to $H_{1}$ while maintaining local bipartiteness gives rise to a unique (up to isomorphism) graph -- $H_{2}$. There is only one way to add a single edge to $H_{2}$ and maintain local bipartiteness -- this gives $\overline{C}_{7}$. $H_{2}^{+}$ is $H_{2}$ with a degree 3 vertex added. $H_{1}^{++}$ is $H_{1}$ with two degree 3 vertices added.
	\item $\overline{C}_{7}$ and $H_{2}^{+}$ are both edge-maximal locally bipartite graphs.
	\item $T_{0}$ is a 7-cycle (the outer cycle) together with two vertices each joined to six of the seven vertices in the outer cycle (with the `seventh' vertices two apart) and finally a vertex of degree three is added.
	\item $W_{7}$, is called the \defn{7-wheel}. More generally, a single vertex joined to all the vertices of a $k$-cycle is called a \defn{$k$-wheel} and is denoted by $W_{k}$. We term any edge from the central vertex to the cycle a \defn{spoke} of the wheel and any edge of the cycle a \defn{rim} of the wheel. Note that a graph is locally bipartite exactly if it does not contain any odd wheel (there is no such nice characterisation for a graph being locally tripartite, locally 4-partite, \ldots).
\end{itemize}

The following observation gives a useful link between local bipartiteness and some of these graphs. We will use it frequently when copies of $H_{0}$, $H_{1}$, $H_{2}$ or $\overline{C}_{7}$ appear.

\begin{remark}\label{rmk:4nbs}
	Any five vertices of $H_{0}$ contain a triangle or a 5-cycle. In particular, if $G$ is a locally bipartite graph, then any vertex can have at most four neighbours in any copy of $H_{0}$ appearing in $G$.
\end{remark}

Containment and homomorphisms between the first seven graphs of \cref{fig:graphs} are summarised in \cref{fig:containment}. Note that a full arrow pointing from $H$ to $G$ signifies that $H$ is a subgraph of $G$ and a dashed arrow from $H$ to $G$ signifies that $H$ is homomorphic to $G$. Furthermore, for two graphs $H$ and $G$ in the diagram, $H$ is homomorphic to $G$ exactly if there is a sequence of arrows starting at $H$ and ending at $G$. We verify \cref{fig:containment} in the appendix.

\begin{figure}[H]
	\centering
	\begin{tikzpicture}
		\node(H21) at (0,0){$H_{2}^{+}$};
		\node(blank) at (-1.7,0){\color{white}H1};
		\node(C7) at (-3.4,0){$\overline{C}_{7}$};
		\node(H12) at (0,-1.7){$H_{1}^{++}$};
		\node(T0) at (1.7,-1.7){$T_{0}$};
		\node(H2) at (-1.7,-1.7){$H_{2}$};
		\node(H1) at (-1.7,-3.4){$H_{1}$};
		\node(H0) at (-1.7,-5.1){$H_{0}$};

		\draw[-Latex] (H2) -- (H21);
		\draw[-Latex] (H2) -- (C7);
		\draw[dashed, -Latex] (H12) -- (H21);
		\draw[dashed, -Latex] (T0) -- (H21);
		\draw[-Latex] (H1) -- (H2);
		\draw[-Latex] (H1) -- (H12);
		\draw[-Latex] (H0) -- (H1);
	\end{tikzpicture}
	\caption{Full arrows denote containment, dashed arrows denote homomorphisms.}\label{fig:containment}
\end{figure}

\subsection{The components of \texorpdfstring{\cref{spec4localbip}}{Theorem 3.1} and initial observations}\label{sec:tight}

We will start by noting which parts of \cref{spec4localbip} are tight. We will need suitable blow-ups of some of the graphs in \cref{fig:graphs}. Below we give these blow-ups which have been chosen so that the ratio between the minimum degree and order of the graph is as large as possible. When we \defn{weight a vertex by $0^{+}$} we are actually giving it some tiny positive weight (we have not deleted the vertices entirely just given them a small weight relative to the rest). 

\begin{figure}[H]
	\centering
	\begin{subfigure}{.3\textwidth}
		\centering
		\begin{tikzpicture}
			\foreach \pt in {0,1,...,6} 
			{
				\tkzDefPoint(\pt*360/7 + 90:1.5){v_\pt}
			} 
			\tkzDefPoint(0,0){u}
			\tkzDrawPolySeg(v_0,v_1,v_2,v_3,v_4,v_5,v_6,v_0)
			\tkzDrawPolySeg(v_1,v_3,v_5,v_0, v_2,v_4,v_6)
			\tkzDrawSegments(u,v_0 u,v_2 u,v_5)
			\tkzDrawPoints(v_0,v_...,v_6)
			\tkzDrawPoint(u)
			\tkzLabelPoint[above](v_0){2}
			\tkzLabelPoint[left](v_1){$0^{+}$}
			\tkzLabelPoint[left](v_2){2}
			\tkzLabelPoint[below](v_3){1}
			\tkzLabelPoint[below](v_4){1}
			\tkzLabelPoint[right](v_5){2}
			\tkzLabelPoint[right](v_6){$0^{+}$}
			\tkzLabelPoint[below](u){1}
		\end{tikzpicture}
		\subcaption*{$H_{2}^{+}$ weighted}
	\end{subfigure}
	\begin{subfigure}{.3\textwidth}
		\centering
		\begin{tikzpicture}
			\foreach \pt in {0,1,...,6} 
			{
				\tkzDefPoint(\pt*360/7 + 90:1.5){v_\pt}
			} 
			\tkzDrawPolySeg(v_0,v_1,v_2,v_3,v_4,v_5,v_6,v_0)
			\tkzDrawPolySeg(v_1,v_3,v_5,v_0, v_2,v_4,v_6)
			\tkzDrawPoints(v_0,v_...,v_6)
			\tkzLabelPoint[above](v_0){3}
			\tkzLabelPoint[left](v_1){1}
			\tkzLabelPoint[left](v_2){2}
			\tkzLabelPoint[below](v_3){1}
			\tkzLabelPoint[below](v_4){1}
			\tkzLabelPoint[right](v_5){2}
			\tkzLabelPoint[right](v_6){1}
		\end{tikzpicture}
		\subcaption*{$H_{2}$ weighted}
	\end{subfigure}
	\begin{subfigure}{.3\textwidth}
		\centering
		\begin{tikzpicture}
			\foreach \pt in {0,1,...,6} 
			{
				\tkzDefPoint(\pt*360/7 + 90:1.5){v_\pt}
			}
			\tkzDefPoint(0,0.5){t}
			\tkzDefPoint(0.7,-0.6){u_1}
			\tkzDefPoint(-0.7,-0.6){u_6}
			\tkzDrawPolySeg(v_0,v_1,v_2,v_3,v_4,v_5,v_6,v_0) 
			\tkzDrawSegments(t,v_0 t,u_1 t,u_6)
			\foreach \pt in {0,2,3,4,5,6}
			{
				\tkzDrawSegment(u_1,v_\pt)
			}
			\foreach \pt in {0,1,2,3,4,5}
			{
				\tkzDrawSegment(u_6,v_\pt)
			}
			\tkzDrawPoints(v_0,v_...,v_6)
			\tkzDrawPoints(t,u_1,u_6)
			\tkzLabelPoint[above](v_0){4}
			\tkzLabelPoint[left](v_1){$0^{+}$}
			\tkzLabelPoint[left](v_2){$0^{+}$}
			\tkzLabelPoint[below](v_3){1}
			\tkzLabelPoint[below](v_4){1}
			\tkzLabelPoint[right](v_5){$0^{+}$}
			\tkzLabelPoint[right](v_6){$0^{+}$}
			\tkzLabelPoint[below](t){1}
			\tkzLabelPoint[below right = -3pt](u_1){3}
			\tkzLabelPoint[below left = -3pt](u_6){3}
		\end{tikzpicture}
		\subcaption*{$T_{0}$ weighted}
	\end{subfigure}
	
	\bigskip
	
	\begin{subfigure}{.3\textwidth}
		\centering
		\begin{tikzpicture}
			\foreach \pt in {0,1,...,6} 
			{
				\tkzDefPoint(\pt*360/7 + 90:1.5){v_\pt}
			}
			\tkzDefPoint(-0.5,0){ul}
			\tkzDefPoint(0.5,0){ur}
			
			\tkzDrawPolySeg(v_0,v_1,v_2,v_3,v_4,v_5,v_6,v_0) 
			\tkzDrawPolySeg(v_3,v_5,v_0,v_2,v_4)
			\tkzDrawPolySeg(v_6,v_1)
			\tkzDrawSegments(ul,v_0 ul,v_2 ul,v_3 ur,v_0 ur,v_5 ur,v_3)
			\tkzDrawPoints(v_0,v_...,v_6)
			\tkzDrawPoints(ul,ur)
			\tkzLabelPoint[above](v_0){5}
			\tkzLabelPoint[left](v_1){$0^{+}$}
			\tkzLabelPoint[left](v_2){3}
			\tkzLabelPoint[below](v_3){2}
			\tkzLabelPoint[below](v_4){$0^{+}$}
			\tkzLabelPoint[right](v_5){3}
			\tkzLabelPoint[right](v_6){$0^{+}$}
			\tkzLabelPoint[right](ul){1}
			\tkzLabelPoint[left](ur){1}
		\end{tikzpicture}
		\subcaption*{$H_{1}^{++}$ weighted}
	\end{subfigure}
	\caption{Weightings of $H_2^{+}$, $H_2$, $T_0$, and $H_{1}^{++}$.}\label{fig:weightings}
\end{figure}

Now we address one-by-one the tightness of the constants in \cref{spec4localbip}. First note that $\overline{C}_{7}$ is locally bipartite, 4-regular with 7 vertices and chromatic number 4. Hence balanced blow-ups of $\overline{C}_{7}$ show that $4/7$ is tight in part~\ref{a}. \Cref{fig:weightings} shows that there are $n$-vertex blow-ups of $H_{2}^{+}$ with minimum degree at least $5/9 \cdot n - \mathcal{O}(1)$. These blow-ups of $H_{2}^{+}$ are 4-chromatic (as $H_{2}^{+}$ is), do not contain $\overline{C}_{7}$ nor are homomorphic to $\overline{C}_{7}$ (as neither of $\overline{C}_{7}$ nor $H_{2}^{+}$ is homomorphic to the other). This gives the tightness of $5/9$ in part~\ref{b}. \Cref{fig:weightings} shows that there are $n$-vertex blow-ups of $H_{2}$ with minimum degree $\lfloor 6/11 \cdot n \rfloor$. These blow-ups are 4-chromatic and contain neither $\overline{C}_{7}$ nor $H_{2}^{+}$, since neither of these is homomorphic to $H_{2}$. Hence, in part~\ref{d}, $6/11$ is tight for the conclusion that $G$ is either 3-colourable or contains $\overline{C}_{7}$ or $H_{2}^{+}$. Tightness is not known for 4-colourability and in fact it seems likely that $\delta_{\chi}(\cF_{1, 2}, 4) < 6/11$. Similar arguments show that the weightings of $T_{0}$ and $H_{1}^{++}$ in \cref{fig:weightings} give the tightness of $7/13$ and $8/15$ in parts~\ref{e} and \ref{f} respectively.

Our second aim in this section is to motivate the proof of parts~\ref{e} and \ref{f} of \cref{spec4localbip}. We will deduce them from the following two theorems.

\begin{restatable}{theorem}{mainforlocalbip}\label{main4localbip}
	Let $G$ be a locally bipartite graph. If $\delta(G) > 8/15 \cdot \abs{G}$, then $G$ is either 3-colourable or contains $H_{0}$ or $T_{0}$. If $\delta(G) > 7/13 \cdot \abs{G}$, then $G$ is either 3-colourable or contains $H_{0}$.
\end{restatable}

\vspace{-10pt}

\begin{theorem}\label{lemma4H}
	Let $G$ be a locally bipartite graph that contains $H_{0}$. If $\delta(G) > 8/15 \cdot \abs{G}$, then $G$ contains $H_{2}$.
\end{theorem}

\begin{Proof}[of parts \ref{e} and \ref{f}]
	Let $G$ be a locally bipartite graph. If $\delta(G) > 7/13 \cdot \abs{G}$, then, by \cref{main4localbip}, $G$ is either 3-colourable or contains $H_{0}$. If $G$ contains $H_{0}$, then \cref{lemma4H} shows it actually contains $H_{2}$. This gives part~\ref{e}. Suppose instead that $\delta(G) > 8/15 \cdot \abs{G}$. By \cref{main4localbip}, $G$ is either 3-colourable, contains $T_{0}$ or contains $H_{0}$. In the last case \cref{lemma4H} shows that $G$ contains $H_{2}$ and so we have part~\ref{f}.
\end{Proof}

We prove \cref{lemma4H} in \cref{sec:H02H2} by counting edges between the copy of $H_{0}$ and $G$. \Cref{main4localbip} is more involved. As motivation, consider a locally bipartite $H_{2}$-free $G$ with $\delta(G) > 7/13 \cdot \abs{G}$ (we will ignore the $8/15$ conclusion here). By \cref{lemma4H}, $G$ is $H_{0}$-free and our aim is to show that $G$ is \mbox{3-colourable}. We may as well assume that $G$ is edge-maximal: the addition of any edge will give a copy of $H_{0}$ or a vertex with non-bipartite neighbourhood. Thus, any non-edge of $G$ is either the missing edge of a $K_{4}$, a missing rim of an odd wheel, a missing spoke of an odd wheel or a missing edge of an $H_{0}$. This motivates a key definition that also appeared in~\cite{Illingworth2022localbipart}. Recall that $G_{u, v}$ denotes the subgraph of $G$ induced by the common neighbourhood of vertices $u$, $v$.

\begin{definition}[dense and sparse]\label{def:dense}
	A pair of \textbf{non-adjacent}, \textbf{distinct} vertices $u, v$ in a graph $G$ is \defn{dense} if $G_{u, v}$ contains an edge and \defn{sparse} if $G_{u, v}$ does not contain an edge.
\end{definition}

First note that every pair of distinct vertices in any graph is exactly one of `adjacent', `dense' or `sparse'. Another way to view being dense is as being the missing edge of a $K_{4}$. Locally bipartite graphs are $K_{4}$-free so any pair of distinct vertices with an edge in their common neighbourhood must be non-adjacent and so dense. Our initial observations above show that each sparse pair of vertices in $G$ is either the missing rim or spoke of an odd wheel or a missing edge of an $H_{0}$. In \cref{sec:oddspokes} we will rule out sparse pairs of vertices being the missing spoke of an odd wheel and in \cref{sec:proof4localbip} we will finish the proof by showing that $G$ is \mbox{3-colourable}.

The distinction between dense and sparse pairs turns out to be crucial. We borrow from~\cite{Illingworth2022localbipart} the following three simple but useful lemmas. The final one hints at the importance of $H_{0}$.

\begin{lemma}\label{lemma4I}
	Let $G$ be a graph with $\delta(G) > 1/2 \cdot \abs{G}$ and let $I$ be any largest independent set in $G$. Then, for every distinct $u, v \in I$, the pair $u, v$ is dense.
\end{lemma}

\begin{lemma}\label{lemma4sparse}
	Let $G$ be a graph with $\delta(G) > 1/2 \cdot \abs{G}$ and suppose $C$ is an induced 4-cycle in $G$. Then at least one of the non-edges of $C$ is a dense pair.
\end{lemma}

\begin{lemma}\label{lemma4dense}
	Let $G$ be a locally bipartite graph which does not contain $H_{0}$. For any vertex $v$ of $G$,
	\begin{equation*}
		D_{v} \coloneqq \set{u \colon \text{the pair } u, v \text{ is dense}}
	\end{equation*}
	is an independent set of vertices.
\end{lemma}

\subsection{From \texorpdfstring{$H_{0}$}{H0} to \texorpdfstring{$H_{2}$}{H2} -- proof of \texorpdfstring{\cref{lemma4H}}{Theorem 3.3}}\label{sec:H02H2}

In this subsection we prove \cref{lemma4H} in two steps. The strategy is to start with a copy of $H_{0}$ and consider edges between it and the rest of $G$. Using the high minimum degree we are able to find a vertex with the correct neighbours in the copy of $H_{0}$ so that a copy of $H_{1}$ is present. We then play the same game to get a copy of $H_{2}$. To aid the reader, we will use \defn{$G[v_{0}, v_{1}, \dotsc, v_{6}]$ is a copy of $H_{1}$} to mean that $G$ has a copy of $H_{1}$ in which $v_{0}$ is the top vertex (as displayed in \cref{fig:graphs}) and $v_{0}$, $v_{1}$, \ldots, $v_{6}$ precede anticlockwise round the figure. We do similarly for copies of $H_{2}$.

\begin{claim}\label{claim:H02H1}
	Let $G$ be a locally bipartite graph containing $H_{0}$. If $\delta(G) > 1/2 \cdot \abs{G}$, then $G$ contains $H_{1}$.
\end{claim}

\begin{Proof}
	We label a copy of $H_{0}$ in $G$ as below and let $X = \set{a_{0}, a_{1}, \dotsc, a_{6}}$.
	\begin{figure}[H]
		\centering
		\begin{tikzpicture}
			\foreach \pt in {0,1,...,6} 
			{
				\tkzDefPoint(\pt*360/7 + 90:1){v_\pt}
			} 
			\tkzDrawPolySeg(v_0,v_1,v_2,v_3,v_4,v_5,v_6,v_0) 
			\tkzDrawPolySeg(v_5,v_0,v_2)
			\tkzDrawPolySeg(v_1,v_3)
			\tkzDrawPolySeg(v_4,v_6)
			\tkzDrawPoints(v_0,v_...,v_6)
			\tkzLabelPoint[above](v_0){$a_{0}$}
			\tkzLabelPoint[left](v_1){$a_{1}$}
			\tkzLabelPoint[left](v_2){$a_{2}$}
			\tkzLabelPoint[below](v_3){$a_{3}$}
			\tkzLabelPoint[below](v_4){$a_{4}$}
			\tkzLabelPoint[right](v_5){$a_{5}$}
			\tkzLabelPoint[right](v_6){$a_{6}$}
		\end{tikzpicture}
	\end{figure}
	
	Let $U_{4}$ be the set of vertices with exactly four neighbours in $X$. \Cref{rmk:4nbs} shows that no vertex has five neighbours in a copy of $H_{0}$, so every non-$U_{4}$ vertex has at most three neighbours in $X$. Hence
	\begin{equation*}
		7/2 \cdot \abs{G} < 7 \delta (G) \leqslant e(X, G) \leqslant 4 \abs{U_{4}} + 3(\abs{G} - \abs{U_{4}}) = 3 \abs{G} + \abs{U_{4}},
	\end{equation*}
	and so 
	\begin{equation*}
		\abs{U_{4}} > 1/2 \cdot \abs{G}.
	\end{equation*}
	Now $\abs{U_{4}} + d(a_{0}) > \abs{G}$ and so some vertex $v$ is adjacent to $a_{0}$ and has exactly four neighbours in $X$. Note that $v$ cannot be adjacent to both $a_{1}$, $a_{2}$ as otherwise $va_{0}a_{1}a_{2}$ is a $K_{4}$, so by symmetry we may assume that $v$ is not adjacent to $a_{1}$. Similarly we may assume that $v$ is not adjacent to $a_{5}$. But $v$ has four neighbours in $X$ so must be adjacent to at least one of $a_{2}$, $a_{6}$ -- by symmetry we may assume $v$ is adjacent to $a_{2}$.
	
	There are two possibilities: $v$ is adjacent to $a_{0}, a_{2}, a_{3}, a_{4}$, or $v$ is adjacent to $a_{0}, a_{2}, a_{6}$ and one of $a_{3}$, $a_{4}$. In the latter case we may assume by symmetry that $v$ is adjacent to $a_{3}$. Hence there are two possibilities for $\Gamma(v) \cap X$: $\set{a_{0}, a_{2}, a_{3}, a_{4}}$ and $\set{a_{0}, a_{2}, a_{3}, a_{6}}$. In both cases $v$ cannot be any $a_{i}$ except for possibly $a_{1}$.
	
	If $\Gamma(v) \cap X = \set{a_{0}, a_{2}, a_{3}, a_{4}}$, then $G[v, a_{3}, a_{4}, a_{5}, a_{6}, a_{0}, a_{2}]$ is a copy of $H_{1}$. If $\Gamma(v) \cap X = \set{a_{0}, a_{2}, a_{3}, a_{6}}$, then $G$ contains the following graph where, in particular, $v$ is not adjacent to $a_{4}$.
	
	\begin{figure}[H]
		\centering
		\begin{tikzpicture}
			\foreach \pt in {0,1,...,6} 
			{
				\tkzDefPoint(\pt*360/7 + 90:1){v_\pt}
			} 
			\tkzDrawPolySeg(v_0,v_1,v_2,v_3,v_4,v_5,v_6,v_0) 
			\tkzDrawPolySeg(v_4,v_6,v_1,v_3)
			\tkzDrawPolySeg(v_5,v_0,v_2)
			\tkzDrawPoints(v_0,v_...,v_6)
			\tkzLabelPoint[above](v_0){$a_{0}$}
			\tkzLabelPoint[left](v_1){$v$}
			\tkzLabelPoint[left](v_2){$a_{2}$}
			\tkzLabelPoint[below](v_3){$a_{3}$}
			\tkzLabelPoint[below](v_4){$a_{4}$}
			\tkzLabelPoint[right](v_5){$a_{5}$}
			\tkzLabelPoint[right](v_6){$a_{6}$}
		\end{tikzpicture}
	\end{figure}
	
	Vertex $a_{6}$ is not adjacent to $a_{3}$ else $G_{a_{6}}$ contains the 5-cycle $a_{0}va_{3}a_{4}a_{5}$. Hence $va_{6}a_{4}a_{3}$ is an induced 4-cycle in $G$. By \cref{lemma4sparse}, at least one of the pairs $v, a_{4}$ and $a_{3}, a_{6}$ is dense. By symmetry we may assume that $v, a_{4}$ is dense: let $a'_{2}a'_{3}$ be an edge in $G_{v, a_{4}}$. Vertex $v$ is not adjacent to $a_{5}$, so $a_{5}$ is neither $a'_{2}$ nor $a'_{3}$. Note that $a_{0}$ is not adjacent to $a_{4}$ (else $a_{0}a_{4}a_{5}a_{6}$ is a $K_{4}$) and so $a_{0}$ is neither $a'_{2}$ nor $a'_{3}$. If $a_{6} = a'_{2}$, then $G_{a_{6}}$ contains the 5-cycle $a'_{3}va_{0}a_{5}a_{4}$, which is impossible. Similarly $a_{6} \neq a'_{3}$. Hence, $a'_{2}, a'_{3}$ are distinct from $a_{4}, a_{5}, a_{6}, a_{0}, v$ and so $G[a_{6}, a_{0}, v, a'_{2}, a'_{3}, a_{4}, a_{5}]$ is a copy of $H_{1}$.
\end{Proof}

\begin{claim}
	Let $G$ be a locally bipartite graph containing $H_{1}$. If $\delta(G) > 8/15 \cdot \abs{G}$, then $G$ contains $H_{2}$.
\end{claim}

\begin{Proof}
	Consider a copy of $H_{1}$ with vertices $X = \set{a_{0}, a_{1}, \dotsc, a_{6}}$. We assign a weighting $\omega \colon X \to \ZZN$ as shown in the diagram below, so, for example $\omega(a_{0}) = 2$ and $\omega(a_{1}) = 1$ (recall this notation from \cref{sec:notation}). We will often use diagrams to give weightings in this way. For each vertex $v \in G$, let $f(v)$ be the total weight of the neighbours of $v$ in $X$.
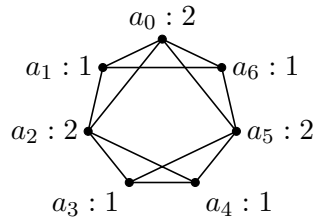
	\begin{figure}[H]
		\centering
		\begin{tikzpicture}
			\foreach \pt in {0,1,...,6} 
			{
				\tkzDefPoint(\pt*360/7 + 90:1){v_\pt}
			} 
			\tkzDrawPolySeg(v_0,v_1,v_2,v_3,v_4,v_5,v_6,v_0) 
			\tkzDrawPolySeg(v_3,v_5,v_0,v_2,v_4)
			\tkzDrawPolySeg(v_6,v_1)
			\tkzDrawPoints(v_0,v_...,v_6)
			\tkzLabelPoint[above](v_0){$a_{0} \colon 2$}
			\tkzLabelPoint[left](v_1){$a_{1} \colon 1$}
			\tkzLabelPoint[left](v_2){$a_{2} \colon 2$}
			\tkzLabelPoint[below left](v_3){$a_{3} \colon 1$}
			\tkzLabelPoint[below right](v_4){$a_{4} \colon 1$}
			\tkzLabelPoint[right](v_5){$a_{5} \colon 2$}
			\tkzLabelPoint[right](v_6){$a_{6} \colon 1$}
		\end{tikzpicture}
	\end{figure}
	
	We will assume that $G$ does not contain $H_{2}$. We first show that any vertex $v$ has $f(v) \leqslant 6$ and further that if $f(v) = 6$, then $\Gamma(v) \cap X = \set{a_{1}, a_{2}, a_{5}, a_{6}}$.
	
	Let $v$ be a vertex with $f(v) \geqslant 6$. No vertex has five neighbours in a copy of $H_{1}$ (as noted in \cref{rmk:4nbs}), so $v$ is adjacent to at most four of the $a_{i}$. Thus, $v$ is adjacent to at least two of the vertices of weight two, that is, to at least two of $a_{0}$, $a_{2}$, $a_{5}$. If $v$ is adjacent to all of $a_{0}, a_{2}, a_{5}$, then $G_{a_{0}}$ contains the odd circuit $va_{5}a_{6}a_{1}a_{2}$. Thus $v$ is adjacent to exactly two of $a_{0}$, $a_{2}$, and $a_{5}$.
	
	Suppose $v$ is adjacent to $a_{0}$. By symmetry we may assume that $v$ is adjacent to $a_{2}$ but not to $a_{5}$. Then $v$ is not adjacent to $a_{1}$, else $va_{0}a_{1}a_{2}$ is a $K_{4}$. Similarly $v$ cannot be adjacent to both $a_{3}$ and $a_{4}$. Hence $v$ is adjacent to $a_{0}$, $a_{2}$, $a_{6}$ and one of $a_{3}$, $a_{4}$. By symmetry, we may assume $v$ is adjacent to $a_{3}$. Now $v$ cannot be $a_{4}$ nor $a_{5}$ as it would then have five neighbours in $X$. Hence $G[a_{5}, a_{6}, a_{0}, v, a_{2}, a_{3}, a_{4}]$ is a copy of $H_{2}$.
	
	Thus $v$ is not adjacent to $a_{0}$ and so is adjacent to both $a_{2}$ and $a_{5}$. Note $v$ cannot be adjacent to both $a_{3}$, $a_{4}$ else $va_{2}a_{3}a_{4}$ is $K_{4}$, so we may assume that $v$ is adjacent to $a_{1}$. If $v$ were adjacent to one of $a_{3}$, $a_{4}$, then we may assume, by symmetry that $v$ is adjacent to $a_{4}$ and not $a_{3}$. Now $v$ cannot be $a_{0}$ nor $a_{6}$ as it would then have five neighbours in $X$. But then $G[a_{5}, a_{6}, a_{0}, a_{1}, a_{2}, v, a_{4}]$ is a copy of $H_{2}$. Therefore $v$ is adjacent to neither $a_{3}$ nor $a_{4}$. But $f(v) \geqslant 6$, so $v$ is adjacent to $a_{1}$, $a_{6}$ and thus $\Gamma(v) \cap X = \set{a_{1}, a_{2}, a_{5}, a_{6}}$.
	
	So every vertex $v$ has $f(v) \leqslant 6$ and furthermore all $v \in \Gamma(a_{0}) \cup \Gamma(a_{3}) \cup \Gamma(a_{4})$ have $f(v) \leqslant 5$. We first claim that there is $i \in \set{3, 4}$ such that $\Gamma(a_{0}, a_{2}, a_{i}) = \emptyset$. If not, there is $a'_{1} \in \Gamma(a_{0}, a_{2}, a_{3})$ and $a''_{1} \in \Gamma(a_{0}, a_{2}, a_{4})$. But then $G_{a_{2}}$ contains the odd circuit $a'_{1} a_{3} a_{4} a''_{1} a_{0}$. Similarly there is $j \in \set{3, 4}$ such that $\Gamma(a_{0}, a_{5}, a_{j}) = \emptyset$.
	
	Next we claim that there is $i \in \set{3, 4}$ with $\Gamma(a_{0}, a_{2}, a_{i}) = \Gamma(a_{0}, a_{5}, a_{i}) = \emptyset$. If not, then without loss of generality there is $a'_{1} \in \Gamma(a_{0}, a_{2}, a_{3})$ and $a'_{6} \in \Gamma(a_{0}, a_{5}, a_{4})$. Certainly, $\Gamma(a'_{1}) \cap X \neq \set{a_{1}, a_{2}, a_{5}, a_{6}}$, so $f(a'_{1}) \leqslant 5$. But $\omega(a_{0}) + \omega(a_{2}) + \omega(a_{3}) = 5$, so $\Gamma(a'_{1}) \cap X = \set{a_{0}, a_{2}, a_{3}}$. Similarly $\Gamma(a'_{6}) \cap X = \set{a_{0}, a_{4}, a_{5}}$. In particular, all of $a_{0}, \dotsc, a_{6}, a'_{1}, a'_{6}$ are distinct. But then $G[a_{0}, a'_{1}, a_{2}, a_{3}, a_{4}, a_{5}, a'_{6}]$ is a copy of $H_{2}$. Thus, without loss of generality, $\Gamma(a_{0}, a_{2}, a_{3}) = \Gamma(a_{0}, a_{5}, a_{3}) = \emptyset$. Then $\Gamma(a_{3})$, $\Gamma(a_{0}, a_{2})$, $\Gamma(a_{0}, a_{5})$ are pairwise disjoint (we already showed that $\Gamma(a_{0}, a_{2}, a_{5}) = \emptyset$). Hence the set $Y = \Gamma(a_{3}) \cup \Gamma(a_{0}, a_{2}) \cup \Gamma(a_{0}, a_{5})$ has size
	\begin{equation*}
		\abs{Y} = d(a_{3}) + d(a_{0}, a_{2}) + d(a_{0}, a_{5}) \geqslant \delta(G) + 2(2\delta(G) - \abs{G}) = 5 \delta(G) - 2 \abs{G}.
	\end{equation*}
	Now all $v \in Y$ have $f(v) \leqslant 5$. We bound $\omega(X, G)$ from both directions (recalling this notation from \cref{sec:notation}) to get
	\begin{equation*}
		10 \delta(G) \leqslant \sum_{x \in X} \omega(x) d(x) = \sum_{v \in G} f(v) \leqslant 5 \abs{Y} + 6(\abs{G} - \abs{Y}) = 6 \abs{G} - \abs{Y} \leqslant 8 \abs{G} - 5 \delta(G),
	\end{equation*}
	which contradicts $\delta(G) > 8/15 \cdot \abs{G}$.
\end{Proof}

These two claims combine to give \cref{lemma4H}.

\subsection{Ruling out sparse pairs being spokes of odd wheels}\label{sec:oddspokes}

In this subsection we make a start on the proof of \cref{main4localbip}, by ruling out the possibility that $G$ contains a sparse pair of vertices which is the spoke of an odd wheel (\cref{cor4oddwheel}). Our next three lemmas show that any sparse pair which is the missing spoke of an odd wheel is, in fact, the missing spoke of a 7-wheel. The following straightforward lemma appeared in~\cite{Illingworth2022localbipart}.

\begin{lemma}\label{lemma45wheel}
	Let $G$ be a locally bipartite graph with $\delta(G) > 1/2 \cdot \abs{G}$ and which does not contain $H_{0}$. Then $G$ does not contain a sparse pair $u, v$ with $uv$ being the missing spoke of a 5-wheel.
\end{lemma}

We will use the following technical lemma on various occasions.

\begin{lemma}\label{lemma42nbs}
	Let $G$ be a locally bipartite graph and let $u, v$ be a sparse pair of vertices in $G$. Suppose that $C$ is the shortest odd cycle which both passes through $v$ and satisfies $C \setminus \set{v} \subset \Gamma(u)$ \textup{(}i.e.\ $G[C \cup \set{u}]$ contains an odd wheel missing the spoke $uv$\textup{)}. Then every neighbour of $u$ has at most two neighbours in $C$ and if two, then they are two apart on $C$. In particular, $C$ is an induced cycle.
\end{lemma}

\begin{Proof}
	Label the configuration as follows and write $v_{0}$ for $v$ (we consider indices modulo $2k + 1$). Note $C = \set{v, v_{1}, \dotsc, v_{2k}}$.
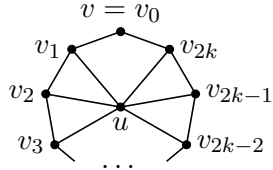
	\begin{figure}[H]
		\centering
		\begin{tikzpicture}
			\foreach \pt in {0,1,2,3,4,5,6,7,8} 
			{
				\tkzDefPoint(\pt*360/9 + 90:1){v_\pt}
			}
			\tkzDefPoint(0,0){u}
			\tkzDefMidPoint(v_3,v_4) \tkzGetPoint{A}
			\tkzDefMidPoint(v_5,v_6) \tkzGetPoint{B}
			\tkzDefMidPoint(v_4,v_5) \tkzGetPoint{C}
			
			\tkzDrawPolySeg(v_6,v_7,v_8,v_0,v_1,v_2,v_3)
			\foreach \pt in {1,2,3,6,7,8} 
			{
				\tkzDrawSegment(u,v_\pt)
			}
			\tkzDrawSegment(v_3,A)
			\tkzDrawSegment(B,v_6)
			
			\tkzDrawPoints(v_0,v_1,v_2,v_3,v_6,v_7,v_8)
			\tkzDrawPoint(u)
			\tkzLabelPoint[above](v_0){$v = v_{0}$}
			\tkzLabelPoint[left](v_1){$v_{1}$}
			\tkzLabelPoint[left](v_2){$v_{2}$}
			\tkzLabelPoint[left](v_3){$v_{3}$}
			\tkzLabelPoint[right](v_6){$v_{2k - 2}$}
			\tkzLabelPoint[right](v_7){$v_{2k - 1}$}
			\tkzLabelPoint[right](v_8){$v_{2k}$}
			\tkzLabelPoint[below](u){$u$}
			\tkzLabelPoint[above](C){$\dotsc$}
		\end{tikzpicture}
	\end{figure}
	Consider a vertex $x$ which is adjacent to $u$. Suppose that $x$ is adjacent to two vertices in $C$ which are not two apart: $x$ is adjacent to $v_{i}$ and $v_{i + r}$ where $r \in \set{1, 3, 4, \dotsc, k}$. Firstly, if $r = 1$, then either $G$ contains the $K_{4}$ $uxv_{i}v_{i + 1}$ or $G_{u, v}$ contains an edge (if one of $v_{i}$ or $v_{i + 1}$ is $v$) which contradicts the sparsity of $u, v$. Secondly, if $r > 1$ is odd, then $C' = x v_{i} v_{i + 1} \dotsb v_{i + r}$ is an odd circuit which is shorter than $C$. But then $C'$ contains an odd cycle $C''$ which is shorter than $C$. Either $C''$ is in $G_{u}$ (if $v \not \in C''$) contradicting the local bipartiteness of $G$ or we have found a shorter odd cycle than $C$ which satisfies the properties of $C$ (if $v \in C''$). Finally, if $r > 2$ is even, then $C' = xv_{i + r}v_{i + r + 1} \dotsb v_{i - 1} v_{i}$ is an odd circuit which is shorter than $C$. Again $C'$ must contain an odd cycle $C''$ which is shorter than $C$. We either obtain an odd cycle in $G_{u}$ or contradict the minimality of $C$. Hence every neighbour of $u$ has at most two neighbours in $C$ and if two, then they are $v_{i}$, $v_{i + 2}$ for some $i$.
	
	All of $v_{1}, \dotsc, v_{2k}$ are neighbours of $u$ so have two neighbours in $C$. Hence $C$ is induced.
\end{Proof}

\begin{lemma}\label{lemma4oddwheel}
	Let $G$ be a locally bipartite graph with $\delta(G) > 8/15 \cdot \abs{G}$ which does not contain $H_{0}$. Any sparse pair in $G$ that is the missing spoke of an odd wheel is the missing spoke of a 7-wheel.
\end{lemma}

\begin{Proof}
	Consider a sparse pair $u, v$ that is the missing spoke of a $(2k + 1)$-wheel. Choose the odd wheel so that $k$ is minimal. Without loss of generality, $u$ is the central vertex and $v$ is in the outer $(2k + 1)$-cycle which we call $C$. \Cref{lemma45wheel} shows that $k > 2$. We are done if we can show that $k = 3$.
	
	By \cref{lemma42nbs}, every neighbour of $u$ has at most two neighbours in $C$. All vertices have at most $2k$ neighbours in $C$ as otherwise $G$ contains a $(2k + 1)$-wheel and so is not locally bipartite. Hence,
	\begin{align*}
		(2k + 1) \delta(G) & \leqslant e(G, C) \leqslant 2 d(u) + 2k(\abs{G} - d(u)) \\
		& = 2k \abs{G} - (2k - 2) d(u) \leqslant 2k \abs{G} - (2k - 2) \delta(G),
	\end{align*}
	so
	\begin{equation*}
		\frac{8}{15} < \frac{\delta(G)}{\abs{G}} \leqslant \frac{2k}{4k - 1},
	\end{equation*}
	which implies that $k < 4$.
\end{Proof}

Thus, to rule out sparse pairs being a missing spoke of an odd wheel we only need to rule out there being a sparse pair which is the missing spoke of a 7-wheel. This is where the graph $T_{0}$ becomes relevant. We will prove the following.

\begin{proposition}\label{lemma4T0}
	Let $G$ be a locally bipartite graph which does not contain $H_{0}$. If either $\delta(G) > 7/13 \cdot \abs{G}$, or $\delta(G) > 8/15 \cdot \abs{G}$ and $G$ does not contain $T_{0}$, then $G$ does not contain a sparse pair that is the missing spoke of a 7-wheel.
\end{proposition}

\begin{Proof}
	Let $G$ be a locally bipartite graph which does not contain $H_{0}$ and either satisfies $\delta(G) > 7/13 \cdot \abs{G}$, or $\delta(G) > 8/15 \cdot \abs{G}$ and $G$ does not contain $T_{0}$. In a slight abuse of notation we will say that a vertex $u$ is \defn{sparse to a cycle $C$} if $u$ is adjacent to $\abs{C} - 1$ vertices of $C$ and is in a sparse pair with the final vertex. We are required to show that there is no vertex $u$ and no 7-cycle $C$ with $u$ sparse to $C$.
	
	By \cref{lemma45wheel}, there is no 5-cycle $C$ and vertex $u$ with $u$ sparse to $C$. Hence, if a vertex $u$ is sparse to a 7-cycle $C$, then, by \cref{lemma42nbs}, $C$ is an induced 7-cycle and any neighbour of $u$ has at most two neighbours in $C$.
	
	\begin{claim}\label{six}
		If a vertex $u$ is sparse to a 7-cycle $C = v_{0}v_{1} \dotsc v_{6}$ with the pair $u, v_{0}$ sparse, then there is some vertex which has six neighbours in $C$ and is adjacent to all of $v_{6}$, $v_{0}$, $v_{1}$.
	\end{claim}
	
	\begin{Proof}[of Claim]
		First consider the induced 4-cycle $uv_{6}v_{0}v_{1}$. The pair $u, v_{0}$ is sparse, so $\Gamma(u, v_{6}, v_{0}) = \Gamma(v_{0}, v_{1}, u) = \emptyset$. Also $\Gamma(v_{1}, u, v_{6}) = \emptyset$ as otherwise $G_{u}$ contains an odd circuit. Hence all $z \not \in \Gamma(v_{6}, v_{0}, v_{1})$ have at most two neighbours in $\set{u, v_{6}, v_{0}, v_{1}}$.
		
		\begin{figure}[H]
			\centering
			\begin{tikzpicture}
				\foreach \pt in {0,1,6} 
				{
					\tkzDefPoint(\pt*360/7 + 90:1){v_\pt}
				} 
				\tkzDefPoint(0,0){u}
				\tkzDrawPolySeg(u,v_6,v_0,v_1,u)
				\tkzDrawPoints(u,v_6,v_0,v_1)
				\tkzLabelPoint[below](u){$u$}
				\tkzLabelPoint[right](v_6){$v_{6}$}
				\tkzLabelPoint[above](v_0){$v_{0}$}
				\tkzLabelPoint[left](v_1){$v_{1}$}
			\end{tikzpicture}
		\end{figure}
		
		Thus,
		\begin{equation*}
			4 \delta(G) \leqslant e(\set{v_{0}, v_{1}, u, v_{6}}, G) \leqslant 3 \abs{\Gamma(v_{6}, v_{0}, v_{1})} + 2(\abs{G} - \abs{\Gamma(v_{6}, v_{0}, v_{1})}),
		\end{equation*}
		so $\abs{\Gamma(v_{6}, v_{0}, v_{1})} \geqslant 4 \delta(G) - 2 \abs{G}$.
		
		If the claim is false, then all vertices in $\Gamma(v_{6}, v_{0}, v_{1})$ have at most five neighbours in $C$. Also note that $\Gamma(u)$ and $\Gamma(v_{6}, v_{0}, v_{1})$ are disjoint and that all vertices have at most six neighbours in $C$ (otherwise $G$ contains a 7-wheel) so,
		\begin{align*}
			7 \delta(G) & \leqslant e(C, G) \leqslant 2 d(u) + 5\abs{\Gamma(v_{6}, v_{0}, v_{1})} + 6(\abs{G} - d(u) - \abs{\Gamma(v_{6}, v_{0}, v_{1})}) \\
			& = 6 \abs{G} - 4d(u) - \abs{\Gamma(v_{6}, v_{0}, v_{1})} \leqslant 6 \abs{G} - 4 \delta(G) - 4 \delta(G) + 2 \abs{G},
		\end{align*}
		which contradicts $\delta(G) > 8/15 \cdot \abs{G}$.
	\end{Proof}
	
	\begin{claim}\label{sixsparse}
		Let $C$ be a 7-cycle such that there is some vertex which is sparse to $C$. Then every vertex with at least six neighbours in $C$ is sparse to $C$.
	\end{claim}
	
	\begin{Proof}[of Claim]
		Let vertex $u$ be sparse to the 7-cycle $C$ and let $x$ be a vertex with at least six neighbours in $C$. Since $G$ is locally bipartite, $x$ has exactly six neighbours in $C$. Let $y$ be the vertex of $C$ to which $x$ is not adjacent. There is a vertex $u'$ which has six neighbours in $C$ and is adjacent to $y$. Indeed, if $u$ is adjacent to $y$ then take $u' = u$ and if $u$ is not adjacent to $y$, then \cref{six} gives the desired $u'$.
		
		Now $\Gamma(u',x)$ contains two consecutive vertices of $C$, so the pair $u'$, $x$ is dense. But $u'$ is adjacent to $y$ so, by \cref{lemma4dense}, the pair $y, x$ cannot be dense. In particular, $x, y$ is sparse and so $x$ is sparse to $C$.
	\end{Proof}
	
	Now fix a 7-cycle $C = v_{0} v_{1} \dotsc v_{6}$ such that there is some vertex which is sparse to $C$. Say vertex $v_{i}$ is \defn{lonely} if there is some vertex $u$ which is adjacent to all of $C \setminus \set{v_{i}}$ -- by the previous claim $u$ is sparse to $C$ and the pair $u, v_{i}$ is sparse.
	
	\begin{claim}\label{lonely}
		For all $i$, $v_{i}$ and $v_{i + 2}$ are not both lonely.
	\end{claim}
	
	\begin{Proof}[of Claim]
		Suppose for contradiction that $v_{1}$ and $v_{6}$ are both lonely: let $u_{1}$ and $u_{6}$ be sparse to $C$ with both the pairs $u_{1}, v_{1}$ and $u_{6}, v_{6}$ sparse. \Cref{lemma42nbs} shows that any neighbour of $u_{1}$ (or $u_{6}$) has at most two neighbours in $C$. Let $X = \set{u_{1}, u_{6}, v_{0}, \dotsc, v_{6}}$ and give a weighting $\omega$ to the vertices in $X$ as shown below.
		
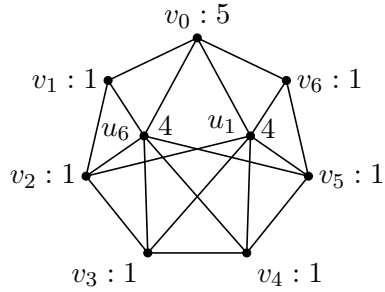
\begin{figure}[H]
			\centering
			\begin{tikzpicture}
				\foreach \pt in {0,1,...,6} 
				{
					\tkzDefPoint(\pt*360/7 + 90:1.5){v_\pt}
				}
				\tkzDefPoint(0.7,0.2){u_1}
				\tkzDefPoint(-0.7,0.2){u_6}
				\tkzDrawPolySeg(v_0,v_1,v_2,v_3,v_4,v_5,v_6,v_0)
				\foreach \pt in {0,2,3,4,5,6}
				{
					\tkzDrawSegment(u_1,v_\pt)
				}
				\foreach \pt in {0,1,2,3,4,5}
				{
					\tkzDrawSegment(u_6,v_\pt)
				}
				\tkzDrawPoints(v_0,v_...,v_6)
				\tkzDrawPoints(u_1,u_6)
				\tkzLabelPoint[above](v_0){$v_{0} \colon 5$}
				\tkzLabelPoint[left](v_1){$v_{1} \colon 1$}
				\tkzLabelPoint[left](v_2){$v_{2} \colon 1$}
				\tkzLabelPoint[below left](v_3){$v_{3} \colon 1$}
				\tkzLabelPoint[below right](v_4){$v_{4} \colon 1$}
				\tkzLabelPoint[right](v_5){$v_{5} \colon 1$}
				\tkzLabelPoint[right](v_6){$v_{6} \colon 1$}
				\tkzLabelPoint[shift={(-0.3,0.4)}](u_1){$u_{1}$}
				\tkzLabelPoint[shift={(0.3,0.4)}](u_1){$4$}
				\tkzLabelPoint[shift={(-0.3,0.4)}](u_6){$u_{6}$}
				\tkzLabelPoint[shift={(0.3,0.4))}](u_6){$4$}
			\end{tikzpicture}
		\end{figure}
		
		For each vertex $v \in G$, let $f(v)$ be the total weight of the neighbours of $v$ in $X$. We shall show that if $v \not \in \Gamma(u_{6}, u_{1}, v_{0})$, then $f(v) \leqslant 10$ and if $v \in \Gamma(u_{6}, u_{1}, v_{0})$, then $f(v) = 13$. Let $v$ be a vertex with $f(v) \geqslant 11$. It suffices to show that $v$ is adjacent to $u_{6}$, $u_{1}$, $v_{0}$ and none of $v_{1}$, \ldots, $v_{6}$.
		
		If $v$ is adjacent to neither $u_{1}$ nor $u_{6}$, then $f(v) \leqslant 11$ with equality only if $v$ is adjacent to all of $C$ which would give a 7-wheel so, in fact, $f(v) \leqslant 10$. Thus, we may assume $v$ is adjacent to $u_{1}$. But then $v$ must have at most two neighbours in $C$. If neither of these is $v_{0}$, then $f(v) \leqslant 4 + 4 + 1 + 1 = 10$. Hence we may assume $v$ is adjacent to both $v_{0}$ and $u_{1}$.
		
		Since $f(v) \geqslant 11$ and $v$ has at most two neighbours in $C$, $v$ must be adjacent to $u_{6}$ as well. Now, $v$ is adjacent to both $u_{1}$ and $v_{0}$ so, by \cref{lemma42nbs}, the only other possible neighbour of $v$ in $C$ is one of $v_{2}$ and $v_{5}$. However if $v$ is adjacent to $v_{2}$, then $G_{u_{1}}$ contains the odd circuit $v_{0}vv_{2}v_{3}\dotsc v_{6}$ while if $v$ is adjacent to $v_{5}$, then $G_{u_{6}}$ contains the odd circuit $v_{0}v_{1}\dotsc v_{5}v$. In conclusion, $v$ is adjacent to $u_{1}$, $u_{6}$, $v_{0}$, and no other vertices of $C$.
		
		Thus
		\begin{equation*}
			19 \delta(G) \leqslant \omega(X, G) = \sum_{v \in G} f(v) \leqslant 13 \abs{\Gamma(v_{0}, u_{1}, u_{6})} + 10 (\abs{G} - \abs{\Gamma(v_{0}, u_{1}, u_{6})}),
		\end{equation*}
		so
		\begin{equation}\label{eq:19}
			3 \abs{\Gamma(v_{0}, u_{1}, u_{6})} \geqslant 19 \delta(G) - 10 \abs{G}.
		\end{equation}
		Any $v \in \Gamma(v_{0}, u_{1}, u_{6})$ satisfies $f(v) = 13$ so has only one neighbour in $C$. Also any neighbour of $u_{1}$ has at most two neighbours in $C$. Hence
		\begin{align*}
			7 \delta(G) & \leqslant e(C, G) \leqslant \abs{\Gamma(v_{0}, u_{1}, u_{6})} + 2(d(u_{1}) - \abs{\Gamma(v_{0}, u_{1}, u_{6})}) + 6(\abs{G} - d(u_{1})) \\
			& = 6 \abs{G} - 4d(u_{1}) - \abs{\Gamma(v_{0}, u_{1}, u_{6})},
		\end{align*}
		so
		\begin{equation}\label{eq:11}
			\abs{\Gamma(v_{0}, u_{1}, u_{6})} \leqslant 6 \abs{G} - 11 \delta(G).
		\end{equation}
		Combining inequalities \eqref{eq:19} and \eqref{eq:11} gives
		\begin{equation*}
			19 \delta(G) - 10 \abs{G} \leqslant 3 \abs{\Gamma(v_{0}, u_{1}, u_{6})} \leqslant 18 \abs{G} - 33 \delta(G),
		\end{equation*}
		so $\delta(G) \leqslant 7/13 \cdot \abs{G}$ and so $G$ is $T_{0}$-free.
		
		Inequality \eqref{eq:19} and $\delta(G) > 8/15 \cdot \abs{G}$ show that $\Gamma(v_{0}, u_{1}, u_{6})$ is non-empty. Let $v$ be a common neighbour of $v_{0}$, $u_{1}$, $u_{6}$. As $G$ is $T_{0}$-free, $v$ must be one of the $v_{i}$. But $C = v_{0}v_{1} \dotsc v_{6}$ is induced, so $v$ must be one of $v_{1}$, $v_{6}$. This means one of the edges $u_{1}v_{1}$, $u_{6}v_{6}$ is present. However these are both sparse pairs giving the required contradiction.
	\end{Proof}
	
	We now finish the proof of \cref{lemma4T0}. By the choice of $C$, some $v_{i}$ is lonely. Without loss of generality, $v_{0}$ is lonely. By \cref{lonely}, neither $v_{2}$ nor $v_{5}$ is lonely. By \cref{six,sixsparse}, at least one of $v_{3}$ and $v_{4}$ is lonely. By symmetry, we may assume $v_{3}$ is lonely. By \cref{lonely}, $v_{1}$ is not lonely and at most one of $v_{4}$, $v_{6}$ is. Again, by symmetry, we may assume $v_{6}$ is not lonely. In conclusion, $v_{1}$, $v_{2}$, $v_{5}$, $v_{6}$ are all not lonely, $v_{0}$ and $v_{3}$ are lonely, and $v_{4}$ may or may not be.
	
	Let $U_{6}$ be the set of vertices with six neighbours in $C$. By \cref{sixsparse}, any vertex in $U_{6}$ is sparse to $C$ so cannot be adjacent to all of $v_{0}$, $v_{3}$, $v_{4}$ (else some other $v_{i}$ is lonely). In particular,
	\begin{equation*}
		U_{6} \subset \overline{\Gamma(v_{0})} \cup \overline{\Gamma(v_{3})} \cup \overline{\Gamma(v_{4})}.
	\end{equation*}
	As $v_{0}$ is lonely, there is a vertex $u$ that is sparse to $C$ with $u, v_{0}$ sparse. No two of $v_{0}, v_{3}, v_{4}$ are two apart on $C$ and so, by \cref{lemma42nbs}, any neighbour of $u$ is in at least two of $\overline{\Gamma(v_{0})}$, $\overline{\Gamma(v_{3})}$, $\overline{\Gamma(v_{4})}$ and is not in $U_{6}$. Hence
	\begin{equation*}
		\abs{U_{6}} + 2d(u) \leqslant \lvert \overline{\Gamma(v_{0})} \rvert + \lvert \overline{\Gamma(v_{3})} \rvert + \lvert \overline{\Gamma(v_{4})} \rvert \leqslant 3 \abs{G} - 3 \delta(G),
	\end{equation*}
	and so $\abs{U_{6}} \leqslant 3 \abs{G} - 5 \delta(G)$. Since every neighbour of $u$ has at most two neighbours in $C$,
	\begin{align*}
		7 \delta(G) & \leqslant e(C, G) \leqslant 6 \abs{U_{6}} + 2d(u) + 5(\abs{G} - \abs{U_{6}} - d(u)) = 5 \abs{G} + \abs{U_{6}} - 3d(u) \\
		& \leqslant 5 \abs{G} + 3 \abs{G} - 5 \delta(G) - 3 \delta(G) = 8 \abs{G} - 8 \delta(G),
	\end{align*}
	which contradicts $\delta(G) > 8/15 \cdot \abs{G}$.
\end{Proof}

\Cref{lemma45wheel,lemma4oddwheel} and \cref{lemma4T0} together give the result we need.

\begin{corollary}\label{cor4oddwheel}
	Let $G$ be a locally bipartite graph which does not contain $H_{0}$. If either $\delta(G) > 7/13 \cdot \abs{G}$, or $\delta(G) > 8/15 \cdot \abs{G}$ and $G$ does not contain $T_{0}$, then $G$ does not contain a sparse pair which is the missing spoke of an odd wheel.
\end{corollary}

\subsection{The proof of \texorpdfstring{\cref{main4localbip}}{Theorem 3.2}}\label{sec:proof4localbip}

Here we will prove \cref{main4localbip} which we restate for convenience.

\mainforlocalbip*

We start with a locally bipartite graph $G$ which does not contain $H_{0}$ and satisfies either $\delta(G) > 7/13 \cdot \abs{G}$, or $\delta(G) > 8/15 \cdot \abs{G}$ and $G$ does not contain $T_{0}$. We are required to show that $G$ is \mbox{3-colourable}. We may assume that $G$ is edge-maximal: for any sparse pair $u, v$ of $G$, the addition of $uv$ to $G$ introduces an odd wheel, a copy of $H_{0}$ or a copy of $T_{0}$. By \cref{lemma4H}, the addition of $uv$ to $G$ introduces an odd wheel, a copy of $H_{2}$ (since $G$ itself does not contain $H_{2}$) or a copy of $T_{0}$. 

Firstly, if the addition of $uv$ introduces an odd wheel, then, by \cref{cor4oddwheel}, $uv$ must be a rim of that wheel -- this case is depicted in \cref{fig:Woddrimuv} below. Secondly, if the addition of $uv$ introduces a copy of $H_{2}$, then that copy of $H_{2}$ less the edge $uv$ must not contain $H_{0}$ -- this case is depicted in \cref{fig:H2uv01,fig:H2uv12,fig:H2uv13,fig:H2uv23,fig:H2uv34} below. Finally, suppose the addition of $uv$ introduces a copy of $T_{0}$ (but not an odd wheel nor a copy of $H_{0}$). Label this copy of $T_{0}$ in $G + uv$ as follows. 

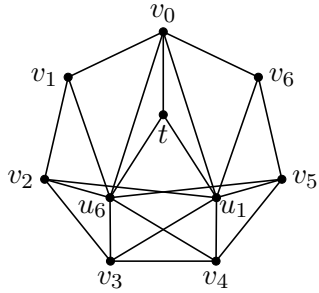
\begin{figure}[H]
	\centering
	\begin{tikzpicture}
		\foreach \pt in {0,1,...,6} 
		{
			\tkzDefPoint(\pt*360/7 + 90:1.6){v_\pt}
		}
		\tkzDefPoint(0,0.5){t}
		\tkzDefPoint(0.7,-0.6){u_1}
		\tkzDefPoint(-0.7,-0.6){u_6}
		\tkzDrawPolySeg(v_0,v_1,v_2,v_3,v_4,v_5,v_6,v_0) 
		\tkzDrawSegments(t,v_0 t,u_1 t,u_6)
		\foreach \pt in {0,2,3,4,5,6}
		{
			\tkzDrawSegment(u_1,v_\pt)
		}
		\foreach \pt in {0,1,2,3,4,5}
		{
			\tkzDrawSegment(u_6,v_\pt)
		}
		\tkzDrawPoints(v_0,v_...,v_6)
		\tkzDrawPoints(t,u_1,u_6)
		\tkzLabelPoint[above](v_0){$v_{0}$}
		\tkzLabelPoint[left](v_1){$v_{1}$}
		\tkzLabelPoint[left](v_2){$v_{2}$}
		\tkzLabelPoint[below](v_3){$v_{3}$}
		\tkzLabelPoint[below](v_4){$v_{4}$}
		\tkzLabelPoint[right](v_5){$v_{5}$}
		\tkzLabelPoint[right](v_6){$v_{6}$}
		\tkzLabelPoint[below](t){$t$}
		\tkzLabelPoint[below right = -3pt](u_1){$u_{1}$}
		\tkzLabelPoint[below left = -3pt](u_6){$u_{6}$}
	\end{tikzpicture}
\end{figure}

Note that $G + uv$ is locally bipartite and does not contain $H_{0}$ so, by \cref{lemma4dense}, for any vertex $x$, $D_{x} = \set{y \colon \text{the pair } x, y \text{ is dense}}$ is an independent set. In $G + uv$, $t \in D_{v_{1}}$ and $t$ is adjacent to $u_{1}$, so the pair $u_{1}, v_{1}$ is not dense. Also $u_{1}v_{1}$ is not an edge in $G + uv$ (else $G + uv$ contains a 7-wheel centred at $u_{1}$), so the pair $u_{1}, v_{1}$ is sparse in $G + uv$. Therefore, $u_{1}, v_{1}$ is a sparse pair in $G$. Now, by \cref{cor4oddwheel}, $G$ does not contain an odd wheel missing a sparse spoke so $uv$ must either be one of the edges $v_{i}v_{i + 1}$ or one of the edges $u_{1}v_{i}$. Similarly $uv$ must either be one of the edges $v_{i}v_{i + 1}$ or one of the edges $u_{6}v_{i}$. Thus, in fact, $uv$ must be one of the edges $v_{i}v_{i + 1}$ and by symmetry we may take $i = 0, 1, 2, 3$ -- this case is depicted in \cref{fig:T0uv01,fig:T0uv12,fig:T0uv23,fig:T0uv34} below.

Thus, in $G$, any sparse pair $u, v$ must appear in one of the following configurations (with the labels of $u$ and $v$ possibly swapped).

\begin{figure}[H]
	\centering
	\begin{subfigure}{.24\textwidth}
		\centering
		\begin{tikzpicture}
			\foreach \pt in {0,1,...,10} 
			{
				\tkzDefPoint(\pt*360/11 + 90:1){v_\pt}
			}
			\tkzDefPoint(0,0){u}
			\tkzDefMidPoint(v_1,v_2) \tkzGetPoint{A}
			\tkzDefMidPoint(v_9,v_10) \tkzGetPoint{B}
			\tkzDefMidPoint(v_1,v_10) \tkzGetPoint{C}
			
			\tkzDrawPolySeg(v_2,v_3,v_4,v_5)
			\tkzDrawPolySeg(v_6,v_7,v_8,v_9)
			\foreach \pt in {2,3,...,9} 
			{
				\tkzDrawSegment(u,v_\pt)
			}
			\tkzDrawSegment(v_2,A)
			\tkzDrawSegment(B,v_9)
			
			\tkzDrawPoints(v_2,v_...,v_9)
			\tkzDrawPoint(u)
			\tkzLabelPoint[below](v_5){$u$}
			\tkzLabelPoint[below](v_6){$v$}
			\tkzLabelPoint[below](C){$\dotsc$}
		\end{tikzpicture}
		\subcaption{}\label{fig:Woddrimuv}
	\end{subfigure}
	\begin{subfigure}{.24\textwidth}
		\centering
		\begin{tikzpicture}
			\foreach \pt in {0,1,...,6} 
			{
				\tkzDefPoint(\pt*360/7 + 90:1){v_\pt}
			} 
			\tkzDrawPolySeg(v_1,v_2,v_3,v_4,v_5,v_6,v_0)
			\tkzDrawPolySeg(v_1,v_3,v_5,v_0, v_2,v_4,v_6)
			\tkzDrawPoints(v_0,v_...,v_6)
			\tkzLabelPoint[above](v_0){$u$}
			\tkzLabelPoint[above left](v_1){$v$}
		\end{tikzpicture}
		\subcaption{}\label{fig:H2uv01}
	\end{subfigure}
	\begin{subfigure}{.24\textwidth}
		\centering
		\begin{tikzpicture}
			\foreach \pt in {0,1,...,6} 
			{
				\tkzDefPoint(\pt*360/7 + 90:1){v_\pt}
			} 
			\tkzDrawPolySeg(v_2,v_3,v_4,v_5,v_6,v_0,v_1)
			\tkzDrawPolySeg(v_1,v_3,v_5,v_0, v_2,v_4,v_6)
			\tkzDrawPoints(v_0,v_...,v_6)
			\tkzLabelPoint[above left](v_1){$u$}
			\tkzLabelPoint[left](v_2){$v$}
		\end{tikzpicture}
		\subcaption{}\label{fig:H2uv12}
	\end{subfigure}
	\begin{subfigure}{.24\textwidth}
		\centering
		\begin{tikzpicture}
			\foreach \pt in {0,1,...,6} 
			{
				\tkzDefPoint(\pt*360/7 + 90:1){v_\pt}
			} 
			\tkzDrawPolySeg(v_3,v_4,v_5,v_6,v_0,v_1,v_2)
			\tkzDrawPolySeg(v_1,v_3,v_5,v_0, v_2,v_4,v_6)
			\tkzDrawPoints(v_0,v_...,v_6)
			\tkzLabelPoint[left](v_2){$u$}
			\tkzLabelPoint[below](v_3){$v$}
		\end{tikzpicture}
		\subcaption{}\label{fig:H2uv23}
	\end{subfigure}
	
	\medskip
	
	\begin{subfigure}{.24\textwidth}
		\centering
		\begin{tikzpicture}
			\foreach \pt in {0,1,...,6} 
			{
				\tkzDefPoint(\pt*360/7 + 90:1){v_\pt}
			} 
			\tkzDrawPolySeg(v_4,v_5,v_6,v_0,v_1,v_2,v_3)
			\tkzDrawPolySeg(v_1,v_3,v_5,v_0, v_2,v_4,v_6)
			\tkzDrawPoints(v_0,v_...,v_6)
			\tkzLabelPoint[below](v_3){$u$}
			\tkzLabelPoint[below](v_4){$v$}
		\end{tikzpicture}
		\subcaption{}\label{fig:H2uv34}
	\end{subfigure}
	\begin{subfigure}{.24\textwidth}
		\centering
		\begin{tikzpicture}
			\foreach \pt in {0,1,...,6} 
			{
				\tkzDefPoint(\pt*360/7 + 90:1){v_\pt}
			} 
			\tkzDrawPolySeg(v_0,v_1,v_2,v_3,v_4,v_5,v_6,v_0)
			\tkzDrawPolySeg(v_3,v_5,v_0, v_2,v_4,v_6)
			\tkzDrawPoints(v_0,v_...,v_6)
			\tkzLabelPoint[above left](v_1){$u$}
			\tkzLabelPoint[below](v_3){$v$}
		\end{tikzpicture}
		\subcaption{}\label{fig:H2uv13}
	\end{subfigure}
	\begin{subfigure}{.24\textwidth}
		\centering
		\begin{tikzpicture}
			\foreach \pt in {0,1,...,6} 
			{
				\tkzDefPoint(\pt*360/7 + 90:1){v_\pt}
			}
			\tkzDefPoint(0,0.33){t}
			\tkzDefPoint(0.46,-0.4){u_1}
			\tkzDefPoint(-0.46,-0.4){u_6}
			\tkzDrawPolySeg(v_1,v_2,v_3,v_4,v_5,v_6,v_0) 
			\tkzDrawSegments(t,v_0 t,u_1 t,u_6)
			\foreach \pt in {0,2,3,4,5,6}
			{
				\tkzDrawSegment(u_1,v_\pt)
			}
			\foreach \pt in {0,1,2,3,4,5}
			{
				\tkzDrawSegment(u_6,v_\pt)
			}
			\tkzDrawPoints(v_0,v_...,v_6)
			\tkzDrawPoints(t,u_1,u_6)
			\tkzLabelPoint[above](v_0){$u$}
			\tkzLabelPoint[above left](v_1){$v$}
		\end{tikzpicture}
		\subcaption{}\label{fig:T0uv01}
	\end{subfigure}
	\begin{subfigure}{.24\textwidth}
		\centering
		\begin{tikzpicture}
			\foreach \pt in {0,1,...,6} 
			{
				\tkzDefPoint(\pt*360/7 + 90:1){v_\pt}
			}
			\tkzDefPoint(0,0.33){t}
			\tkzDefPoint(0.46,-0.4){u_1}
			\tkzDefPoint(-0.46,-0.4){u_6}
			\tkzDrawPolySeg(v_2,v_3,v_4,v_5,v_6,v_0,v_1) 
			\tkzDrawSegments(t,v_0 t,u_1 t,u_6)
			\foreach \pt in {0,2,3,4,5,6}
			{
				\tkzDrawSegment(u_1,v_\pt)
			}
			\foreach \pt in {0,1,2,3,4,5}
			{
				\tkzDrawSegment(u_6,v_\pt)
			}
			\tkzDrawPoints(v_0,v_...,v_6)
			\tkzDrawPoints(t,u_1,u_6)
			\tkzLabelPoint[above left](v_1){$u$}
			\tkzLabelPoint[left](v_2){$v$}
		\end{tikzpicture}
		\subcaption{}\label{fig:T0uv12}
	\end{subfigure}
	
	\medskip
	
	\begin{subfigure}{.24\textwidth}
		\centering
	\end{subfigure}
	\begin{subfigure}{.24\textwidth}
		\centering
		\begin{tikzpicture}
			\foreach \pt in {0,1,...,6} 
			{
				\tkzDefPoint(\pt*360/7 + 90:1){v_\pt}
			}
			\tkzDefPoint(0,0.33){t}
			\tkzDefPoint(0.46,-0.4){u_1}
			\tkzDefPoint(-0.46,-0.4){u_6}
			\tkzDrawPolySeg(v_3,v_4,v_5,v_6,v_0,v_1,v_2) 
			\tkzDrawSegments(t,v_0 t,u_1 t,u_6)
			\foreach \pt in {0,2,3,4,5,6}
			{
				\tkzDrawSegment(u_1,v_\pt)
			}
			\foreach \pt in {0,1,2,3,4,5}
			{
				\tkzDrawSegment(u_6,v_\pt)
			}
			\tkzDrawPoints(v_0,v_...,v_6)
			\tkzDrawPoints(t,u_1,u_6)
			\tkzLabelPoint[left](v_2){$u$}
			\tkzLabelPoint[below](v_3){$v$}
		\end{tikzpicture}
		\subcaption{}\label{fig:T0uv23}
	\end{subfigure}
	\begin{subfigure}{.24\textwidth}
		\centering
		\begin{tikzpicture}
			\foreach \pt in {0,1,...,6} 
			{
				\tkzDefPoint(\pt*360/7 + 90:1){v_\pt}
			}
			\tkzDefPoint(0,0.33){t}
			\tkzDefPoint(0.46,-0.4){u_1}
			\tkzDefPoint(-0.46,-0.4){u_6}
			\tkzDrawPolySeg(v_4,v_5,v_6,v_0,v_1,v_2,v_3) 
			\tkzDrawSegments(t,v_0 t,u_1 t,u_6)
			\foreach \pt in {0,2,3,4,5,6}
			{
				\tkzDrawSegment(u_1,v_\pt)
			}
			\foreach \pt in {0,1,2,3,4,5}
			{
				\tkzDrawSegment(u_6,v_\pt)
			}
			\tkzDrawPoints(v_0,v_...,v_6)
			\tkzDrawPoints(t,u_1,u_6)
			\tkzLabelPoint[below](v_3){$u$}
			\tkzLabelPoint[below](v_4){$v$}
		\end{tikzpicture}
		\subcaption{}\label{fig:T0uv34}
	\end{subfigure}
	\begin{subfigure}{.24\textwidth}
		\centering
	\end{subfigure}
	
	\caption{Configurations in which a sparse pair $u, v$ may appear (labels $u$ and $v$ possibly swapped).}\label{fig:cases}
\end{figure}

We will now consider a largest independent set in $G$: an independent set $I$ of size $\alpha(G)$. We will shortly show that all vertices are either in $I$ or adjacent to all of $I$.

By \cref{lemma4I}, for every $u \in I$, $I \subset D_{u} \cup \set{u}$. We now show there is set equality.

\begin{proposition}\label{lemma4DuI}
	For every $u \in I$, $I = D_{u} \cup \set{u}$
\end{proposition}

\begin{Proof}
	By \cref{lemma4dense} and the definition of dense, $D_{u} \cup \set{u}$ is an independent set. However it contains the maximal independent set $I$, so must equal it.
\end{Proof}

The following definition will be helpful.

\begin{definition}[quasidense]
	A pair of vertices $u, v$ is \defn{quasidense} if there is a sequence of vertices $u = d_{1}, d_{2}, \dotsc, d_{k}, d_{k + 1} = v$ such that all pairs $d_{i}, d_{i + 1}$ are dense $(i = 1, 2, \dotsc, k)$.
\end{definition}

\Cref{lemma4DuI} immediately implies that if the pair $u, v$ is quasidense and $u \in I$, then $v \in I$ also.

\begin{proposition}\label{lemma4adjI}
	Every vertex of $G$ is either in $I$ or adjacent to all of $I$.
\end{proposition}

\begin{Proof}
	Fix a vertex $u \in I$ and let $v$ be any other vertex which is not adjacent to $u$. It suffices to show that $v \in I$. If the pair $u, v$ is (quasi)dense, then $v \in I$, so we may assume that $u, v$ is sparse (and not quasidense). Thus $u, v$ appears in one of the configurations given in \cref{fig:cases} (with labels $u$ and $v$ possibly swapped). However in each of \cref{fig:Woddrimuv,fig:H2uv01,fig:H2uv34,fig:T0uv01,fig:T0uv12,fig:T0uv23,fig:T0uv34} the pair $u, v$ is quasidense. Hence we may assume that $u, v$ appear in one of \cref{fig:H2uv12,fig:H2uv23,fig:H2uv13}.
	
	We consider \cref{fig:H2uv12,fig:H2uv23} together. For ease we label some more of the vertices as follows.
	
	
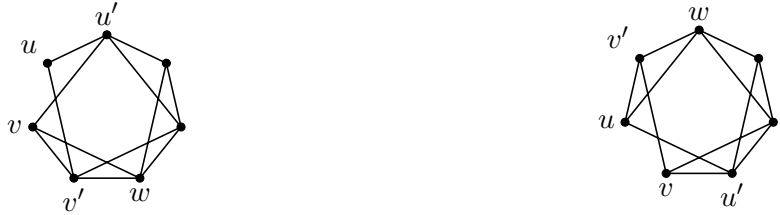
\begin{figure}[H]
		\centering
		\begin{subfigure}{.45\textwidth}
			\centering
			\begin{tikzpicture}
				\foreach \pt in {0,1,...,6} 
				{
					\tkzDefPoint(\pt*360/7 + 90:1){v_\pt}
				} 
				\tkzDrawPolySeg(v_2,v_3,v_4,v_5,v_6,v_0,v_1)
				\tkzDrawPolySeg(v_1,v_3,v_5,v_0, v_2,v_4,v_6)
				\tkzDrawPoints(v_0,v_...,v_6)
				\tkzLabelPoint[above left](v_1){$u$}
				\tkzLabelPoint[left](v_2){$v$}
				\tkzLabelPoint[above](v_0){$u'$}
				\tkzLabelPoint[below](v_3){$v'$}
				\tkzLabelPoint[below](v_4){$w$}
			\end{tikzpicture}
		\end{subfigure}
		\begin{subfigure}{.45\textwidth}
			\centering
			\begin{tikzpicture}
				\foreach \pt in {0,1,...,6} 
				{
					\tkzDefPoint(\pt*360/7 + 90:1){v_\pt}
				} 
				\tkzDrawPolySeg(v_3,v_4,v_5,v_6,v_0,v_1,v_2)
				\tkzDrawPolySeg(v_1,v_3,v_5,v_0, v_2,v_4,v_6)
				\tkzDrawPoints(v_0,v_...,v_6)
				\tkzLabelPoint[left](v_2){$u$}
				\tkzLabelPoint[below](v_3){$v$}
				\tkzLabelPoint[above left](v_1){$v'$}
				\tkzLabelPoint[below](v_4){$u'$}
				\tkzLabelPoint[above](v_0){$w$}
			\end{tikzpicture}
		\end{subfigure}
	\end{figure}
	\addtocounter{figure}{-1}
	In both cases, the pair $u', w$ is dense and so, by \cref{lemma4dense}, the pair $u', v'$ is not dense. However, $u'v'$ is not an edge, as the pair $u, v$ is sparse, and so $u', v'$ is a sparse pair. But then $uu'vv'$ is an induced 4-cycle in which both non-edges are sparse which contradicts \cref{lemma4sparse}.
	
	Finally we consider \cref{fig:H2uv13} which we label as follows. Let $X = \set{u, v, v', v_{1}, v_{2}, v_{3}, v_{4}}$ and give a weighting $\omega$ to the vertices of $X$ as shown.
	
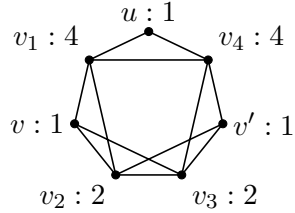
\begin{figure}[H]
		\centering
		\begin{tikzpicture}
			\foreach \pt in {0,1,...,6} 
			{
				\tkzDefPoint((\pt-1)*360/7 + 90:1){v_\pt}
			} 
			\tkzDrawPolySeg(v_0,v_1,v_2,v_3,v_4,v_5,v_6,v_0)
			\tkzDrawPolySeg(v_3,v_5,v_0, v_2,v_4,v_6)
			\tkzDrawPoints(v_0,v_...,v_6)
			\tkzLabelPoint[above](v_1){$u \colon 1$}
			\tkzLabelPoint[left](v_3){$v \colon 1$}
			\tkzLabelPoint[right](v_6){$v' \colon 1$}
			\tkzLabelPoint[above right](v_0){$v_{4} \colon 4$}
			\tkzLabelPoint[above left](v_2){$v_{1} \colon 4$}
			\tkzLabelPoint[below left](v_4){$v_{2} \colon 2$}
			\tkzLabelPoint[below right](v_5){$v_{3} \colon 2$}
		\end{tikzpicture}
	\end{figure}
	The pair $v, v'$ is dense, so, if $u, v'$ is quasidense, then $u, v$ is quasidense, a contradiction. Also, as $G$ is $H_{0}$-free, $uv'$ is not an edge. Hence the pair $u, v'$ is sparse and not quasidense. Now $v, v'$ is dense and so, by \cref{lemma4dense}, the pair $vv_{4}$ is not. But $vv_{4}$ is not an edge (else $u, v$ is dense), so $v, v_{4}$ is a sparse pair. Similarly, $v', v_{1}$ is a sparse pair. To summarise, the pairs $u, v$ and $u, v'$ are sparse and not quasidense  and the pairs $v, v_{4}$ and $v', v_{1}$ are sparse. It follows that $G[X]$ contains no more edges than shown.
	
	If a vertex $x$ has five neighbours in $X$, then it is adjacent to three consecutive vertices round the 7-cycle, so $x$ is either adjacent to a triangle or to all of $u, v_{1}, v$ or to all of $v', v_{4}, u$. The first gives a $K_{4}$ while the latter two contradict the pairs $u, v$ and $u, v'$ being sparse. Hence all vertices have at most four neighbours in $X$.
	
	For each vertex $x$ in $G$, let $f(x)$ be the total weight of the neighbours of $x$ in $X$. Now
	\begin{equation*}
		\sum_{x \in G} f(x) = \omega(X, G) \geqslant 15 \delta(G) > 8 \abs{G},
	\end{equation*}
	so some vertex $x$ has $f(x) \geqslant 9$. All vertices of $X$ have $f$ value at most eight, so $x \not \in X$. As $x$ has at most four neighbours in $X$, either $x$ is adjacent to both $v_{1}$ and $v_{4}$ or $x$ is adjacent to exactly one of $v_{1}$, $v_{4}$ and both of $v_{2}$, $v_{3}$.
	
	First suppose that $x$ is adjacent to both $v_{1}$ and $v_{4}$. As $v, v_{4}$ is sparse, $x$ is not adjacent to $v$. Similarly $x$ is not adjacent to $v'$. If $x$ is adjacent to $v_{2}$, then $x, v$ is dense (the edge $v_{1}v_{2}$). But also $u, x$ is dense (the edge $v_{1}v_{4}$), so $u, v$ is quasidense, a contradiction. Hence $x$ is not adjacent to $v_{2}$. Similarly $x$ is not adjacent to $v_{3}$, and so $f(x) = 8$, a contradiction.
	
	In the second case we may assume, by symmetry, that $x$ is adjacent to $v_{1}$, $v_{2}$, $v_{3}$ but not to $v_{4}$. Then $x, v$ and $x, v'$ are both dense pairs (the edge $v_{2}v_{3}$) so $x$ is adjacent to neither $v$ nor $v'$. Finally, if $x$ is adjacent to $u$, then $x, v_{4}$ is dense (the edge $uv_{1}$). But then the edge $v_{4}v'$ is in $D_{x}$, contradicting \cref{lemma4dense}. Hence $f(x) = 8$, a contradiction.
\end{Proof}

\begin{Proof}[of \cref{main4localbip}]
	Let $u \in I$. \Cref{lemma4adjI} gives $G[V(G) \setminus I] = G_{u}$, so $G[V(G) \setminus I]$ is 2-colourable. Using a third colour for the independent set $I$ gives a 3-colouring of $G$.
\end{Proof}

\section{Locally \texorpdfstring{$b$}{b}-partite graphs}\label{sec:localbpart}

In this section we prove \cref{main4alocalbpart} in the case $a = 1$, showing that, for $b \geqslant 3$, any locally $b$-partite graph $G$ with minimum degree greater than $(1 - 1/(b + 1/7)) \cdot \abs{G}$ is $(b + 1)$-colourable.

The proof of this will be an induction upon $b$, with some ideas from the proof of \cref{main4localbip} persisting. In this introduction to the section we generalise some of our previous ideas to the locally $b$-partite case and then we give a sketch of the proof. We first generalise dense and sparse pairs.
\begin{definition}[$b$-dense and $b$-sparse]
	A pair of \textbf{non-adjacent}, \textbf{distinct} vertices $u, v$ in a graph $G$ is \defn{$b$-dense} if $G_{u, v}$ contains a $b$-clique and is \defn{$b$-sparse} if $G_{u, v}$ does not contain a $b$-clique.
\end{definition}
This extends the notion of dense and sparse given in \cref{def:dense} -- the definitions given there are identical to those of 2-dense and 2-sparse. Note that any pair of distinct vertices is exactly one of `adjacent', `$b$-dense' or `$b$-sparse'. Another way to view being $b$-dense is being the missing edge of a $K_{b + 2}$. Locally $b$-partite graphs are $K_{b + 2}$-free so any pair of distinct vertices with a $b$-clique in their common neighbourhood must be non-adjacent and so $b$-dense. The following lemma will be helpful for lifting results from the locally bipartite case.

\begin{lemma}[lifting]\label{lemma4lifting}
	Let $b$, $s$ be positive integers and $\gamma$ any real with $b + \gamma > s$. Let $G$ be a graph with $\delta(G) > (1 - 1/(b + \gamma)) \cdot \abs{G}$. For any $s$-set $X \subset V(G)$, we have
	\begin{align*}
		\abs{G_{X}} & \geqslant s \delta(G) - (s - 1) \abs{G} > \bigl(1 - \tfrac{s}{b + \gamma}\bigr) \cdot \abs{G}, \text{ and}\\
		\delta(G_{X}) & > \big(1 - \tfrac{1}{b - s + \gamma}\big) \cdot \abs{G_{X}}.
	\end{align*}
\end{lemma}

\begin{Proof}
	Let $X = \set{x_{1}, \dotsc, x_{s}}$. Note that for each $v \in V(G)$
	\begin{equation*}
		\mathds{1}(v \in G_{X}) \geqslant \mathds{1}(vx_{1} \in E(G)) + \dotsb + \mathds{1}(vx_{s} \in E(G)) - (s - 1),
	\end{equation*}
	and summing over $v \in V(G)$ gives 
	\begin{equation*}
		\abs{G_{X}} \geqslant s \delta(G) - (s - 1) \abs{G} > s\bigl(1 - \tfrac{1}{b + \gamma}\bigr) \abs{G} - (s - 1) \abs{G} = \bigl(1 - \tfrac{s}{b + \gamma}\bigr) \cdot \abs{G}.
	\end{equation*}
	Note that $\delta(G_{X}) \geqslant \delta(G) - (\abs{G} - \abs{G_{X}}) = \abs{G_{X}} - (\abs{G} - \delta(G))$ so
	\begin{align*}
		\frac{\delta(G_{X})}{\abs{G_{X}}} & \geqslant 1 - \frac{\abs{G} - \delta(G)}{\abs{G}} \cdot \frac{\abs{G}}{\abs{G_{X}}} > 1 - \biggl(1 - \frac{\delta(G)}{\abs{G}}\biggr) \cdot \frac{1}{1 - s/(b + \gamma)} \\
		& > 1 - \frac{1}{b + \gamma} \cdot \frac{1}{1 - s/(b + \gamma)} = 1 - \frac{1}{b + \gamma - s}.\qedhere
	\end{align*}
\end{Proof}

\begin{remark}\label{remark4lifting}
	If, in \cref{lemma4lifting}, $X$ is an $s$-clique and $G$ is locally $b$-partite, then $G_{X}$ is $(b - s + 1)$-colourable and for any $u, v \in G_{X}$, if the pair $u, v$ is $b$-sparse in $G$, then $u, v$ is $(b - s)$-sparse in $G_{X}$. (The former can be seen by taking a vertex $x \in X$ and noting that $G_{x}$ is $b$-colourable and contains the $(s - 1)$-clique $X - \set{x}$ joined to $G_{X}$; the latter by noting that if vertices $u, v \in G_{X}$ form a $(b - s)$-dense pair in $G_{X}$, then they form a $b$-dense pair in $G$.)
\end{remark}

The next lemma extends \cref{lemma4dense} and partly explains why the situation for locally $b$-partite is simpler for $b \geqslant 3$ than for $b = 2$. \Cref{lemma4dense} said that in any locally bipartite $H_{0}$-free graph, the set of vertices which form a dense pair with a fixed vertex is independent. Here, in place of $H_{0}$-free, we give a minimum degree condition which guarantees this -- for $b \geqslant 3$ this minimum degree condition falls below the chromatic threshold so is, for our purposes, automatic. When $b = 2$, this minimum degree condition is 4/7, which corresponds to $\overline{C}_{7}$.

\begin{lemma}\label{lemma4bdense}
	Let $b \geqslant 2$ be an integer and $G$ be a locally $b$-partite graph with $\delta(G) > 2b/(2b + 3) \cdot \abs{G}$. For each vertex $v$ of $G$,
	\begin{equation*}
		D_{v} \coloneqq \set{u \colon \text{the pair } u, v \text{ is } b\text{-dense}}
	\end{equation*}
	is an independent set of vertices.
\end{lemma}

\begin{Proof}
	Suppose that in fact there are vertices $v$, $u_{1}$, and $u_{2}$ with both pairs $v, u_{1}$ and $v, u_{2}$ $b$-dense as well as $u_{1}$ adjacent to $u_{2}$. Let $Q_{1}$ and $Q_{2}$ be $b$-cliques in $G_{v, u_{1}}$ and $G_{v, u_{2}}$ respectively. Choose $Q_{1}$ and $Q_{2}$ so that $\ell = \abs{V(Q_{1}) \cap V(Q_{2})}$ is maximal.
	
	Firstly if $\ell \geqslant 1$, then fix $y \in V(Q_{1}) \cap V(Q_{2})$. Now $G_{y}$ contains
	
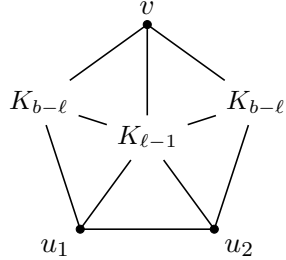
\begin{figure}[H]
		\centering
		\begin{tikzpicture}
			\foreach \pt in {0,1,...,5} 
			{
				\tkzDefPoint(\pt*360/5 + 90:1.5){v_\pt}
			}
			\node(K1) at (v_1){$K_{b - \ell}$};
			\node(K2) at (v_4){$K_{b - \ell}$};
			\node(K) at (0,0){$K_{\ell - 1}$};
			\draw[line width = 0.6pt](v_0) -- (K1);
			\draw[line width = 0.6pt](K1) -- (v_2);
			\tkzDrawSegment(v_2,v_3)
			\draw[line width = 0.6pt](v_3) -- (K2);
			\draw[line width = 0.6pt](K2) -- (v_0);
			\draw[line width = 0.6pt](K) -- (v_0);
			\draw[line width = 0.6pt](K) -- (K1);
			\draw[line width = 0.6pt](K) -- (v_2);
			\draw[line width = 0.6pt](K) -- (v_3);
			\draw[line width = 0.6pt](K) -- (K2);
			\tkzDrawPoints(v_0,v_2,v_3)
			
			\tkzLabelPoint[above](v_0){$v$}
			\tkzLabelPoint[below left](v_2){$u_{1}$}
			\tkzLabelPoint[below right](v_3){$u_{2}$}
		\end{tikzpicture}
	\end{figure}
	The pair $v, u_{1}$ is $(b - 1)$-dense in $G_{y}$ so in any $b$-colouring of $G_{y}$, $v$ and $u_{1}$ are the same colour. Similarly in any $b$-colouring of $G_{y}$, $v$ and $u_{2}$ are the same colour. In particular, $G_{y}$ is not $b$-colourable which contradicts the local $b$-colourability of $G$.
	
	Hence $\ell = 0$. Let $X = V(Q_{1}) \cup V(Q_{2}) \cup \set{v, u_{1}, u_{2}}$ which is a set of $2b + 3$ vertices. As $\delta(G) > 2b/(2b + 3) \cdot \abs{G}$, some vertex has at least $2b + 1$ neighbours in $X$ -- call this vertex $x$. As $G_{x}$ is $K_{b + 1}$-free, $x$ has a non-neighbour $x_{1} \in V(Q_{1}) \cup \set{u_{1}}$ and a non-neighbour $x_{2} \in V(Q_{2}) \cup \set{u_{2}}$. These must be the only non-neighbours of $x$ in $X$. In particular, $x$ is adjacent to $v$ and so, as $G_{x}$ is $K_{b + 1}$-free, $x_{1}$ must be in $V(Q_{1})$ and $x_{2}$ must be in $V(Q_{2})$. In particular, $Q'_{1} = Q_{1} - \set{x_{1}} + \set{x}$ is a $b$-clique in $G_{v, u_{1}}$ and $Q'_{2} = Q_{2} - \set{x_{2}} + \set{x}$ is a $b$-clique in $G_{v, u_{2}}$. But $\abs{V(Q'_{1}) \cap V(Q'_{2})} = 1$ which contradicts the maximality of $\ell$.
\end{Proof}

Note that
\begin{equation*}
	1 - \tfrac{1}{b} \geqslant \tfrac{2b}{2b + 3},
\end{equation*}
for all $b \geqslant 3$. We are only interested in locally $b$-partite graphs with $\delta(G) > (1 - 1/b) \cdot \abs{G}$ (the chromatic threshold) and the conclusion of \cref{lemma4bdense} holds for all such graphs.

We are now ready to give a sketch of \cref{main4alocalbpart} for locally $b$-partite graphs. Let $G$ be a locally $b$-partite graph ($b \geqslant 3$) with $\delta(G) > (1 - 1/(b + 1/7)) \abs{G}$. We aim to show that $G$ is $(b + 1)$-colourable. We may assume that $G$ is edge-maximal: for any missing edge $uv$, $G' = G + uv$ is not locally $b$-partite. Using induction on $b$ and the lifting lemma we will show there is some $(b - 2)$-clique $K$ with $G'_{K}$ not 3-colourable.

We will rule out the configurations where at least one of $u, v$ is in $K$ and so $u, v \not\in K$. Since $G$ is locally $b$-partite, $G_{K}$ is 3-colourable and also, using the lifting lemma, we will obtain $\delta(G_{K}) > 8/15 \cdot \abs{G_{K}}$. As $G_{K}$ is 3-colourable while $G'_{K}$ is not, we must have $G'_{K} = G_{K} + uv$. Further, by \cref{main4localbip}, the addition of $uv$ to $G_{K}$ introduces an odd wheel, a copy of $H_{2}$ or a copy of $T_{0}$. We are then in a position to use \cref{lemma4bdense} and finish in similar (although more involved) way to \cref{sec:proof4localbip}.

The ruling out of configurations is done in \cref{sec:bconfiguration} while the rest of the proof is presented in \cref{sec:proof4localbpart}.

\subsection{Ruling out configurations}\label{sec:bconfiguration}

In this section, we will rule out various configurations from locally $b$-partite graphs in a similar vein to \cref{sec:oddspokes}. From now on we will use $C_{\text{odd}}$ to denote any odd cycle.

\begin{proposition}\label{config:uvCodd}
	Fix an integer $b \geqslant 3$, let $G$ be a locally $b$-partite graph with $\delta(G) > (1 - 1/(b + 1/7)) \cdot \abs{G}$ and let $u, v$ be a $b$-sparse pair in $G$. Then $G + uv$ does not contain a $K_{b - 1} + C_{\text{odd}}$ where at least one of $u$, $v$ is in the $K_{b - 1}$.
\end{proposition}

\begin{Proof}
	We will prove this for $b = 3$ and then use the lifting lemma for larger $b$.
	
	For $b = 3$, $G$ is a locally tripartite graph with $\delta(G) > 15/22 \cdot \abs{G} > 2/3 \cdot \abs{G}$ and $u, v$ is a 3-sparse pair in $G$. Suppose the conclusion does not hold. Each neighbourhood of $G$ is 3-colourable so does not contain an odd wheel. In particular, $G$ does not contain $K_{2} + C_{\text{odd}}$, and so $G$ contains one of the following two configurations (corresponding to whether only one of $u, v$ is in the clique or they both are). In the left figure, $x$ is adjacent to all vertices inside the ring, and in the right figure, $u$ and $v$ are adjacent to all vertices inside the ring. In the future, we will use rings in this way.
	
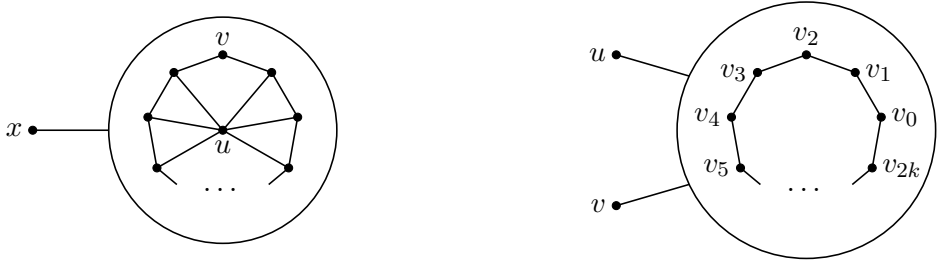
\begin{figure}[H]
		\centering
		\begin{subfigure}{.45\textwidth}
			\centering
			\begin{tikzpicture}
				\foreach \pt in {0,1,2,3,4,5,6,7,8} 
				{
					\tkzDefPoint(\pt*360/9 + 90:1){v_\pt}
				}
				\tkzDefPoint(0,0){u}
				\tkzDefMidPoint(v_3,v_4) \tkzGetPoint{A}
				\tkzDefMidPoint(v_5,v_6) \tkzGetPoint{B}
				\tkzDefMidPoint(v_4,v_5) \tkzGetPoint{C}
				
				\tkzDrawPolySeg(v_6,v_7,v_8,v_0,v_1,v_2,v_3)
				\foreach \pt in {1,2,3,6,7,8} 
				{
					\tkzDrawSegment(u,v_\pt)
				}
				\tkzDrawSegment(v_3,A)
				\tkzDrawSegment(B,v_6)
				
				\tkzDrawPoints(v_0,v_1,v_2,v_3,v_6,v_7,v_8)
				\tkzDrawPoint(u)
				\tkzLabelPoint[above](v_0){$v$}
				\tkzLabelPoint[below](u){$u$}
				\tkzLabelPoint[above](C){$\dotsc$}
				
				\tkzDefShiftPoint[u](180:1.5cm){x'}
				\tkzDrawCircle[black, line width = 0.6pt](u,x')
				\tkzDefPoint(-2.5,0){x}
				\tkzDrawSegment(x,x')
				\tkzDrawPoint(x)
				\tkzLabelPoint[left](x){$x$}
			\end{tikzpicture}
		\end{subfigure}
		\begin{subfigure}{.45\textwidth}
			\centering
			\begin{tikzpicture}
				\foreach \pt in {0,1,2,3,4,5,6,7,8} 
				{
					\tkzDefPoint(\pt*360/9 + 90:1){v_\pt}
				}
				\tkzDefMidPoint(v_3,v_4) \tkzGetPoint{A}
				\tkzDefMidPoint(v_5,v_6) \tkzGetPoint{B}
				\tkzDefMidPoint(v_4,v_5) \tkzGetPoint{C}
				
				\tkzDrawPolySeg(v_6,v_7,v_8,v_0,v_1,v_2,v_3)
				\tkzDrawSegment(v_3,A)
				\tkzDrawSegment(B,v_6)
				
				\tkzDrawPoints(v_0,v_1,v_2,v_3,v_6,v_7,v_8)
				\tkzLabelPoint[above](C){$\dotsc$}
				
				\tkzDefPoint(0,0){O}
				\tkzDefShiftPoint[u](155:1.7cm){u'}
				\tkzDefShiftPoint[u](205:1.7cm){v'}
				\tkzDrawCircle[black, line width = 0.6pt](O,u')
				\tkzDefPoint(-2.5,1){u}
				\tkzDefPoint(-2.5,-1){v}
				\tkzDrawSegments(u,u' v,v')
				\tkzDrawPoints(u,v)
				\tkzLabelPoint[left](u){$u$}
				\tkzLabelPoint[left](v){$v$}
				
				\tkzLabelPoint[above](v_0){$v_{2}$}
				\tkzLabelPoint[left](v_1){$v_{3}$}
				\tkzLabelPoint[left](v_2){$v_{4}$}
				\tkzLabelPoint[left](v_3){$v_{5}$}
				\tkzLabelPoint[right](v_6){$v_{2k}$}
				\tkzLabelPoint[right](v_7){$v_{0}$}
				\tkzLabelPoint[right](v_8){$v_{1}$}
			\end{tikzpicture}
		\end{subfigure}
	\end{figure}
	\addtocounter{figure}{-1}
	We deal with the left-hand configuration first. Odd wheels, $H_{0}$, and $T_{0}$ are not 3-colourable and so $G_{x}$ is locally bipartite, $H_{0}$-free and $T_{0}$-free. Applying \cref{lemma4lifting} with $X = \set{x}$ and $\gamma = 1/7$ gives $\delta(G_{x}) > (1 - 1/(2 + 1/7)) \cdot \abs{G_{x}} = 8/15 \cdot \abs{G_{x}}$ so, by \cref{cor4oddwheel}, $G_{x}$ does not contain a sparse pair which is the missing spoke of an odd wheel. In particular, the pair $u, v$ is dense in $G_{x}$ so must be 3-dense in $G$, a contradiction.
	
	Now consider the right-hand configuration. We may assume that $k$ is minimal. As the pair $u, v$ is 3-sparse, $k > 1$. Let $X = \set{u, v, v_{0}, v_{1}, \dotsc, v_{2k}}$. We consider indices modulo $2k + 1$.
	
	First consider a vertex $x$ adjacent to both $u$ and $v$. We claim that $x$ is adjacent to at most two of the $v_{i}$ so has at most four neighbours in $X$. For $r \not \in \set{0, \pm 2}$, $x$ cannot be adjacent to both $v_{i}$ and $v_{i + r}$. Indeed if $r = \pm 1$, then $xv_{i}v_{i + r}$ is a triangle in $G_{u,v}$ while other $r$ give a shorter odd cycle in $G_{u,v}$.
	
	Now consider any other vertex $x$: this is not adjacent to at least one of $u$ or $v$. If this is the only non-neighbour of $x$ in $X$ then $G_{x}$ contains an odd-wheel, which is not 3-colourable. Thus all vertices have at most $2k + 1$ neighbours in $X$. Hence,
	\begin{align*}
		(2k + 3) \delta(G) & \leqslant e(G, X) \leqslant 4 \abs{G_{u, v}} + (2k + 1)(\abs{G} - \abs{G_{u, v}}) \\
		& \leqslant (2k + 1) \abs{G} - (2k - 3)(2 \delta(G) - \abs{G}) = (4k -2) \abs{G} - (4k - 6) \delta(G),
	\end{align*}
	which contradicts $\delta(G) > 2/3 \cdot \abs{G}$.
	
	Suppose now that $b \geqslant 4$ and that $G + uv$ does contain a $K_{b - 1} + C_{\text{odd}}$ where at least one of $u, v$ is in the $K_{b - 1}$. The $(b - 1)$-clique contains a $(b - 3)$-clique, $L$, with $u, v \not \in L$. Thus $L$ is a $(b - 3)$-clique in $G$.
	
	Now $G_{L}$ is a 4-colourable (so locally tripartite) graph in which $u, v$ is 3-sparse (by \cref{remark4lifting}). Applying \cref{lemma4lifting} with $X = L$, $\gamma = 1/7$ gives $\delta(G_{L}) > (1 - 1/(b - (b - 3) + 1/7)) \cdot \abs{G_{L}} = 15/22 \cdot \abs{G_{L}}$. Finally, $G_{L} + uv$ contains a $K_{2} + C_{\text{odd}}$ where at least one of $u$, $v$ is in the 2-clique. This contradicts the result for $b = 3$.
\end{Proof}

\begin{proposition}\label{config:uvH0}
	Fix an integer $b \geqslant 3$, let $G$ be a locally $b$-partite graph with $\delta(G) > (1 - 1/(b + 1/7)) \cdot \abs{G}$ and let $u, v$ be a $b$-sparse pair in $G$. Then $G + uv$ does not contain a $K_{b - 2} + H_{0}$ where at least one of $u$, $v$ is in the $K_{b - 2}$.
\end{proposition}

\begin{Proof}
	We split into two cases depending upon whether only one of $u$, $v$ is in the $(b - 2)$-clique or both are. In each case we will prove the result for small $b$ and then use the lifting lemma for larger $b$. 
	
	Suppose only one of $u$, $v$ is in the $(b - 2)$-clique. We first prove the result for $b = 3$: $G$ is locally tripartite, the pair $u, v$ is 3-sparse, and $\delta(G) > 15/22 \cdot \abs{G} > 2/3 \cdot \abs{G}$. Suppose the conclusion does not hold. Since $G$ does not contain $K_{1} + H_{0}$, $G$ contains the configuration shown in \cref{fig:H0uv3} where $u$ is adjacent to all of the copy of $H_{0}$ except for one vertex ($v$) to which is it 3-sparse.
	\begin{figure}[H]
		\centering
		\begin{subfigure}{.45\textwidth}
			\centering
			\begin{tikzpicture}
				\foreach \pt in {0,1,...,6} 
				{
					\tkzDefPoint(\pt*360/7 + 90:1){v_\pt}
				} 
				\tkzDrawPolySeg(v_0,v_1,v_2,v_3,v_4,v_5,v_6,v_0)
				\tkzDrawPolySeg(v_5,v_0,v_2)
				\tkzDrawPolySeg(v_1,v_3)
				\tkzDrawPolySeg(v_4,v_6)
				\tkzDrawPoints(v_0,v_...,v_6)
				
				\tkzDefPoint(0,0){O}
				
				\tkzDefShiftPoint[O](180:1.2cm){u'}
				\tkzDrawCircle[black, line width = 0.6pt](O,u')
				\tkzDefPoint(-2.5,0){u}
				\tkzDrawSegment[densely dashed](u,u')
				\tkzDrawPoint(u)
				\tkzLabelPoint[left](u){$u$}
			\end{tikzpicture}
			\subcaption{}\label{fig:H0uv3}
		\end{subfigure}
		\begin{subfigure}{.45\textwidth}
			\centering
			\begin{tikzpicture}
				\foreach \pt in {0,1,...,6} 
				{
					\tkzDefPoint(\pt*360/7 + 90:1){v_\pt}
				} 
				\tkzDrawPolySeg(v_0,v_1,v_2,v_3,v_4,v_5,v_6,v_0) 
				\tkzDrawPolySeg(v_5,v_0,v_2)
				\tkzDrawPolySeg(v_1,v_3)
				\tkzDrawPolySeg(v_4,v_6)
				\tkzDrawPoints(v_0,v_...,v_6)
				
				\tkzDefPoint(0,0){O}
				
				\tkzDefShiftPoint[O](145:1.2cm){u'}
				\tkzDefShiftPoint[O](215:1.2cm){v'}
				\tkzDrawCircle[black, line width = 0.6pt](O,u')
				\tkzDefPoint(-2.5,1){u}
				\tkzDefPoint(-2.5,-1){v}
				\tkzDrawSegments(u,u' v,v')
				\tkzDrawPoints(u,v)
				\tkzLabelPoint[left](u){$u$}
				\tkzLabelPoint[left](v){$v$}
			\end{tikzpicture}
			\subcaption{}\label{fig:H0uv4}
		\end{subfigure}
		\caption{Configurations in \cref{config:uvH0}.}
	\end{figure}
	Let $X$ be the vertex set of the $H_{0}$. First consider a vertex $x$ adjacent to $u$: $x$ cannot be adjacent to a triangle or 5-cycle in $G[X]$ otherwise $G + uv$ contains a 2-clique $ux$ joined to an odd cycle which contradicts \cref{config:uvCodd}. In particular, any neighbour of $u$ has at most four neighbours in $X$.
	
	Consider any vertex $x$: $\chi(G[X]) = 4$ so $x$ has at most six neighbours in $X$. Hence
	\begin{equation*}
		7 \delta(G) \leqslant e(G,X) \leqslant 4 d(u) + 6(\abs{G} - d(u)) \leqslant 6 \abs{G} - 2\delta(G),
	\end{equation*}
	which contradicts $\delta(G) > 15/22 \cdot \abs{G} > 2/3 \cdot \abs{G}$.
	
	Now let $b \geqslant 4$ and suppose $G + uv$ does contain a $K_{b - 2} + H_{0}$ where exactly one of $u$, $v$ (say $u$) is in the $(b - 2)$-clique. Graph $G$ is locally $b$-partite so does not contain $K_{b - 2} + H_{0}$, so $v$ is in the copy of $H_{0}$. Let $L$ be the $(b - 2)$-clique without $u$: $L$ is a $(b - 3)$-clique in $G$ and so $G_{L}$ is 4-colourable and so locally tripartite. Also, by \cref{lemma4lifting}, $\delta(G_{L}) > 15/22 \cdot \abs{G_{L}}$ and, by \cref{remark4lifting}, $u, v$ is a 3-sparse pair in $G_{L}$. Finally $G_{L} + uv$ contains $u + H_{0}$, which contradicts the result we just proved for $b = 3$.
	
	Now consider the second case, where both $u$ and $v$ are in the $(b - 2)$-clique: this means that $b \geqslant 4$. We first prove the result for $b = 4$: $G$ is locally 4-partite, the pair $u, v$ is 4-sparse, and $\delta(G) > 22/29 \cdot \abs{G} > 3/4 \cdot \abs{G}$. If the result is false, then $G$ contains the configuration shown in \cref{fig:H0uv4}.
	
	Let $X = V(H_{0}) \cup \set{u}$ -- a set of eight vertices. First consider a vertex $x$ adjacent to both $u$ and $v$: $x$ cannot be adjacent to a triangle or 5-cycle in $V(H_{0})$ as this would contradict \cref{config:uvCodd}. Hence $x$ has at most four neighbours in $V(H_{0})$ and so at most five in $X$.
	
	Since $G[X]$ is not 4-colourable, all vertices have at most seven neighbours in $X$. Hence
	\begin{align*}
		8 \delta(G) & \leqslant 5 \abs{G_{u, v}} + 7(\abs{G} - \abs{G_{u, v}}) \\
		& \leqslant 7 \abs{G} - 2(2 \delta(G) - \abs{G}) = 9 \abs{G} - 4 \delta(G),
	\end{align*}
	which contradicts $\delta(G) > 22/29 \cdot \abs{G} > 3/4 \cdot \abs{G}$.
	
	Finally suppose $b \geqslant 5$ and $G + uv$ does contain a $K_{b - 2} + H_{0}$ where both of $u$, $v$ are in the $(b - 2)$-clique. Let $L$ be the $(b - 2)$-clique without $u$, $v$: $L$ is a $(b - 4)$-clique in $G$ so, by \cref{remark4lifting}, $G_{L}$ is 5-colourable and so locally 4-partite. Also, by \cref{lemma4lifting}, $\delta(G_{L}) > 22/29 \cdot \abs{G_{L}}$ and $u, v$ is a 4-sparse pair in $G_{L}$. Finally, $G_{L} + uv$ contains $uv + H_{0}$, which contradicts the result we obtained for $b = 4$.
\end{Proof}

\begin{proposition}\label{config:uvT0}
	Fix an integer $b \geqslant 3$, let $G$ be a locally $b$-partite graph with $\delta(G) > (1 - 1/(b + 1/7)) \cdot \abs{G}$ and let $u,v$ be a $b$-sparse pair in $G$. Then $G + uv$ does not contain a $K_{b - 2} + T_{0}$ where at least one of $u$, $v$ is in the $(b - 2)$-clique.
\end{proposition}

\begin{Proof}
	Again we split into two cases depending upon whether only one of $u$, $v$ is in the $(b - 2)$-clique or both are and in each case we will prove the result for small $b$ and then use the lifting lemma for larger $b$.
	
	Suppose only one of $u$, $v$ is in the $(b - 2)$-clique. We first prove the result for $b = 3$: $G$ is locally tripartite, the pair $u, v$ is 3-sparse, and $\delta(G) > 15/22 \cdot \abs{G} > 2/3 \cdot \abs{G}$. Suppose the conclusion does not hold. Since $G$ does not contain $K_{1} + T_{0}$, $G$ contains the configuration shown in \cref{fig:T0uv3} (labels have been added for convenience) where $u$ is adjacent to all of the copy of $T_{0}$ except for one vertex ($v$) to which it is 3-sparse.
	\begin{figure}[H]
		\centering
		\begin{subfigure}{.45\textwidth}
			\centering
			\begin{tikzpicture}
				\foreach \pt in {0,1,...,6} 
				{
					\tkzDefPoint(\pt*360/7 + 90:1.6){v_\pt}
				}
				\tkzDefPoint(0,0.5){t}
				\tkzDefPoint(0.7,-0.6){u_1}
				\tkzDefPoint(-0.7,-0.6){u_6}
				\tkzDrawPolySeg(v_0,v_1,v_2,v_3,v_4,v_5,v_6,v_0) 
				\tkzDrawSegments(t,v_0 t,u_1 t,u_6)
				\foreach \pt in {0,2,3,4,5,6}
				{
					\tkzDrawSegment(u_1,v_\pt)
				}
				\foreach \pt in {0,1,2,3,4,5}
				{
					\tkzDrawSegment(u_6,v_\pt)
				}
				\tkzDrawPoints(v_0,v_...,v_6)
				\tkzDrawPoints(t,u_1,u_6)
				\tkzLabelPoint[above](v_0){$v_{0}$}
				\tkzLabelPoint[left](v_1){$v_{1}$}
				\tkzLabelPoint[left](v_2){$v_{2}$}
				\tkzLabelPoint[below](v_3){$v_{3}$}
				\tkzLabelPoint[below](v_4){$v_{4}$}
				\tkzLabelPoint[right](v_5){$v_{5}$}
				\tkzLabelPoint[right](v_6){$v_{6}$}
				\tkzLabelPoint[below](t){$t$}
				\tkzLabelPoint[below right = -3pt](u_1){$u_{1}$}
				\tkzLabelPoint[below left = -3pt](u_6){$u_{6}$}
				
				\tkzDefPoint(0,0){O}
				
				\tkzDefShiftPoint[O](180:2.2cm){u'}
				\tkzDrawCircle[black, line width = 0.6pt](O,u')
				\tkzDefPoint(-3.5,0){u}
				\tkzDrawSegment[densely dashed](u,u')
				\tkzDrawPoint(u)
				\tkzLabelPoint[left](u){$u$}
			\end{tikzpicture}
			\subcaption{}\label{fig:T0uv3}
		\end{subfigure}
		\begin{subfigure}{.45\textwidth}
			\centering
			\begin{tikzpicture}
				\foreach \pt in {0,1,...,6} 
				{
					\tkzDefPoint(\pt*360/7 + 90:1.6){v_\pt}
				}
				\tkzDefPoint(0,0.5){t}
				\tkzDefPoint(0.7,-0.6){u_1}
				\tkzDefPoint(-0.7,-0.6){u_6}
				\tkzDrawPolySeg(v_0,v_1,v_2,v_3,v_4,v_5,v_6,v_0) 
				\tkzDrawSegments(t,v_0 t,u_1 t,u_6)
				\foreach \pt in {0,2,3,4,5,6}
				{
					\tkzDrawSegment(u_1,v_\pt)
				}
				\foreach \pt in {0,1,2,3,4,5}
				{
					\tkzDrawSegment(u_6,v_\pt)
				}
				\tkzDrawPoints(v_0,v_...,v_6)
				\tkzDrawPoints(t,u_1,u_6)
				\tkzLabelPoint[above](v_0){$v_{0}$}
				\tkzLabelPoint[left](v_1){$v_{1}$}
				\tkzLabelPoint[left](v_2){$v_{2}$}
				\tkzLabelPoint[below](v_3){$v_{3}$}
				\tkzLabelPoint[below](v_4){$v_{4}$}
				\tkzLabelPoint[right](v_5){$v_{5}$}
				\tkzLabelPoint[right](v_6){$v_{6}$}
				\tkzLabelPoint[below](t){$t$}
				\tkzLabelPoint[below right = -3pt](u_1){$u_{1}$}
				\tkzLabelPoint[below left = -3pt](u_6){$u_{6}$}
				
				\tkzDefPoint(0,0){O}
				
				\tkzDefShiftPoint[O](160:2.2cm){u'}
				\tkzDefShiftPoint[O](200:2.2cm){v'}
				\tkzDrawCircle[black, line width = 0.6pt](O,u')
				\tkzDefPoint(-3.5,1){u}
				\tkzDefPoint(-3.5,-1){v}
				\tkzDrawSegments(u,u' v,v')
				\tkzDrawPoints(u,v)
				\tkzLabelPoint[left](u){$u$}
				\tkzLabelPoint[left](v){$v$}
			\end{tikzpicture}
			\subcaption{}\label{fig:T0uv4}
		\end{subfigure}
		\caption{Configurations in \cref{config:uvT0}.}
	\end{figure}	
	We first show that at least one of the pairs $u_{1}, v_{1}$ and $u_{6}, v_{6}$ is 3-sparse. Neither of these is an edge as otherwise $G + uv$ contains $u$ joined to a 7-wheel which contradicts \cref{config:uvCodd}. If $v \not \in \set{u_{1}, u_{6}, v_{3}, v_{4}}$, then $u_{1}, u_{6}$ is 3-dense (triangle $uv_{3}v_{4}$) and so $u_{1}, v_{1}$ is 3-sparse, by \cref{lemma4bdense}. On the other hand, suppose $v \in \set{u_{1}, u_{6}, v_{3}, v_{4}}$ -- by symmetry we may assume $v \neq u_{6}$. Then $v_{1}, t$ is 3-dense (triangle $uv_{0}u_{6}$) and so $u_{1}, v_{1}$ is 3-sparse, by \cref{lemma4bdense}. From now on, we will assume that the pair $u_{1}, v_{1}$ is 3-sparse in $G$.
	
	Let $G' = G + uv$ -- we work in $G'$. From \cref{config:uvCodd}, $G'_{u}$ contains no odd wheel, i.e.\ is locally bipartite, and from \cref{config:uvH0}, $G'_{u}$ is $H_{0}$-free. Now $u_{1}v_{1}$ is not an edge (else there is a 7-wheel in $G'_{u}$) and the pair $v_{1}, t$ is 2-dense in $G'_{u}$ so $u_{1}, v_{1}$ is 2-sparse in $G'_{u}$, by \cref{lemma4dense}. Note that $\delta(G') > 2/3 \cdot \abs{G'}$ and so applying \cref{lemma4lifting} with $b = 3$, $\gamma = 0$, $s = 1$, and $X = \set{u}$ gives
	\begin{equation*}
		\delta(G'_{u}) > \big(1 - \tfrac{1}{3 - 1}) \cdot \lvert G'_{u} \rvert = 1/2 \cdot \lvert G'_{u} \rvert.
	\end{equation*}
	Hence we may apply \cref{lemma45wheel} to show that $u_{1}v_{1}$ is not a missing spoke of a 5-wheel in $G'_{u}$. Hence, by \cref{lemma42nbs}, any neighbour of $u_{1}$ in $G'_{u}$ is adjacent to at most two of the $v_{i}$. In particular, any common neighbour of $u$ and $u_{1}$ in $G'$ is adjacent to at most two of the $v_{i}$.
	
	Let $X = \set{u,u_{1},v_{0},v_{1}, \dotsc, v_{6}}$. What we have just shown is that any common neighbour of $u$ and $u_{1}$ in $G'$ has at most four neighbours in $X$. Consider a vertex $x$ which is not adjacent to both $u$ and $u_{1}$: $x$ cannot be adjacent to all of $X \setminus \set{u_{1}}$ otherwise $G'_{u, x}$ contains a 7-cycle which contradicts \cref{config:uvCodd}. Also $x$ cannot be adjacent to all of $X \setminus \set{u}$ as otherwise $G'_{x} = G_{x}$ contains a 7-wheel missing the spoke $u_{1}v_{1}$ which is 3-sparse in $G$. This again contradicts \cref{config:uvCodd}. Hence, any vertex has at most seven neighbours in $X$. Thus,
	\begin{align*}
		9 \delta(G) & \leqslant 9 \delta(G') \leqslant 4 \lvert G'_{u, u_{1}} \rvert + 7(\abs{G'} - \lvert G'_{u, u_{1}} \rvert) = 7 \abs{G} - 3 \lvert G'_{u, u_{1}} \rvert \\
		& \leqslant 7 \abs{G} - 3 \abs{G_{u, u_{1}}} \leqslant 7 \abs{G} - 3(2\delta(G) - \abs{G}) = 10 \abs{G} - 6 \delta(G),
	\end{align*}
	which contradicts $\delta(G) > 2/3 \cdot \abs{G}$.
	
	Now let $b \geqslant 4$ and suppose $G + uv$ does contain a $K_{b - 2} + T_{0}$ where exactly one of $u$, $v$ (say $u$) is in the $(b - 2)$-clique. Graph $G$ is locally $b$-partite so does not contain $K_{b - 2} + T_{0}$ so $v$ is in the copy of $T_{0}$. Let $L$ be the $(b - 2)$-clique without $u$. We make use of \cref{remark4lifting}: $L$ is a $(b - 3)$-clique in $G$ and so $G_{L}$ is 4-colourable and so locally tripartite. Also, by \cref{lemma4lifting}, $\delta(G_{L}) > 15/22 \cdot \abs{G_{L}}$ and $u, v$ is a 3-sparse pair in $G_{L}$. Finally $G_{L} + uv$ contains $u + T_{0}$, which contradicts the result we just proved for $b = 3$.
	
	Now consider the second case, where both $u$ and $v$ are in the $(b - 2)$-clique: this means that $b \geqslant 4$. We first prove the result for $b = 4$: $G$ is locally 4-partite, the pair $u, v$ is 4-sparse, and $\delta(G) > 22/29 \cdot \abs{G} > 3/4 \cdot \abs{G}$. If the result is false, then $G$ contains the configuration shown in \cref{fig:T0uv4} (labels have been added for convenience).
	
	By \cref{config:uvCodd,config:uvH0}, $G_{u, v}$ is locally bipartite and $H_{0}$-free. By \cref{lemma4lifting}, $\delta(G_{u, v}) > 1/2 \cdot \abs{G_{u, v}}$. Since $G_{u, v}$ does not contain an odd wheel, $u_{1}v_{1}$ is not an edge. Also $v_{1}, t$ is a 2-dense pair, $tu_{1}$ is an edge and $G_{u, v}$ is $H_{0}$-free, so, by \cref{lemma4dense}, $u_{1}, v_{1}$ is a 2-sparse pair in $G_{u, v}$. Similarly $u_{6}, v_{6}$ is 2-sparse in $G_{u, v}$.
	
	Now let $X = \set{u, u_{1}, u_{6}, v_{0}, \dotsc, v_{6}}$ (note that this does not contain $t$ or $v$). Within $G_{u, v}$, $u_{1}$ together with the $v_{i}$ form a 7-wheel missing the spoke $u_{1}v_{1}$ which is a sparse pair. Since $G_{u,v}$ is locally bipartite with $\delta(G_{u, v}) > 1/2 \cdot \abs{G_{u, v}}$, \cref{lemma45wheel} implies that $G_{u, v}$ does not contain any 5-wheels missing a sparse spoke. Hence, by \cref{lemma42nbs}, any neighbour of $u_{1}$ in $G_{u, v}$ is adjacent to at most two of the $v_{i}$. Thus, any neighbour of $u$, $v$, $u_{1}$ has at most five neighbours in $X$ (two amongst $v_{i}$ together with possibly $u_{1}, u_{6}, u$). Similarly any neighbour of $u$, $v$, $u_{6}$ has at most five neighbours in $X$. Next consider a vertex $x$ adjacent to both $u$, $v$ but to neither $u_{1}$ nor $u_{6}$. As $G_{u,v}$ is locally bipartite, $x$ is adjacent to at most six of the $v_{i}$ so $x$ has at most seven neighbours in $X$. Finally $\chi(G[X]) = 5$ so all vertices have at most nine neighbours in $X$. Hence
	\begin{align*}
		10 \delta(G) & \leqslant 5\abs{\Gamma(u, v, u_{1}) \cup \Gamma(u, v, u_{6})} + 7(\abs{G_{u, v}} - \abs{\Gamma(u, v, u_{1}) \cup \Gamma(u, v, u_{6})}) + 9(\abs{G} - \abs{G_{u, v}}) \\
		& = 9 \abs{G} - 2 \abs{G_{u, v}} - 2 \abs{\Gamma(u, v, u_{1}) \cup \Gamma(u, v, u_{6})} \leqslant 9 \abs{G} - 2 \abs{G_{u, v}} - 2 \abs{G_{u, v, u_{1}}} \\
		& \leqslant 9 \abs{G} - 2(2 \delta(G) - \abs{G}) - 2(3 \delta(G) - 2 \abs{G}) = 15 \abs{G} - 10 \delta(G),
	\end{align*}
	which contradicts $\delta(G) > 3/4 \cdot \abs{G}$.
	
	Now let $b \geqslant 5$ and suppose $G + uv$ does contain a $K_{b - 2} + T_{0}$ where both of $u$, $v$ are in the $(b - 2)$-clique. Let $L$ be the $(b - 2)$-clique without $u$, $v$: $L$ is a $(b - 4)$-clique in $G$ and so $G_{L}$ is 5-colourable and so locally 4-partite. Also, by \cref{lemma4lifting}, $\delta(G_{L}) > 22/29 \cdot \abs{G_{L}}$ and, by \cref{remark4lifting}, the pair $u, v$ is 4-sparse in $G_{L}$. Finally $G_{L} + uv$ contains $uv + T_{0}$, which contradicts the result just proved for $b = 4$.
\end{Proof}

\subsection{Finishing the proof}\label{sec:proof4localbpart}

Here we will prove \cref{main4alocalbpart} for locally $b$-partite graphs.

\begin{Proof}
	Take an edge-maximal locally $b$-partite graph $G$ with $\delta(G) > (1 - 1/(b + 1/7)) \cdot \abs{G}$. We need to show that $G$ is $(b + 1)$-colourable. We may assume by induction that the theorem holds for all $b'$ with $3 \leqslant b' < b$ (if there are any).
	
	We first show that for any $b$-sparse pair $u, v$ of $G$, $G' = G + uv$ contains a $(b - 2)$-clique $K$ with $G'_{K}$ not 3-colourable. Indeed, for $b = 3$, $G'$ is not locally tripartite (by edge-maximality) so there is a vertex  in $G'$ whose neighbourhood is not 3-colourable. Take $K$ to be this vertex. For $b > 3$, $G'$ is not locally $b$-partite so contains a vertex $w_{1}$ with $G'_{w_{1}}$ not $b$-colourable. Applying \cref{lemma4lifting} with $X = \set{w_{1}}$ and $\gamma = 1/7$ gives
	\begin{equation*}
		\delta(G'_{w_{1}}) > \biggl(1 - \frac{1}{b - 1 + 1/7}\biggr) \cdot \lvert G'_{w_{1}} \rvert.
	\end{equation*}
	By the induction hypothesis, if $G'_{w_{1}}$ was locally $(b - 1)$-partite, then it would be $b$-colourable. In particular, $G'_{w_{1}}$ is not locally $(b - 1)$-partite and so there is a vertex $w_{2}$ in $G'_{w_{1}}$ with $G'_{w_{1}, w_{2}}$ not $(b - 1)$-colourable. Repeating this argument gives a $(b - 2)$-clique $K$ with $G'_{K}$ not 3-colourable.
	
	Now, applying \cref{lemma4lifting} with $X = K$ and $\gamma = 1/7$ gives
	\begin{equation*}
		\delta(G'_{K}) > \biggl(1 - \frac{1}{b - (b - 2) + 1/7}\biggr) \cdot \lvert G'_{K}\rvert = 8/15 \cdot \lvert G'_{K} \rvert.
	\end{equation*}
	By \cref{main4localbip} and \cref{lemma4H}, $G'_{K}$ contains either an odd wheel, a copy of $H_{2}$, or a copy of $T_{0}$. Hence $G'$ contains either $K_{b - 2} + W_{\textnormal{odd}} = K_{b - 1} + C_{\textnormal{odd}}$, $K_{b - 2} + H_{2}$ or $K_{b - 2} + T_{0}$. Note that $G$ cannot contain any of these so $uv$ is a missing edge from one of these configurations. \cref{config:uvCodd,config:uvH0,config:uvT0} mean that both $u$ and $v$ lie in the $C_{\textnormal{odd}}$, the $H_{2}$, or the $T_{0}$.
	
	In particular, $u, v \not \in K$ so $K$ is a $(b - 2)$-clique in $G$ and $V(G_{K}) = V(G'_{K})$. We have the following facts.
	\begin{itemize}[noitemsep]
		\item By \cref{remark4lifting}, $G_{K}$ is 3-colourable and so locally bipartite.
		\item By \cref{remark4lifting}, $u, v$ is a 2-sparse pair in $G_{K}$.
		\item Applying \cref{lemma4lifting} with $X = K$ and $\gamma = 1/7$ gives $\delta(G_{K}) > 8/15 \cdot \abs{G_{K}}$.
		\item The graph $G_{K}$ contains no odd wheel, $H_{0}$ or $T_{0}$ ($G_{K}$ is 3-colourable) but the addition of $uv$ introduces an odd wheel, a copy of $H_{2}$, or a copy of $T_{0}$.
	\end{itemize}
	Using the argument at the start of \cref{sec:proof4localbip}, we deduce that, within $G_{K}$, there must be one of the configurations appearing in \cref{fig:cases} (with labels $u$ and $v$ possibly swapped). Note in that proof we only used edge-maximality to show that $uv$ was the missing edge of an odd wheel, a copy of $H_{2}$ or a copy of $T_{0}$ (and so here we do not need $G_{K}$ to be an edge-maximal locally bipartite graph).
	
	We now mimic the remainder of the proof of \cref{main4localbip}. Let $I$ be a largest independent set in $G$: $\abs{I} = \alpha(G)$.
	
	\begin{proposition}\label{lemma4bI}
		For all distinct $u, v \in I$, the pair $u, v$ is $b$-dense and furthermore every $u \in I$ has $I = D_{u} \cup \set{u}$.
	\end{proposition}
	
	\begin{Proof}
		Fix distinct $u, v \in I$. We will first show that $G_{u, v}$ is not $(b - 1)$-colourable. Note that $\Gamma(u), \Gamma(v) \subset V(G) \setminus I$ so $\abs{\Gamma(u) \cup \Gamma(v)} \leqslant \abs{G} - \abs{I}$. Also $I \subset V(G) \setminus \Gamma(u)$, so $\abs{I} \leqslant \abs{G} - d(u) \leqslant \abs{G} - \delta(G)$. Hence,
		\begin{align*}
			\abs{\Gamma(u, v)} & = d(u) + d(v) - \abs{\Gamma(u) \cup \Gamma(v)} \geqslant 2 \delta(G) + \abs{I} - \abs{G} \\
			& = b \delta(G) - (b - 1) \abs{G} + (b - 2)(\abs{G} - \delta(G)) + \abs{I} \\
			& \geqslant b \delta(G) - (b - 1) \abs{G} + (b - 1) \abs{I} > (b - 1) \abs{I},
		\end{align*}
		where we used $\delta(G) > (1 - 1/b) \cdot \abs{G}$ in the final inequality. But $I$ was a largest independent set in $G$ so $G_{u, v}$ is not $(b - 1)$-colourable.
		
		Now we will show that $u, v$ is $b$-dense. Suppose not and so they form a $b$-sparse pair. If $b = 3$, then $G_{u, v}$ is not bipartite and so contains an odd cycle. This contradicts \cref{config:uvCodd}. For $b \geqslant 4$ we will find a $(b - 4)$-clique $K$ in $G_{u, v}$ with $G_{u, v, K}$ not 3-colourable. For $b = 4$, we take $K = \emptyset$ and this suffices. For larger $b$, we note that, by \cref{lemma4lifting}, $\delta(G_{u, v}) > (1 - 1/(b - 2 + 1/7)) \cdot \abs{G_{u, v}}$. By the induction hypothesis, if $G_{u, v}$ were locally $(b - 2)$-partite, then $G_{u, v}$ would be $(b - 1)$-colourable, which it is not. Thus, there is $w_{1} \in G_{u, v}$ with $G_{u, v, w_{1}}$ not $(b - 2)$-colourable. Repeating this argument we obtain a $(b - 4)$-clique $K$ with $G_{u, v, K}$ not 3-colourable. Applying \cref{lemma4lifting} with $X = \set{u, v} \cup K$ and $\gamma = 1/7$ gives
		\begin{equation*}
			\delta(G_{u, v, K}) > (1 - 1/(b - (b - 2) + 1/7)) \cdot \abs{G_{u, v, K}} = 8/15 \cdot \abs{G_{u, v, K}}.
		\end{equation*}
		Then \cref{main4localbip} gives that $G_{u, v, K}$ either contains an odd wheel, a copy of $H_{0}$ or a copy of $T_{0}$. These contradict \cref{config:uvCodd,config:uvH0,config:uvT0} (applied to $G$). We have shown that $u, v$ must be $b$-dense. Thus, $I \subset D_{u} \cup \set{u}$.
		
		On the other hand, by the definition of density and \cref{lemma4bdense}, $D_{u} \cup \set{u
		}$ is an independent set. It contains the maximal independent set $I$ so must equal it.
	\end{Proof}
	
	\begin{definition}[$b$-quasidense]
		A pair of vertices $u, v$ is \defn{$b$-quasidense} if there is a sequence of vertices $u = d_{1}, d_{2}, \dotsc, d_{k}, d_{k + 1} = v$ such that all pairs $d_{i}, d_{i + 1}$ are $b$-dense ($i = 1, 2, \dotsc, k$).
	\end{definition}
	
	\cref{lemma4bI} immediately implies that if $u, v$ is $b$-quasidense and $u \in I$, then $v \in I$ as well. Now we can finish the proof. It suffices to show that every vertex is either in $I$ or is adjacent to all of $I$. Indeed, we may then fix $u \in I$ and note that $G[V(G) \setminus I] = G_{u}$ so $G[V(G) \setminus I]$ is $b$-colourable. Using a further colour for the independent set $I$ gives a $(b + 1)$-colouring of $G$.
	
	Suppose instead there is $u \in I$ and $v \not \in I$ with $u$ not adjacent to $v$. In particular, the pair $u, v$ cannot be $b$-quasidense and so is $b$-sparse. Thus from our remarks just preceding \cref{lemma4bI}, there is a $(b - 2)$-clique $K$ in $G$ such that $G_{K}$ contains one of the configurations appearing in \cref{fig:cases} (with labels $u$ and $v$ possibly swapped) and the pair $u, v$ is 2-sparse in $G_{K}$.
	
	Focus on $G_{K}$ -- this graph is 3-colourable so locally bipartite, $H_{0}$-free and $T_{0}$-free and the pair $u, v$ is 2-sparse in $G_{K}$. Also, by \cref{lemma4lifting}, $\delta(G_{K}) > (1 - 1/(b - (b - 2) + 1/7)) \cdot \abs{G_{K}} = 8/15 \cdot \abs{G_{K}}$. In the proof of \cref{lemma4adjI}, we used these facts alone to show that $u, v$ is quasidense in every configuration appearing in \cref{fig:cases}. Hence the pair $u, v$ is quasidense in $G_{K}$ so is $b$-quasidense in $G$. This is our required contradiction.
\end{Proof}

\section{\texorpdfstring{$a$}{a}-locally \texorpdfstring{$b$}{b}-partite graphs}\label{sec:alocalbpart}

In this section we relate the chromatic profile of $a$-locally $b$-partite graphs to the chromatic profile of locally $b$-partite graphs, making precise our comment in the introduction that to understand $a$-locally $b$-partite graphs it seems to be enough to understand locally $b$-partite graphs. This is elucidated at the end of \cref{sec:alocalbpart:first} and just before \cref{spec4alocalbip}. Along the way we will prove \cref{main4alocalbpart,main4alocalbipart}.

\subsection{The first threshold -- proving \texorpdfstring{\cref{main4alocalbpart,main4alocalbipart}}{Theorems 1.3 and 1.4}}\label{sec:alocalbpart:first}

As noted in the introduction, the first interesting threshold is $\delta_{\chi}(\cF_{a,b}, a + b)$ -- what values of $c$ guarantee that any $a$-locally $b$-partite graph with $\delta(G) \geqslant c \abs{G}$ is $(a + b)$-colourable? We already know $\delta_{\chi}(\cF_{1,2}, 3) = 4/7$ and $\delta_{\chi}(\cF_{1,b}, b + 1) \leqslant 1 - 1/(b + 1/7)$, and will extend these to all values of $a$. To simplify the statements of our results and make comparisons between different values of $a$ and $b$, it is helpful to write
\begin{equation*}
	\delta_{\chi}(\cF_{a,b}, a + b) = 1 - \frac{1}{a + b - 1 + \gamma_{a, b}},
\end{equation*}
and to focus our attention on the $\gamma_{a, b}$. As $\delta_{\chi}(\cF_{a,b}, a + b) \geqslant \delta_{\chi}(\cF_{a,b})$ we have, from \cref{chromthreFab},
\begin{equation*}
	\gamma_{a, b} \geqslant 0.
\end{equation*}
We collect some other basic properties of the $\gamma_{a, b}$.

\begin{lemma}\label{lemma4gammaproperties}
	For all positive integers $a$ and $b$ the following hold.
	\begin{itemize}[noitemsep]
		\item $\delta_{\chi}(\cF_{a, b+1}, a + b) \leqslant \delta_{\chi}(\cF_{a + 1, b}, a + b)$ and so
		\begin{equation}\label{eq:ab}
			\gamma_{a, b + 1} \leqslant \gamma_{a + 1, b}.
		\end{equation}
		\item $1/(2 - \delta_{\chi}(\cF_{a, b}, a + b)) \leqslant \delta_{\chi}(\cF_{a + 1, b}, a + b + 1)$ and so
		\begin{equation}\label{eq:a}
			\gamma_{a, b} \leqslant \gamma_{a + 1, b}.
		\end{equation}
	\end{itemize}
	Also $\gamma_{1,2} = 1/3$ and $\gamma_{1, b} \leqslant 1/7$ for all $b \geqslant 3$.
\end{lemma}

\begin{Proof}
	In~\cite{Illingworth2022localbipart} it was shown that $\delta_{\chi}(\cF_{1,2},3) = 4/7$ giving $\gamma_{1,2} = 1/3$ while, for $b \geqslant 3$, \cref{sec:localbpart} showed $\delta_{\chi}(\cF_{1,b},b + 1) \leqslant 1 - 1/(b + 1/7)$ so $\gamma_{1, b} \leqslant 1/7$.
	
	Now $\cF_{a, b + 1} \subset \cF_{a + 1, b}$ from which $\delta_{\chi}(\cF_{a, b + 1}, a + b) \leqslant \delta_{\chi}(\cF_{a + 1, b}, a + b)$ immediately follows. This gives inequality \eqref{eq:ab}.
	
	Finally let $d < \delta_{\chi}(\cF_{a, b}, a + b)$: there is an $a$-locally $b$-partite graph $G$ with $\delta(G) \geqslant d \abs{G}$ and $\chi(G) > a + b$. Let $G'$ be $G$ joined to an independent set of size $\abs{G} - \delta(G)$, that is,
	\begin{equation*}
		G' = K_{1}(\abs{G} - \delta(G)) + G.
	\end{equation*}
	Since $G$ is $a$-locally $b$-partite, it is also $(a + 1)$-locally $(b - 1)$-partite. From both of these it follows that $G'$ is $(a + 1)$-locally $b$-partite. Also, $\chi(G') = \chi(G) + 1 > a + b + 1$, and
	\begin{align*}
		\frac{\delta(G')}{\abs{G'}} = \frac{\abs{G}}{2 \abs{G} - \delta(G)} = \frac{1}{2 - \delta(G) \abs{G}^{-1}} \geqslant \frac{1}{2 - d},
	\end{align*}
	so $\delta_{\chi}(\cF_{a + 1, b}, a + b + 1) \geqslant 1/(2 - d)$. This holds for all $d < \delta_{\chi}(\cF_{a, b}, a + b)$, so $\delta_{\chi}(\cF_{a + 1, b}, a + b + 1) \geqslant 1/(2 - \delta_{\chi}(\cF_{a, b}, a + b))$. Thus
	\begin{equation*}
		1 - \frac{1}{a + b + \gamma_{a + 1, b}} \geqslant \frac{1}{1 + \frac{1}{a + b - 1 + \gamma_{a, b}}} = \frac{a + b - 1 + \gamma_{a, b}}{a + b + \gamma_{a, b}} = 1 - \frac{1}{a + b + \gamma_{a, b}},
	\end{equation*}
	and so $\gamma_{a + 1, b} \geqslant \gamma_{a, b}$, as required.
\end{Proof}

Inequality \eqref{eq:a} gives a lower bound for $\gamma_{a + 1, b}$ in terms of $\gamma_{a, b}$. The next lemma, which lies at the heart of our analysis, gives an upper bound.

\begin{lemma}\label{lemma4gammamain0}
	For all positive integers $a$ and $b$,
	\begin{equation}\label{eq:gammamain0}
		\gamma_{a, b} \leqslant \gamma_{a + 1, b} \leqslant \max \set{\gamma_{a, b}, \gamma_{1, a + b}}.
	\end{equation}
\end{lemma}

\begin{Proof}
	The left-hand inequality is just inequality \eqref{eq:a}. Let $\gamma = \max \set{\gamma_{a, b}, \gamma_{1, a + b}}$. Let $G$ be an $(a + 1)$-locally $b$-partite graph with
	\begin{equation*}
		\delta(G) > \biggl(1 - \frac{1}{a + b + \gamma}\biggr) \cdot \abs{G}.
	\end{equation*}
	It suffices to show that $\chi(G) \leqslant a + b + 1$ as then 
	\begin{equation*}
		1 - \frac{1}{a + b + \gamma} \geqslant \delta_{\chi}(\cF_{a + 1, b}, a + b + 1) = 1 - \frac{1}{a + b + \gamma_{a + 1, b}}.
	\end{equation*}
	Fix any $u \in V(G)$ and consider $G_{u}$: $G_{u}$ is an $a$-locally $b$-partite graph with
	\begin{equation*}
		\delta(G_{u}) > \biggl(1 - \frac{1}{a + b - 1 + \gamma}\biggr) \cdot \abs{G_{u}},
	\end{equation*}
	by the lifting lemma, \cref{lemma4lifting}. But $\gamma \geqslant \gamma_{a, b}$ so
	\begin{equation*}
		\delta(G_{u}) > \delta_{\chi}(\cF_{a, b}, a + b) \cdot \abs{G_{u}},
	\end{equation*}
	and hence $G_{u}$ is $(a + b)$-colourable. Thus the graph $G$ is locally $(a + b)$-partite. Also $\gamma \geqslant \gamma_{1, a + b}$, so
	\begin{equation*}
		\delta(G) > \delta_{\chi}(\cF_{1, a + b}, a + b + 1) \cdot \abs{G}.
	\end{equation*}
	Thus $G$ is $(a + b + 1)$-colourable.
\end{Proof}

From this one can immediately deduce \cref{main4alocalbipart,main4alocalbpart}.

\begin{corollary}[\cref{main4alocalbipart,main4alocalbpart}]\label{cor4gamma0}
	For all positive integers $a$ and for all $b \geqslant 3$,
	\begin{equation*}
		\gamma_{a, 2} = 1/3, \quad \gamma_{a, b} \leqslant 1/7,
	\end{equation*}
	and so
	\begin{equation*}
		\delta_{\chi}(\cF_{a, 2}, a + 2) = 1 - \frac{1}{a + 1 + 1/3}, \qquad \delta_{\chi}(\cF_{a, b}, a + b) \leqslant 1 - \frac{1}{a + b - 1 + 1/7}.
	\end{equation*}
\end{corollary}

\begin{Proof}
	From \cref{lemma4gammaproperties}, $\gamma_{1, 2} = 1/3$ and $\gamma_{1, b} \leqslant 1/7$ for any $b \geqslant 3$. By \cref{lemma4gammamain0}, for any $a$ and any $b \geqslant 2$,
	\begin{equation*}
		\gamma_{a, b} \leqslant \gamma_{a + 1, b} \leqslant \max \set{\gamma_{a, b}, \gamma_{1, a + b}} \leqslant \max \set{\gamma_{a, b}, 1/7}.
	\end{equation*}
	An easy induction gives $\gamma_{a, 2} = 1/3$ for all $a$ and $\gamma_{a, b} \leqslant 1/7$ for any $b \geqslant 3$.
\end{Proof}

Note that inequalities \eqref{eq:ab} and \eqref{eq:a} give
\begin{equation*}
	\gamma_{a, b} \geqslant \gamma_{1, b'},
\end{equation*}
for all $b \leqslant b' \leqslant a + b - 1$. In line with our inductive arguments in \cref{sec:localbpart}, we believe that in fact $\gamma_{1, b} \geqslant \gamma_{1, b + 1}$ for all $b$. If this were true, then $\gamma_{a, b} \geqslant \gamma_{1, a + b}$, and so \cref{lemma4gammamain0} would give $\gamma_{a, b} = \gamma_{a + 1, b}$ for all $a, b$. Of course, this implies $\gamma_{a, b} = \gamma_{1, b}$ and so $\delta_{\chi}(\cF_{a, b}, a + b)$ would be determined by $\delta_{\chi}(\cF_{1, b}, b + 1)$ -- a particular manifestation of our aforementioned belief that to understand $a$-locally $b$-partite graphs, we should focus on locally $b$-partite graphs. It also highlights the following question.

\begin{question}
	Is the sequence $\gamma_{1, b}$ non-increasing in $b$?
\end{question}

\subsection{\texorpdfstring{$a$}{a}-locally bipartite graphs}\label{sec:alocalbpart:second}

One could replicate the elementary approach of the previous section to try to evaluate $\delta_{\chi}(\cF_{a, b}, k)$ for $k > a + b$. Indeed one might define $\gamma_{a, b, m}$ by
\begin{equation*}
	\delta_{\chi}(\cF_{a, b}, a + b + m) = 1 - \frac{1}{a + b - 1 + \gamma_{a, b, m}},
\end{equation*}
so that $\gamma_{a, b, 0} = \gamma_{a, b}$. Many of the properties of the $\gamma_{a, b}$ pass over: the $\gamma_{a, b, m}$ are non-negative (and, in fact, $\lim_{m \to \infty} \gamma_{a, b, m} = 0$) and both inequalities \eqref{eq:ab} and \eqref{eq:a} extend easily ($\gamma_{a, b + 1, m} \leqslant \gamma_{a + 1, b, m}$ and $\gamma_{a, b, m} \leqslant \gamma_{a + 1, b, m}$). However, there seems to be no argument to produce inequality \eqref{eq:gammamain0} or anything similar. A more involved approach would be required.

The next threshold to consider is $\delta_{\chi}(\cF_{a, b}, a + b + 1)$. For locally bipartite graphs, we showed $\delta_{\chi}(\cF_{1, 2}, 4) \allowbreak \leqslant 6/11$ and had many structural results (some of which we will extend). For $b \geqslant 3$, we know very little about $\delta_{\chi}(\cF_{1, b}, b + 2)$ beyond it being at least $\delta_{\chi}(\cF_{1, b}) = 1 - 1/b$ and at most $\delta_{\chi}(\cF_{1, b}, b + 1) \leqslant 1 - 1/(b + 1/7)$. The question for $b = 3$ is of particular interest. Tantalisingly, the complement of the 9-cycle is locally tripartite, 5-chromatic and has minimum degree $6 = 2/3 \cdot 9$.

\begin{question}
	Is there a locally tripartite graph $G$ with minimum degree greater than $2/3 \cdot \abs{G}$ which is not 4-colourable?
\end{question}

We now focus on $a$-locally bipartite graphs. The following theorem, which will be essential for extending the Andr\'{a}sfai-Erd\H{o}s-S\'{o}s theorem~\cite{Illingworth2023mindegstab}, should be compared to \cref{spec4localbip} -- again we see that the key to understanding $a$-locally bipartite graphs is to understand locally bipartite ones. The proof is an induction combining our results for locally bipartite (\cref{spec4localbip}) and locally $b$-partite (\cref{main4alocalbpart}) graphs.

\begin{theorem}[$a$-locally bipartite graphs]\label{spec4alocalbip}
	Let $G$ be an $a$-locally bipartite graph.
	\begin{enumerate}[noitemsep, label=\textit{\alph{*}}., ref=\textit{\alph{*}}]
		\item If $\delta(G) > (1 - 1/(a + 4/3)) \cdot \abs{G}$, then $G$ is $(a + 2)$-colourable. Suitable blow-ups of $K_{a - 1} + \overline{C}_{7}$ show that this is tight.
		\item If $\delta(G) > (1 - 1/(a + 5/4)) \cdot \abs{G}$, then $G$ is either $(a + 2)$-colourable or contains $K_{a - 1} + \overline{C}_{7}$.
		\item If $\delta(G) > (1 - 1/(a + 6/5)) \cdot \abs{G}$, then $G$ is either $(a + 2)$-colourable or contains $K_{a - 1} + \overline{C}_{7}$ or $K_{a - 1} + H_{2}^{+}$.
		\item If $\delta(G) > (1 - 1/(a + 7/6)) \cdot \abs{G}$, then $G$ is either $(a + 2)$-colourable or contains $K_{a - 1} + H_{2}$.
		\item If $\delta(G) > (1 - 1/(a + 8/7)) \cdot \abs{G}$, then $G$ is either $(a + 2)$-colourable or contains $K_{a - 1} + H_{2}$ or $K_{a - 1} + T_{0}$. \label{e2}
	\end{enumerate}
\end{theorem}

\begin{Proof}
	The graph $K_{a - 1} + \overline{C}_{7}$ is $a$-locally bipartite with chromatic number $a + 3$. Hence $K_{a - 1}(3) + \overline{C}_{7}$ also has these properties and is, furthermore, $(3a + 1)$-regular with $3a + 4$ vertices. Thus balanced blow-ups of $K_{a - 1}(3) + \overline{C}_{7}$ give the tightness of the first bullet point. Proving everything else is a simple induction on $a$ (with \cref{spec4localbip} covering the base case). Indeed we will just demonstrate it for part~\ref{e2}. Let $G$ be an $a$-locally bipartite graph with $\delta(G) > (1 - 1/(a + 8/7)) \cdot \abs{G}$. Fix any vertex $u$ of $G$ and consider $G_{u}$: this is $(a - 1)$-locally bipartite and by \cref{lemma4lifting}, $\delta(G_{u}) > (1 - 1/(a - 1 + 8/7)) \cdot \abs{G_{u}}$. By induction, either $G_{u}$ contains one of $K_{a - 2} + H_{2}$, $K_{a - 2} + T_{0}$ or is $(a + 1)$-colourable. If there is some vertex $u$ with $G_{u}$ not $(a + 1)$-colourable, then $G$ contains one of $K_{a - 1} + H_{2}$, $K_{a - 1} + T_{0}$. Otherwise, $G$ is locally $(a + 1)$-partite. But, by \cref{main4alocalbpart},
	\begin{equation*}
		\delta_{\chi}(\cF_{1, a + 1}, a + 2) \leqslant 1 - \frac{1}{a + 8/7},
	\end{equation*}
	so $G$ is $(a + 2)$-colourable, as required.
\end{Proof}

\section*{Acknowledgement}

It is a pleasure to thank Andrew Thomason for many helpful discussions. I am grateful to the anonymous referees for their careful reading and excellent suggestions for improving the presentation.

{
\fontsize{11pt}{12pt}
\selectfont
	
\hypersetup{linkcolor={red!70!black}}
\setlength{\parskip}{2pt plus 0.3ex minus 0.3ex}

\newcommand{\etalchar}[1]{$^{#1}$}

}

\begin{appendix}

\section{Verifying \texorpdfstring{\cref{fig:containment}}{Figure 3}}\label{sec:containment}

In this appendix we verify \cref{fig:containment} which, for convenience, we display again here.

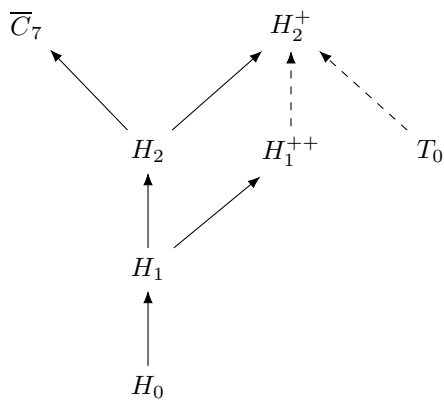
\begin{figure}[H]
	\centering
	\begin{tikzpicture}
		\node(H21) at (0,0){$H_{2}^{+}$};
		\node(blank) at (-1.7,0){\color{white}H1};
		\node(C7) at (-3.4,0){$\overline{C}_{7}$};
		\node(H12) at (0,-1.7){$H_{1}^{++}$};
		\node(T0) at (1.7,-1.7){$T_{0}$};
		\node(H2) at (-1.7,-1.7){$H_{2}$};
		\node(H1) at (-1.7,-3.4){$H_{1}$};
		\node(H0) at (-1.7,-5.1){$H_{0}$};

		\draw[-Latex] (H2) -- (H21);
		\draw[-Latex] (H2) -- (C7);
		\draw[dashed, -Latex] (H12) -- (H21);
		\draw[dashed, -Latex] (T0) -- (H21);
		\draw[-Latex] (H1) -- (H2);
		\draw[-Latex] (H1) -- (H12);
		\draw[-Latex] (H0) -- (H1);
	\end{tikzpicture}
\end{figure}

The reader will recall that full arrows represent containment and dashed arrows represent homomorphisms. Furthermore, $H$ is homomorphic to $G$ in the diagram if there is a sequence of arrows starting at $H$ and ending at $G$.

All the containments are clear. The following figure gives a homomorphism from $H_{1}^{++}$ to $H_{2}^{+}$: the left diagram is a labelling of the vertices of $H_{1}^{++}$ and the right diagram shows the images of those vertices under the map.

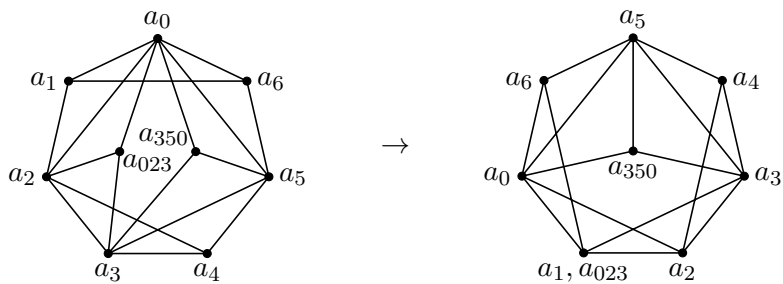
\begin{figure}[H]
	\centering
	\begin{subfigure}{.33\textwidth}
		\centering
		\begin{tikzpicture}
			\foreach \pt in {0,1,...,6} 
			{
				\tkzDefPoint(\pt*360/7 + 90:1.5){v_\pt}
			}
			\tkzDefPoint(-0.5,0){ul}
			\tkzDefPoint(0.5,0){ur}
			
			\tkzDrawPolySeg(v_0,v_1,v_2,v_3,v_4,v_5,v_6,v_0) 
			\tkzDrawPolySeg(v_3,v_5,v_0,v_2,v_4)
			\tkzDrawPolySeg(v_6,v_1)
			\tkzDrawSegments(ul,v_0 ul,v_2 ul,v_3 ur,v_0 ur,v_5 ur,v_3)
			\tkzDrawPoints(v_0,v_...,v_6)
			\tkzDrawPoints(ul,ur)
			\tkzLabelPoint[above](v_0){$a_{0}$}
			\tkzLabelPoint[left](v_1){$a_{1}$}
			\tkzLabelPoint[left](v_2){$a_{2}$}
			\tkzLabelPoint[below](v_3){$a_{3}$}
			\tkzLabelPoint[below](v_4){$a_{4}$}
			\tkzLabelPoint[right](v_5){$a_{5}$}
			\tkzLabelPoint[right](v_6){$a_{6}$}
			\tkzLabelPoint[below right = -3pt](ul){$a_{023}$}
			\tkzLabelPoint[above left = -2pt](ur){$a_{350}$}
		\end{tikzpicture}
	\end{subfigure}
	\begin{subfigure}{.05\textwidth}
		\centering
		\begin{tikzpicture}
			\tikz\draw[thick,black,->] (0,-1.7) -- (0.9\textwidth,-1.7);
		\end{tikzpicture}
	\end{subfigure}
	\begin{subfigure}{.33\textwidth}
		\centering
		\begin{tikzpicture}
			\foreach \pt in {0,1,...,6} 
			{
				\tkzDefPoint(\pt*360/7 + 90:1.5){v_\pt}
			} 
			\tkzDefPoint(0,0){u}
			\tkzDrawPolySeg(v_0,v_1,v_2,v_3,v_4,v_5,v_6,v_0)
			\tkzDrawPolySeg(v_1,v_3,v_5,v_0, v_2,v_4,v_6)
			\tkzDrawSegments(u,v_0 u,v_2 u,v_5)
			\tkzDrawPoints(v_0,v_...,v_6)
			\tkzDrawPoint(u)
			\tkzLabelPoint[above](v_0){$a_{5}$}
			\tkzLabelPoint[left](v_1){$a_{6}$}
			\tkzLabelPoint[left](v_2){$a_{0}$}
			\tkzLabelPoint[below](v_3){$a_{1}, a_{023}$}
			\tkzLabelPoint[below](v_4){$a_{2}$}
			\tkzLabelPoint[right](v_5){$a_{3}$}
			\tkzLabelPoint[right](v_6){$a_{4}$}
			\tkzLabelPoint[below](u){$a_{350}$}
		\end{tikzpicture}
	\end{subfigure}
\end{figure}
\addtocounter{figure}{-1}

The following figures gives a homomorphism from $T_{0}$ to $H_{2}^{+}$: the left diagram is a labelling of the vertices of $T_{0}$ and the right diagram shows the images of those vertices under the map.

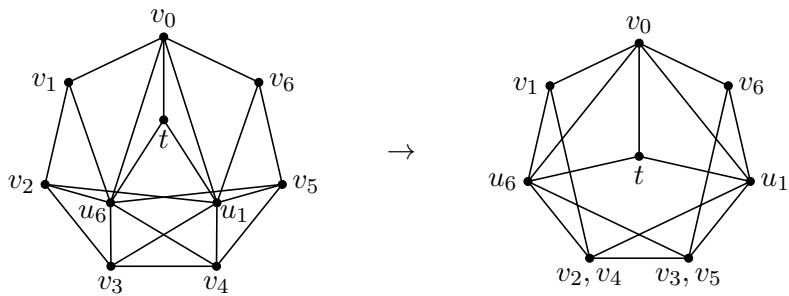
\begin{figure}[H]
	\centering
	\begin{subfigure}{.33\textwidth}
		\centering
		\begin{tikzpicture}
			\foreach \pt in {0,1,...,6} 
			{
				\tkzDefPoint(\pt*360/7 + 90:1.6){v_\pt}
			}
			\tkzDefPoint(0,0.5){t}
			\tkzDefPoint(0.7,-0.6){u_1}
			\tkzDefPoint(-0.7,-0.6){u_6}
			\tkzDrawPolySeg(v_0,v_1,v_2,v_3,v_4,v_5,v_6,v_0) 
			\tkzDrawSegments(t,v_0 t,u_1 t,u_6)
			\foreach \pt in {0,2,3,4,5,6}
			{
				\tkzDrawSegment(u_1,v_\pt)
			}
			\foreach \pt in {0,1,2,3,4,5}
			{
				\tkzDrawSegment(u_6,v_\pt)
			}
			\tkzDrawPoints(v_0,v_...,v_6)
			\tkzDrawPoints(t,u_1,u_6)
			\tkzLabelPoint[above](v_0){$v_{0}$}
			\tkzLabelPoint[left](v_1){$v_{1}$}
			\tkzLabelPoint[left](v_2){$v_{2}$}
			\tkzLabelPoint[below](v_3){$v_{3}$}
			\tkzLabelPoint[below](v_4){$v_{4}$}
			\tkzLabelPoint[right](v_5){$v_{5}$}
			\tkzLabelPoint[right](v_6){$v_{6}$}
			\tkzLabelPoint[below](t){$t$}
			\tkzLabelPoint[below right = -3pt](u_1){$u_{1}$}
			\tkzLabelPoint[below left = -3pt](u_6){$u_{6}$}
		\end{tikzpicture}
	\end{subfigure}
	\begin{subfigure}{.05\textwidth}
		\centering
		\begin{tikzpicture}
			\tikz\draw[thick,black,->] (0,-1.7) -- (0.9\textwidth,-1.7);
		\end{tikzpicture}
	\end{subfigure}
	\begin{subfigure}{.33\textwidth}
		\centering
		\begin{tikzpicture}
			\foreach \pt in {0,1,...,6} 
			{
				\tkzDefPoint(\pt*360/7 + 90:1.5){v_\pt}
			} 
			\tkzDefPoint(0,0){u}
			\tkzDrawPolySeg(v_0,v_1,v_2,v_3,v_4,v_5,v_6,v_0)
			\tkzDrawPolySeg(v_1,v_3,v_5,v_0, v_2,v_4,v_6)
			\tkzDrawSegments(u,v_0 u,v_2 u,v_5)
			\tkzDrawPoints(v_0,v_...,v_6)
			\tkzDrawPoint(u)
			\tkzLabelPoint[above](v_0){$v_{0}$}
			\tkzLabelPoint[left](v_1){$v_{1}$}
			\tkzLabelPoint[left](v_2){$u_{6}$}
			\tkzLabelPoint[below](v_3){$v_{2}, v_{4}$}
			\tkzLabelPoint[below](v_4){$v_{3}, v_{5}$}
			\tkzLabelPoint[right](v_5){$u_{1}$}
			\tkzLabelPoint[right](v_6){$v_{6}$}
			\tkzLabelPoint[below](u){$t$}
		\end{tikzpicture}
	\end{subfigure}
\end{figure}
\addtocounter{figure}{-1}

In particular, all arrows in \cref{fig:containment} are correct. We need to show that further arrows could not be added. There are some subtleties in our notation that we now elucidate. Given a homomorphism $\varphi \colon H \to G$, we say $\varphi$ is surjective or injective if the map $\varphi \colon V(H) \to V(G)$ is surjective or injective, respectively. Note that $\varphi$ being injective implies that $H$ is actually a subgraph of $G$. By $\varphi(H)$ we mean the graph on vertex set $\varphi(V(H))$ and edge set $\varphi(E(H))$. In particular, this is a spanning subgraph of $G[\varphi(V(H))]$ but it may not have all the edges of $G[\varphi(V(H))]$. We make frequent use of the fact that $\chi(\varphi(H)) \geqslant \chi(H)$.

We first deal with left-hand side ($H_{0} \to H_{1} \to H_{2} \to \overline{C}_{7}$) of the figure: we need to show that $\overline{C}_{7} \nrightarrow H_{2}$, $H_{2} \nrightarrow H_{1}$, and $H_{1} \nrightarrow H_{0}$. The arguments are very similar, making use of the fact that $H_{0}$, $H_{1}$, $H_{2}$, and $\overline{C}_{7}$ are all vertex-critical 4-chromatic graphs on seven vertices, so we only give the explicit proof for $H_{2} \nrightarrow H_{1}$.

\begin{proposition}
	The graph $H_{2}$ is not homomorphic to $H_{1}$.
\end{proposition}

\begin{Proof}
	Suppose there is a homomorphism $\varphi \colon H_{2} \to H_{1}$. Then $\chi(\varphi(H_{2})) \geqslant \chi(H_{2}) = 4$. Now $H_{1}$ is a vertex-critical 4-chromatic graph, so $\varphi$ is surjective. But $H_{1}$ and $H_{2}$ both have seven vertices so $\varphi$ is injective. That is, $H_{1}$ must contain a copy of $H_{2}$, which is absurd as $e(H_{1}) < e(H_{2})$.
\end{Proof}

We now know that the left-hand side of \cref{fig:containment} is correct and consider how $H_{2}^{+}$ relates to it. It suffices to show that $H_{2}^{+} \nrightarrow \overline{C}_{7}$ and $\overline{C}_{7} \nrightarrow H_{2}^{+}$ (note that $H_{2}^{+} \nrightarrow \overline{C}_{7}$ implies $H_{2}^{+} \nrightarrow H_{2}, H_{1}, H_{0}$). That $H_{2}^{+}$ is not homomorphic to $\overline{C}_{7}$ and vice versa follows from the next lemma, which appeared in~\cite{Illingworth2022localbipart} -- both $\overline{C}_{7}$ and $H_{2}^{+}$ are edge-maximal locally bipartite graphs and neither is a subgraph of the other ($\overline{C}_{7}$ has fewer vertices than $H_{2}^{+}$ and $H_{2}^{+}$ does not have seven vertices all of degree at least four).

\begin{lemma}\label{homscores}
	Let $F$ be an edge-maximal locally bipartite graph in which no two neighbourhoods are the same. Let $F$ be homomorphic to a locally bipartite graph $G$. Then $F$ is an induced subgraph of $G$.
\end{lemma}

Next we relate $H_{1}^{++}$ to the diagram. It suffices to show that $H_{1}^{++} \nrightarrow \overline{C}_{7}$ and $H_{2} \nrightarrow H_{1}^{++}$ (note that $H_{1}^{++} \nrightarrow \overline{C}_{7}$ implies $H_{1}^{++} \nrightarrow H_{2}, H_{1}, H_{0}$ while $H_{2} \nrightarrow H_{1}^{++}$ implies $\overline{C}_{7} \nrightarrow H_{1}^{++}$ and $H_{2}^{+} \nrightarrow H_{1}^{++}$).

\begin{proposition}\label{H12homC72}
	The graph $H_{1}^{++}$ is not homomorphic to $\overline{C}_{7}$.
\end{proposition}

\begin{Proof}
	Suppose there is a homomorphism $\varphi \colon H_{1}^{++} \to \overline{C}_{7}$. Label the copy of $H_{1}^{++}$ as shown below and let $A = \set{a_{0}, a_{1}, \dotsc, a_{6}}$ so $H_{1}^{++}[A]$ is a copy of $H_{1}$. Note that $\chi(\varphi(H_{1}^{++}[A])) \geqslant \chi(H_{1}^{++}[A]) = 4$ and $\overline{C}_{7}$ is a vertex-critical 4-chromatic graph, so the restriction of $\varphi$ to $A$ is a surjection onto $\overline{C}_{7}$ and so $\varphi(a_{0}), \varphi(a_{1}), \dotsc, \varphi(a_{6})$ are all distinct.
	
	
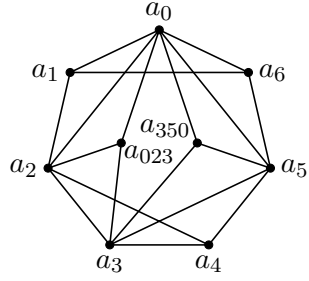
\begin{figure}[H]
		\centering
		\begin{tikzpicture}
			\foreach \pt in {0,1,...,6} 
			{
				\tkzDefPoint(\pt*360/7 + 90:1.5){v_\pt}
			}
			\tkzDefPoint(-0.5,0){ul}
			\tkzDefPoint(0.5,0){ur}
			
			\tkzDrawPolySeg(v_0,v_1,v_2,v_3,v_4,v_5,v_6,v_0) 
			\tkzDrawPolySeg(v_3,v_5,v_0,v_2,v_4)
			\tkzDrawPolySeg(v_6,v_1)
			\tkzDrawSegments(ul,v_0 ul,v_2 ul,v_3 ur,v_0 ur,v_5 ur,v_3)
			\tkzDrawPoints(v_0,v_...,v_6)
			\tkzDrawPoints(ul,ur)
			\tkzLabelPoint[above](v_0){$a_{0}$}
			\tkzLabelPoint[left](v_1){$a_{1}$}
			\tkzLabelPoint[left](v_2){$a_{2}$}
			\tkzLabelPoint[below](v_3){$a_{3}$}
			\tkzLabelPoint[below](v_4){$a_{4}$}
			\tkzLabelPoint[right](v_5){$a_{5}$}
			\tkzLabelPoint[right](v_6){$a_{6}$}
			\tkzLabelPoint[below right = -3pt](ul){$a_{023}$}
			\tkzLabelPoint[above left = -2pt](ur){$a_{350}$}
		\end{tikzpicture}
	\end{figure}
	
	Now $a_{0}$ has degree 6 while $\varphi(a_{0}) \in \overline{C}_{7}$ only has degree 4. The four neighbours of $\varphi(a_{0})$ are $\varphi(a_{1})$, $\varphi(a_{2})$, $\varphi(a_{5})$, $\varphi(a_{6})$ and so $\varphi(a_{023})$ is one of these. Also $a_{2}$ has degree 5 while $\varphi(a_{2})$ only has degree 4. The four neighbours of $\varphi(a_{2})$ are $\varphi(a_{0}), \varphi(a_{1}), \varphi(a_{3}), \varphi(a_{4})$ and so $\varphi(a_{023})$ is one of these. Hence $\varphi(a_{023}) = \varphi(a_{1})$. Similarly, considering the neighbourhoods of $a_{0}$ and $a_{5}$ shows that $\varphi(a_{350}) = \varphi(a_{6})$. Then $\Gamma(\varphi(a_{3}))$ contains
	\begin{equation*}
		\set{\varphi(a_{2}), \varphi(a_{023}), \varphi(a_{350}), \varphi(a_{4}), \varphi(a_{5})} = \set{\varphi(a_{2}), \varphi(a_{1}), \varphi(a_{6}), \varphi(a_{4}), \varphi(a_{5})},
	\end{equation*}
	which has size 5. This contradicts the 4-regularity of $\overline{C}_{7}$.
\end{Proof}

\begin{proposition}
	The graph $H_{2}$ is not homomorphic to $H_{1}^{++}$.
\end{proposition}

\begin{Proof}
	Suppose there is a homomorphism $\varphi \colon H_{2} \to H_{1}^{++}$. Now $\chi(\varphi(H_{2})) \geqslant \chi(H_{2}) = 4$ and any 6-vertex subgraph of $H_{1}^{++}$ is 3-colourable (it is homomorphic to some 6-vertex subgraph of $H_{2}^{+}$) so $\varphi$ must be injective. Thus $H_{1}^{++}$ contains $H_{2}$. But $H_{1}^{++}$ only has 4 vertices of degree at least 4 while $H_{2}$ has 5 vertices of degree 4.
\end{Proof}

Finally we relate $T_{0}$ to the diagram. It suffices to show that $H_{0} \nrightarrow T_{0}$, $T_{0} \nrightarrow \overline{C}_{7}$, and $T_{0} \nrightarrow H_{1}^{++}$ (note that $H_{0} \nrightarrow T_{0}$ implies that no other graph in the diagram is homomorphic to $T_{0}$ while $T_{0} \nrightarrow \overline{C}_{7}$ implies that $T_{0} \nrightarrow H_{0}, H_{1}, H_{2}$). We use the following labelling of the copy of $T_{0}$ in all three proofs.

\begin{figure}[H]
	\centering
	\begin{tikzpicture}
		\foreach \pt in {0,1,...,6} 
		{
			\tkzDefPoint(\pt*360/7 + 90:1.6){v_\pt}
		}
		\tkzDefPoint(0,0.5){t}
		\tkzDefPoint(0.7,-0.6){u_1}
		\tkzDefPoint(-0.7,-0.6){u_6}
		\tkzDrawPolySeg(v_0,v_1,v_2,v_3,v_4,v_5,v_6,v_0) 
		\tkzDrawSegments(t,v_0 t,u_1 t,u_6)
		\foreach \pt in {0,2,3,4,5,6}
		{
			\tkzDrawSegment(u_1,v_\pt)
		}
		\foreach \pt in {0,1,2,3,4,5}
		{
			\tkzDrawSegment(u_6,v_\pt)
		}
		\tkzDrawPoints(v_0,v_...,v_6)
		\tkzDrawPoints(t,u_1,u_6)
		\tkzLabelPoint[above](v_0){$v_{0}$}
		\tkzLabelPoint[left](v_1){$v_{1}$}
		\tkzLabelPoint[left](v_2){$v_{2}$}
		\tkzLabelPoint[below](v_3){$v_{3}$}
		\tkzLabelPoint[below](v_4){$v_{4}$}
		\tkzLabelPoint[right](v_5){$v_{5}$}
		\tkzLabelPoint[right](v_6){$v_{6}$}
		\tkzLabelPoint[below](t){$t$}
		\tkzLabelPoint[below right = -3pt](u_1){$u_{1}$}
		\tkzLabelPoint[below left = -3pt](u_6){$u_{6}$}
	\end{tikzpicture}
\end{figure}

\begin{proposition}
	The graph $H_{0}$ is not homomorphic to $T_{0}$.
\end{proposition}

\begin{Proof}
	We first claim that any 7-vertex subgraph of $T_{0}$ is 3-colourable. Let $F$ be a 7-vertex subgraph of $T_{0}$. If $F$ contains all the $v_{i}$, then $F$ is a subgraph of a 7-cycle and so is 3-colourable. Otherwise $F$ is a subgraph of $T_{0} - v_{i}$ for some $i$. This graph is 3-colourable: 2-colour the remaining $v_{j}$ with colours 1 and 2, give $u_{1}$ and $u_{6}$ colour 3 and then give $t$ colour 1 or 2 (opposite to the colour of $v_{0}$ if it is present).
	
	Suppose there is a homomorphism $\varphi \colon H_{0} \to T_{0}$. Then $\varphi(H_{0})$ is a subgraph of $T_{0}$ with at most 7 vertices, so is 3-colourable. But then $3 \geqslant \chi(\varphi(H_{0})) \geqslant \chi(H_{0}) = 4$.
\end{Proof}

\begin{proposition}
	The graph $T_{0}$ is not homomorphic to $\overline{C}_{7}$.
\end{proposition}

\begin{Proof}
	Suppose $\varphi \colon T_{0} \to \overline{C}_{7}$ is a homomorphism. Label the copy of $\overline{C}_{7}$ as follows.
	
	
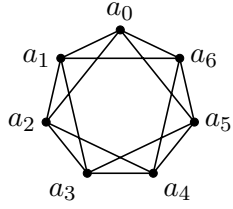
\begin{figure}[H]
		\centering
		\begin{tikzpicture}
			\foreach \pt in {0,1,...,6} 
			{
				\tkzDefPoint(\pt*360/7 + 90:1){v_\pt}
			} 
			\tkzDrawPolySeg(v_0,v_1,v_2,v_3,v_4,v_5,v_6,v_0) 
			\tkzDrawPolySeg(v_0,v_2,v_4,v_6,v_1,v_3,v_5,v_0)
			\tkzDrawPoints(v_0,v_...,v_6)
			\tkzLabelPoint[above](v_0){$a_{0}$}
			\tkzLabelPoint[left](v_1){$a_{1}$}
			\tkzLabelPoint[left](v_2){$a_{2}$}
			\tkzLabelPoint[below left](v_3){$a_{3}$}
			\tkzLabelPoint[below right](v_4){$a_{4}$}
			\tkzLabelPoint[right](v_5){$a_{5}$}
			\tkzLabelPoint[right](v_6){$a_{6}$}
		\end{tikzpicture}
	\end{figure}
	
	Without loss of generality we may assume $\varphi(u_{1}) = a_{0}$. The common neighbourhood $\Gamma(u_{1}, u_{6})$ contains the edge $tv_{0}$, so $\Gamma(\varphi(u_{1}), \varphi(u_{6}))$ contains an edge and so $\varphi(u_{6}) \in \set{a_{0}, a_{3}, a_{4}}$. By symmetry, we may assume $\varphi(u_{6}) \in \set{a_{0}, a_{3}}$.
	
	First suppose that $\varphi(u_{6}) = a_{0}$. Then
	\begin{equation*}
		\varphi(\set{v_{0}, v_{1}, \dotsc, v_{6}}) \subset \varphi(\Gamma(u_{1}) \cup \Gamma(u_{6})) \subset \Gamma(\varphi(u_{1})) \cup \Gamma(\varphi(u_{6})) = \Gamma(a_{0}).
	\end{equation*}
	However, $v_{0}v_{1}\dotsc v_{6}$ form a 7-cycle which is 3-chromatic, while $\Gamma(a_{0})$ is a path of length 3 (which is bipartite).
	
	Now suppose that $\varphi(u_{6}) = a_{3}$. The edge $tv_{0}$ is in $\Gamma(u_{1}, u_{6})$ so $\varphi(t)\varphi(v_{0})$ must be an edge in $\Gamma(a_{0}, a_{3})$. In particular, $\set{\varphi(t), \varphi(v_{0})} = \set{a_{1}, a_{2}}$. By symmetry we may assume that $\varphi(v_{0}) = a_{1}$. Now $v_{1} \in \Gamma(v_{0}, u_{6})$, so $\varphi(v_{1}) \in \Gamma(a_{1}, a_{3})$ and so $\varphi(v_{1}) = a_{2}$. Next $v_{2} \in \Gamma(u_{1}, u_{6}, v_{1})$, so $\varphi(v_{2}) \in \Gamma(a_{0}, a_{3}, a_{2})$ and so $\varphi(v_{2}) = a_{1}$. Working in this way round the outer 7-cycle gives $\varphi(v_{3}) = a_{2}$, $\varphi(v_{4}) = a_{1}$, and $\varphi(v_{5}) = a_{2}$. Finally $v_{6} \in \Gamma(u_{1}, v_{0}, v_{5})$ and so $\varphi(v_{6}) \in \Gamma(a_{0}, a_{1}, a_{2}) = \emptyset$, which is a contradiction.
\end{Proof}

\begin{proposition}
	The graph $T_{0}$ is not homomorphic to $H_{1}^{++}$.
\end{proposition}

\begin{Proof}
	Suppose $\varphi \colon T_{0} \to H_{1}^{++}$ is a homomorphism. If $x, y$ is a dense pair (see \cref{def:dense}) of vertices in $T_{0}$, then $\Gamma(x, y)$ contains an edge, so $\Gamma(\varphi(x), \varphi(y))$ contains an edge and so either $\varphi(x) = \varphi(y)$ or $\varphi(x), \varphi(y)$ is a dense pair in $H_{1}^{++}$.
	
	For a graph $G$, let \defn{$\cD G$} be the graph with vertex set $V(G)$ and with vertices $x$ and $y$ adjacent if $x, y$ is a dense pair in $G$. The previous paragraph shows that $\varphi$ maps a connected set of vertices in $\cD T_{0}$ to a connected set in $\cD H_{1}^{++}$. The graphs $\cD T_{0}$ and $\cD H_{1}^{++}$ are displayed below (we have used the same labelling of the vertices of $H_{1}^{++}$ as in \cref{H12homC72}).
	
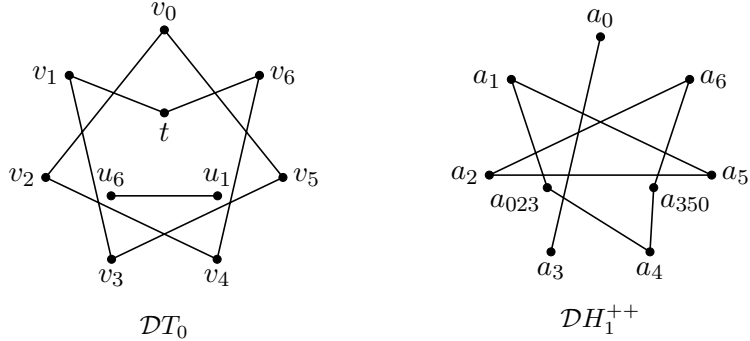
	\begin{figure}[H]
		\centering
		\begin{subfigure}{.33\textwidth}
			\centering
			\begin{tikzpicture}
				\foreach \pt in {0,1,...,6} 
				{
					\tkzDefPoint(\pt*360/7 + 90:1.6){v_\pt}
				}
				\tkzDefPoint(0,0.5){t}
				\tkzDefPoint(0.7,-0.6){u_1}
				\tkzDefPoint(-0.7,-0.6){u_6}
				\tkzDrawPolySeg(v_0,v_2,v_4,v_6,t,v_1,v_3,v_5,v_0) 
				\tkzDrawSegment(u_1,u_6)
				
				\tkzDrawPoints(v_0,v_...,v_6)
				\tkzDrawPoints(t,u_1,u_6)
				\tkzLabelPoint[above](v_0){$v_{0}$}
				\tkzLabelPoint[left](v_1){$v_{1}$}
				\tkzLabelPoint[left](v_2){$v_{2}$}
				\tkzLabelPoint[below](v_3){$v_{3}$}
				\tkzLabelPoint[below](v_4){$v_{4}$}
				\tkzLabelPoint[right](v_5){$v_{5}$}
				\tkzLabelPoint[right](v_6){$v_{6}$}
				\tkzLabelPoint[below](t){$t$}
				\tkzLabelPoint[above](u_1){$u_{1}$}
				\tkzLabelPoint[above](u_6){$u_{6}$}
			\end{tikzpicture}
			\caption*{$\cD T_{0}$}
		\end{subfigure}
		\begin{subfigure}{.33\textwidth}
			\centering
			\begin{tikzpicture}
				\foreach \pt in {0,1,...,6} 
				{
					\tkzDefPoint(\pt*360/7 + 90:1.5){v_\pt}
				}
				\tkzDefPoint(-0.7,-0.5){ul}
				\tkzDefPoint(0.7,-0.5){ur}
				
				\tkzDrawPolySeg(v_1,ul,v_4,ur,v_6,v_2,v_5,v_1) 
				\tkzDrawPolySeg(v_0,v_3)
				
				\tkzDrawPoints(v_0,v_...,v_6)
				\tkzDrawPoints(ul,ur)
				\tkzLabelPoint[above](v_0){$a_{0}$}
				\tkzLabelPoint[left](v_1){$a_{1}$}
				\tkzLabelPoint[left](v_2){$a_{2}$}
				\tkzLabelPoint[below](v_3){$a_{3}$}
				\tkzLabelPoint[below](v_4){$a_{4}$}
				\tkzLabelPoint[right](v_5){$a_{5}$}
				\tkzLabelPoint[right](v_6){$a_{6}$}
				\tkzLabelPoint[below left = -2pt](ul){$a_{023}$}
				\tkzLabelPoint[below right = -2pt](ur){$a_{350}$}
			\end{tikzpicture}
			\caption*{$\cD H_{1}^{++}$}
		\end{subfigure}
	\end{figure}
	\addtocounter{figure}{-1}
	
	Let $C$ be the 8-cycle $v_{0}v_{2}v_{4}v_{6}tv_{1}v_{3}v_{5}$ of $\cD T_{0}$ and so $\varphi(C)$ is connected in $\cD H_{1}^{++}$. If $\varphi(C)$ meets $\set{a_{0}, a_{3}}$, then $\abs{\varphi(C)} \leqslant 2$ and so $\abs{\varphi(T_{0})} \leqslant 4$ while $\chi(\varphi(T_{0})) \geqslant \chi(T_{0}) = 4$ so $\varphi(T_{0})$ is a 4-clique which is absurd as $H_{1}^{++}$ is $K_{4}$-free. Hence, $\varphi(C) \subset H_{1}^{++} - \set{a_{0}, a_{3}}$.
	
	Now both $H_{1}^{++} - a_{0}$ and $H_{1}^{++} - a_{3}$ are 3-colourable (they are both 2-degenerate) and $\chi(\varphi(T_{0})) \geqslant 4$, so $a_{0}, a_{3} \in \varphi(T_{0})$. In particular, $\varphi(\set{u_{1}, u_{6}}) = \set{a_{0}, a_{3}}$. By symmetry we may assume that $\varphi(u_{1}) = a_{0}$ and $\varphi(u_{6}) = a_{3}$.
	
	In $T_{0}$, the path $v_{2}v_{3}v_{4}v_{5}$ lies in the common neighbourhood of $u_{1}$ and $u_{6}$. While, in $H_{1}^{++}$, the common neighbourhood of $a_{0} = \varphi(u_{1})$, $a_{3} = \varphi(u_{6})$ consists of two disconnected edges $a_{2}a_{023}$ and $a_{5}a_{350}$. Thus $\set{\varphi(v_{2}), \varphi(v_{5})}$ is either $\set{a_{2}, a_{023}}$ or $\set{a_{5}, a_{350}}$.
	
	Back in $\cD T_{0}$, $v_{2}v_{0}v_{5}$ is a path, so $\varphi(v_{2})$, $\varphi(v_{5})$ are within distance two of each other in $\cD H_{1}^{++}$. But this is inconsistent with $\set{\varphi(v_{2}), \varphi(v_{5})}$ being either $\set{a_{2}, a_{023}}$ or $\set{a_{5}, a_{350}}$.
\end{Proof}
	
\end{appendix}

\end{document}